
\mag=\magstep1
\documentstyle{amsppt}
\input amsppt1
\pageheight{23.5true cm}
\pagewidth{15.5true cm}
\parindent=4mm
\baselineskip=13pt plus.1pt
\parskip=3pt plus1pt minus.5pt
\nologo
\NoRunningHeads
\NoBlackBoxes

\def  \Rh{$\rho^{\text{an}}
(X \! \! \supset \! \! E / \, Y \! \! \ni \! \! P) \, $}
\def  \NE{$\overline{NE}^{\, \text{an}}
(X \! \! \supset \! \! E / \, Y \! \! \ni \! \! P) \, $}

\topmatter

\title Extremal contractions from 4-dimensional manifolds to 3-folds
\endtitle
\author Yasuyuki Kachi \endauthor
\affil Department of Mathematical Sciences \\
              University of Tokyo \\
            Hongo, Tokyo 113, Japan \\
   E-mail : kachi\@ \! \! clover.math.s.u-tokyo.ac.jp \endaffil

\abstract
{Let $g : X \to Y$ be the contraction of an extremal ray
of a smooth projective 4-fold $X$ such that $\dim Y=3$.
Then $g$ may have a finite number of 2-dimensional fibers.
We shall classify those fibers.
Especially we shall prove that
any two points of such a fiber is
joined by a chain of rational curves of
length at most 2 with respect to $-K_X$, and that
$|-K_X|$ is $g\text{-free}$.}

\toc
\subhead
 \S 0.\, \, Introduction
 \endsubhead
\subhead
 \S 1.\, \, The universal family of extremal
 rational curves of length 1
 (In the case $l_E(R)=1$)
 \endsubhead
\subhead
 \S 2.\, \, The universal family of extremal
 rational curves of length 2
 \endsubhead
\subhead
 \S 3.\, \, Rational 2-chain connectedness of $E$
 \endsubhead
\subhead
 \S 4.\, \, Relative base-point-freeness of $|-K_X|$
 \endsubhead
\subhead
 \S 5.\, \, The case $l_{E_i}(R)=2$ for some $i$
 \endsubhead
\subhead
 \S 6.\, \, Glueing chain of rational curves
\endsubhead
\subhead
 \S 7.\, \, Proof of \Rh $\leq 2$
 \endsubhead
\subhead
 \S 8.\, \, The classification of $E$ in the case
 $l_E(R)=1$ and \Rh$ \, = \, 2$
 \endsubhead
\subhead
 \S 9.\, \, Deformation loci of extremal rational curves of length 1
 (the case \Rh$=1$)
 \endsubhead
\subhead
 \S 10.\, \, The classification of \, $E$ \, in the case \, $l_E(R)=1$
 \, and $\rho^{\text{an}} (X \! \! \supset$
 \endsubhead
\subhead
 \, \, \, $E/Y \! \! \ni \! \! P)=1$
 \endsubhead
\subhead
 \S 11.\, \, Examples
 \endsubhead
\endtoc

\endtopmatter

\document
We will work over $\Bbb C$,
the complex number field.

\newpage

\head \S 0.\ Introduction \endhead

In the classification theory of
higher dimensional algebraic varieties,
it is believed that every algebraic variety
is transformed, via successive birational
maps determined by extremal rays, either to
a minimal model or a model which admits
a structure of Mori fiber space, {\it i.e.\/}
a fiber space determined by an extremal ray.
This is so-called the Minimal Model Conjecture
(see for {\it e.g.\/} [KaMaMa]),
which was classical in dimension 2,
and was solved affirmatively also
in dimension 3 by Mori [Mo4], with [Mo2,3],
Reid [R2], Kawamata [Ka1,2,3], Shokurov [Sh2].
(See also [Utah], Koll\'ar-Mori [KoMo],
Shokurov [Sh3] and Kawamata [Ka6] for
further developments).
In the first place,
it was Mori [Mo2] who introduced
the notion of extremal rays.
He proved the existence of the
contraction morphism associated to
any extremal ray in the case of
smooth projective 3-folds, and
completely classified those structures
[loc.cit].
Kawamata [Ka1] then generalized
the existence of the contraction morphism
to the case of singular 3-folds
which was necessary for the
Minimal Model Conjecture.
After then
Kawamata [Ka2] and Shokurov [Sh2]
furthermore generalized it to
an arbitrary dimension,
even though the Minimal Model Conjecture
itself is still unsolved in dimension
greater than or equal to 4.

Thus, it is worth trying to investigate
the structures of those contractions
also for dimension greater than or equal to 4.
There are several results known so far
in this direction (Ando [A], Beltrametti [Bel],
Fujita [F], Kawamata [Ka4], and
Andreatta-Wi\'sniewski [AW2]).
Especially, Kawamata [Ka4] proved the
existence of flips from smooth 4-folds,
which, together with the Termination Theorem
[KaMaMa] after [Sh2],
should be considered as the first step
of generalizing the Minimal Model Conjecture
to dimension 4.

Our interest is the extremal contractions
from smooth projective 4-folds to 3-folds.
More generally, let $g: X \to Y$ be
the contraction of an extremal ray
of smooth projective $n\text{-folds}$ $X$
such that $\dim Y = n-1$. For $n=3$,
Mori [Mo2] proved that $g$ is a
conic bundle (see Beauville [Bea],
{\it cf.\/} Sarkisov [Sa], Ando [A]).
Moreover for an arbitrary $n \geq 4$,
if we assume that $g$ is equi-dimensional,
then Ando [A] proved that $g$ is a
conic bundle also in this case.
In $n \geq 4$, however, the situation is
more complicated, namely, it is no longer
true in general that $g$ is equi-dimensional.
Actually, even in the case
$n=4$, $g$ may admit a finite number of
2-dimensional fibers (Beltrametti [Bel]
Example 3.6, Mukai, and Reid, we shall
give essentially the same example as theirs
in Example 11.1). On the other hand, it can
be shown that $g$ is still a conic bundle
elsewhere, by modifying the argument
of [Mo2], [A] (See also Proposition 2.2
below).

The purpose of this paper is to
describe the structure of
2-dimensional fibers of $g : X \to Y$,
when $n=4$. The classification results
will be given in
Theorem 0.6, 0.7 and 0.8 below.
In the course of determining
those local structures, we necessarily
need the results on other types of
contractions from 4-folds, such as flips
[Ka4], divisorial contractions [A], [AW2],
and flops.

We put no assumptions on the singularities of $Y$,
although we assume that $X$ is smooth.
The general theories (Kawamata-Matsuda-Matsuki
[KaMaMa] and Koll\'ar [Ko1,2]) tell us that
$Y$ has only isolated $\Bbb Q\text{-factorial}$
rational singularities.

Throughout this paper,
we fix the following notation unless otherwise stated:

\definition{Notation 0.1}\ \
Let $X$ be a 4-dimensional
smooth projective variety, $R$ an extremal
ray of $X$, and $g : X \to Y$
the contraction morphism associated to $R$.
Assume that $\dim Y =3$. Let $E$ be any
2-dimensional fiber of $g$,
$E = \bigcup\limits_{i=1}^n E_i$
the irreducible decomposition of $E$,
and $\nu_i : \widetilde{E}_i \to E_i$
the normalization of $E_i$.
Let $V$ be a sufficiently small
analytic neighborhood of $P := g(E)$ in $Y$,
let $U := g^{-1}(V)$, and
$g_U := g|_U : U \to V$.
\enddefinition

To state our main result, we shall prepare
some terminologies and notations.

\definition{Definition 0.2} \ \
Let
$$
\align
l_E(R):= \text{ Min } \{ (-K_X \,. \, \, C) \quad | \quad
C \text{ is an irreducible rational curve}& \\
\text{contained in } E& \},
\endalign
$$
and call the {\it length \/} of $R$ at $E$.
It is easily seen that
$l_E(R)=1$ or $2$.

Also, for any curve (or 1-cycle) $C$
of $X$, we call $(-K_X \, . \, C)$
the {\it length \/} of $C$.
\enddefinition

\definition{Definition 0.3}\ \
Under the Notation 0.1,
$g_U|_{U-E} : U \! - \! E \to V \! - \! P$
is a proper flat morphism
whose fiber is 1-dimensional.
Furthermore, this is actually a conic bundle
by the results of Ando [A].
(See also Proposition 2.2 below.)
A 1-dimensional closed subscheme $C$ of $X$ is
called a {\it limit conic\/}
if there is a curve $A$ in $(V,P)$ passing through $P$
such that $C$ is the fiber at $P$
of the 1-parameter family $\{X_t\}_{t \in A}$
induced from $g_U|_{g^{-1}(A) - E}
: g^{-1}(A) \! - \! E \to A \! - \! P$
by taking the closure
(See also Notation 2.1.)

We note that $(-K_X \, . \, C)=2$,
since a general fiber $( \simeq \Bbb P^1)$ of $g$
has length 2.
In particular, $C$ satisfies either one of the
followings:
$$
\cases
\text{(0.3.1)} \, \, C \text{ is irreducible,
generically reduced, and } (-K_X \, . \, C)=2. \\
\text{(0.3.2)} \, \, C \text{ is generically reduced and }
C_{\text{red}} = l \cup l' \text{ for irreducible rational curves} \\
\qquad \, \, \, \, \, l, l' \subset E \text{ with }
(-K_X \, . \, l) = (-K_X \, . \, l') = 1. \\
\text{(0.3.3)} \, \, C \text{ is not generically reduced, }
(-K_X \, . \, C)=2, \text{ and }
(-K_X \, . \, C_{\text{red}}) = 1.
\endcases
$$
\enddefinition

\definition{Definition 0.4}\ \
We define
the {\it analytically-local relative cone of curves \/}
\NE
and the {\it analytically-local relative Picard number \/}
\Rh
{\it of \/} $g$ {\it along \/} $E$
as
$\overline{NE}(U/V)$ and
$\rho(U/V)$ for a sufficiently small $V \ni P$ and $U:=g^{-1}(V)$,
respectively ([R3], [Ka2,3], [N]).
\enddefinition

Now we come to stating our main results.
They consist of four parts: Theorem 0.5 (general properties),
Theorem 0.6, 0.7 and 0.8 (classifications).

\proclaim{Theorem 0.5}\ \
Let $X$ be a 4-dimensional smooth projective variety,
$R$ an extremal ray of $X$, and $g : X \to Y$
the contraction morphism associated to $R$.
Assume that $\dim Y=3$. Then the followings hold:

\flushpar
(1) \, \, \, $|-K_X|$ is $g\text{-free}$.

\flushpar
(2) \, \, \, Let $E$ be any 2-dimensional fiber
of $g$, and let $P := g(E)$. Then

\flushpar
\qquad \qquad \qquad \qquad \qquad \qquad \quad \, \, \, \,
\Rh$\leq 2$.

\flushpar
(3) \, \, \, Either one of the followings holds;

\flushpar
(CL) \, \, \, (Connected by limit conics)

\flushpar
For any two points $x$, $y$ of $E$,
there exists a limit conic in $E$
which passes through both $x$ and $y$, or

\flushpar
(MW) \, \, \, (Mukai-Wi\'sniewski type)

\flushpar
$E=E_1 \cup E_2$, $E_1 \simeq E_2 \simeq \Bbb P^2$,
$E_1 \cap E_2$ is a point, and
$N_{E_1/X} \simeq N_{E_2/X} \simeq \Cal O_{\Bbb P^2}(-1)^{\oplus 2}$.
$g_U|_{U-E} : U \! - \! E \to V \! - \! P$
is a $\Bbb P^1\text{-bundle}$.
In this case, for any points $x \in E_1$ and $y \in E_2$,
there exists a limit conic in $E$ which passes through
both $x$ and $y$.

In particular, $E$ is rationally chain connected,
in either case.
(See Definition 0.10.)
\endproclaim

According to Definition 0.2 and Theorem 0.5 (2),
there are exactly three possible types for
2-dimensional fibers $E$ of $g$:
$$
\cases
\text{Type (A)} \, : \quad l_E(R)=2, \\
\text{Type (B)} \, : \quad
l_E(R)=1 \text{ and } \text{\Rh}=2, \text{ or} \\
\text{Type (C)} \, : \quad
l_E(R)=1 \text{ and } \text{\Rh}=1.
\endcases
\tag 0.5.1
$$

Let us give a classification result for each type
of $E$ in the followings:

\proclaim{Theorem 0.6}\ \ (For the Type (A))

Assume

\flushpar
\qquad \qquad \qquad \qquad \qquad \qquad \qquad \qquad
$l_E(R)=2$.

\flushpar
Then

\flushpar
(A) \, \, \, $E$ is irreducible and is isomorphic to $\Bbb P^2$,
$N_{E/X} \simeq \Omega_{\Bbb P^2}^1(1)$,
and $Y$ is smooth at $P$.
In particular, \Rh$=1$, and
$\Cal O_{E}(-K_X) \simeq \Cal O_{\Bbb P^2}(2)$.

Furthermore, the followings hold:

\flushpar
(local elementary transformation)

Let $U$, $V$, and $g_U : U \to V$ be as in Notation 0.1.
Let $x \in E$ be an arbitrary point. Then there exists
a smooth surface $S_x \subset U$ proper over $V$
such that $g|_{S_x - x} : S_x \! - \! x \to g(S_x) \! - \! P$
is an isomorphism, and that
$S_x \cap E=\{x\}$
intersecting transversally.
Let $\varphi : \overline{U} \to U$
be the blow-up with center $S_x$.
Then $-K_{\overline{U}}$ is
$(g|_U) \circ \varphi\text{-ample}$,
and $\varphi^{-1}(E) \simeq \Sigma_1$.
Let $\varphi^+ : \overline{U} \to U^+$
be the contraction
associated to the extremal ray
of $\overline{NE}(\overline{U}/V)$
other than $\varphi$.
Then $U^+ \simeq V \times \Bbb P^1$.
Let $g^+ : U^+ \to V$ be the first projection.
Then

\flushpar
\qquad \qquad \qquad \qquad \qquad \qquad \qquad
$(g|_U) \circ \varphi = g^+ \circ \varphi^+$,

\flushpar
and $\varphi^+$ is the blow-up with center $S^+$
which is a smooth surface proper over $V$ such that
$S^+ \supset g^{+ \, -1}(P) \simeq \Bbb P^1$
with the normal bundle
$N_{g^{+ \, -1}(P)/S^+} \simeq \Cal O_{\Bbb P^1}(-1)$.
\endproclaim

\proclaim{Theorem 0.7}\ \ (For the Type (B))

Assume that

\flushpar
\qquad \qquad \qquad \qquad \qquad \, \,
$l_E(R)=1$ and \Rh$=2$.

\flushpar
Then $E$ is one of the followings:

\flushpar
(B-0) \, \, (=(MW) in Theorem 0.5) \, \, \, (Mukai-Wi\'sniewski type)

\flushpar
$E=E_1 \cup E_2$, $E_1 \simeq E_2 \simeq \Bbb P^2$,
$E_1 \cap E_2$ is a point, and
$N_{E_1/X} \simeq N_{E_2/X} \simeq \Cal O_{\Bbb P^2}(-1)^{\oplus 2}$.
In particular

\flushpar
\qquad \qquad \qquad \qquad \qquad \quad \, \,
$\Cal O_{E_i}(-K_X) \simeq \Cal O_{\Bbb P^2}(1)$ \, \, $(i=1,2)$.

\flushpar
(B-1) \, \, \, $E$ is irreducible and is isomorphic to
$\Bbb P^1 \times \Bbb P^1$,

\flushpar
\qquad \qquad \qquad \qquad \qquad \qquad \, \, \,
$\Cal O_E(-K_X) \simeq \Cal O_{\Bbb P^1 \times \Bbb P^1}(1,1)$.

\flushpar
(B-2) \, \, \, $E$ is irreducible and is isomorphic to
$\Sigma_1$,

\flushpar
\qquad \qquad \qquad \qquad \qquad \qquad \, \, \,
$\Cal O_E(-K_X) \simeq \Cal O_{\Sigma_1}(M + 2l)$

\flushpar
(Notation 0.14).

\flushpar
(B-3) \, \, \, $E=E_1 \cup E_2$, $E_1 \simeq \Bbb P^2$ and
$E_2 \simeq \Bbb P^1 \times \Bbb P^1$.
$E_1 \cap E_2$ is a line of $E_1$
and is a ruling of $E_2$.

\flushpar
\qquad \qquad \qquad \,
$\Cal O_{E_1}(-K_X) \simeq \Cal O_{\Bbb P^2}(1)$, \, and \,
$\Cal O_{E_2}(-K_X) \simeq \Cal O_{\Bbb P^1 \times \Bbb P^1}(1,1)$.

\flushpar
(B-4) \, \, \, $E=E_1 \cup E_2$, $E_1 \simeq \Bbb P^2$ and
$E_2 \simeq \Sigma_1$.
$E_1 \cap E_2$ is a line of $E_1$
and is the negative section of $E_2$.

\flushpar
\qquad \qquad \qquad \, \,
$\Cal O_{E_1}(-K_X) \simeq \Cal O_{\Bbb P^2}(1)$, \, and \,
$\Cal O_{E_2}(-K_X) \simeq \Cal O_{\Sigma_1}(M + 2l)$.

In (B-0), $P$ is an ordinary double singular point of $Y$,
and $g_U|_{U-E} : U \! - \! E \to V \! - \! P$
is a $\Bbb P^1\text{-bundle}$.

In (B-1) $\sim$ (B-4), $P$ is a smooth point of $Y$,
and $g_U|_{U-E} : U \! - \! E \to V \! - \! P$ is
a conic bundle with an irreducible
discriminant divisor.
\endproclaim

\proclaim{Theorem 0.8}\ \ (For the Type (C))

Assume that

\flushpar
\qquad \qquad \qquad \qquad \qquad \, \,
$l_E(R)=1$ and \Rh$=1$.

\flushpar
Then $E$ is one of the followings:

\flushpar
(C-1) \, \, \, $E$ is irreducible,
$E \simeq \Bbb P^2$, and
$\Cal O_E(-K_X) \simeq \Cal O_{\Bbb P^2}(1)$.

\flushpar
(C-2) \, \, \, $E=E_1 \cup E_2$,
$E_1 \simeq E_2 \simeq \Bbb P^2$,
and $E_1 \cap E_2$ is a line of both $E_i$.

\flushpar
\qquad \qquad \qquad \qquad \qquad \quad \, \,
$\Cal O_{E_i}(-K_X) \simeq \Cal O_{\Bbb P^2}(1) \, \, \, \, (i=1,2)$.

\flushpar
(C-3) \, \, \, $E$ is irreducible.
Let $\nu : \widetilde{E} \to E$ be the normalization.
Then $\widetilde{E} \simeq S_m$ and

\flushpar
\qquad \qquad \qquad \qquad \qquad \qquad \quad
$\nu^* \Cal O_E(-K_X) \simeq \Cal O_{S_m}(1)$

\flushpar
(Notation 0.14).

\flushpar
(C-4) \, \, \, $E = E_1 \cup E_2$.
Let $\nu_i : \widetilde{E}_i \to E_i$ be
the normalization of $E_i$ $(i=1,2)$. Then
$\widetilde{E}_i \simeq S_{m_i}$ for some
$m_i \geq 2$, and

\flushpar
\qquad \qquad \qquad \qquad \qquad \qquad \, \, \,
$\nu_i^* \Cal O_{E_i}(-K_X) \simeq \Cal O_{S_{m_i}}(1)$.

\flushpar
Moreover, let $v_i$
be the vertex of $\widetilde{E}_i$, then
$\nu_1(v_1) = \nu_2(v_2) =: Q$ in $E$.
$E_1 \cap E_2$ is a ruling of both $E_i$.

In each case (C-1) $\sim$ (C-4),
$g_U|_{U - E} : U \! -  \! E \to V \! - \! P$
is a conic bundle with a non-empty
discriminant divisor.
\endproclaim

We do not know at present any example for
(C-3) with $m \geq 4$, or (C-4) above.
Note that if $E$ is of type (C-3) with $m \geq 4$,
then $E$ is necessarily non-normal, and similarly
for $E_i\text{'s}$ of type (C-4).

By the way,
as for the property of fibers of extremal contractions in
a more general setting,
Kawamata showed the following as an application of
his adjunction ([Ka5] Lemma) and
Miyaoka-Mori's Theorem ([MiMo]) :

\proclaim {Remark 0.9}\ (Kawamata [Ka5], Theorem 2)\ \
Let $X$ be a projective variety of an
arbitrary dimension, and
$g : X \to Y$ the contraction of an extremal ray of $X$.
Then any non-trivial fiber of $g$
is covered by rational curves.
\endproclaim

In the context of this direction,
we obtain a stronger result (Corollary 0.11 below)
when we concentrate on $\dim X \leq 4$,
by virtue of
[Mo2], [A], [Bel], [F], [Ka4], [KoMiMo1], [C], [AW2],
together with our Theorem 0.5 (3).
To state the result, we recall
the notion of rational chain-connectedness
and rational connectedness
(Koll\'ar-Miyaoka-Mori [KoMiMo1,2,3], Campana [C]):

\definition{Definition 0.10}\ \
A scheme $S$ proper over $\Bbb C$ is
{\it rationally chain connected\/}
if for arbitrary two closed points $P$,$Q$ of $S$, there exists
finitely many rational curves $C_1, \dots, C_r$ on $S$
such that
$\bigcup\limits_{i=1}^r C_i$ is connected, $P \in C_1$,
and $Q \in C_r$.

$S$ is {\it rationally connected\/} if for
two general closed points $P,Q$ of $S$, there exists
an irreducible rational curve on $S$ which passes through
both $P$ and $Q$.
\enddefinition

\proclaim{Corollary 0.11}\ \
Let $X$ be a smooth projective variety with $\dim X \leq 4$,
and $g:X \to Y$ the contraction of an extremal ray of $X$.
Then any non-trivial fiber of $g$ is
rationally chain connected.
\endproclaim

\proclaim{Question 0.12}\ \
Does the same hold when $\dim X$ is arbitrary?
\endproclaim

\proclaim{Remark 0.13}\ \
Under the same assumption as in Corollary 0.11,
it is not true in general that
any irreducible component of
any fiber of $g$ is rationally connected, even
in the case $\dim X = 4$
(Example 11.9, inspired by Hidaka-Oguiso).
\endproclaim

This paper is organized as follows:

First in \S 1, we recall the construction of
the universal family of extremal rational
curves of length 1, following Ionescu [Io],
which is a generalization of Mori [Mo1].
We shall prove Proposition 1.4, which roughly gives
the structure of those $\widetilde{E}_i\text{'s}$
which are covered by rational curves
of length 1.
The idea is mostly indebted to Kawamata [Ka4].

In \S 2, we recall the construction of the
universal family containing a general fiber
of $g$, due to Mori [Mo2].
Especially in Proposition 2.2,
we prove that $g$ is a conic bundle
outside any 2-dimensional fiber of $g$,
by the argument of Ando [A] after [Mo2].
As a consequence of this construction,
we obtain several informations about
$E$. More precisely, we shall prove that
$E$ contains no 1-dimensional
irreducible component
(Proposition 2.4), and that each $E_i$
is a rational surface
(Poroposition 2.6).
These, together with Proposition 1.4,
give us a coarse classification of
$\widetilde{E}_i\text{'s}$
(Theorem 2.8).

In \S 3, we shall prove Theorem 3.1,
which asserts that either
(CL) Any two points of $E$ are
joined by a limit conic, or
(MW) $E$ is of Mukai-Wi\'sniewski type,
as described in Theorem 0.5 (3).
We are inspired this by
Campana [C], Koll\'ar-Miyaoka-Mori
[KoMiMo1]. The proof will be
divided into three  parts;
$$
\cases
\text{\Rh}=1 \, \, \Longrightarrow \, (CL), \\
\text{\Rh}\geq 2, \,
g_U|_{U - E} : U \! - \! E \to V \! - \! P
\text{ is a } \Bbb P^1\text{-bundle}
\, \, \Longrightarrow \, (MW), \text{ and} \\
\text{\Rh}\geq 2, \,
g_U|_{U - E} : U \! - \! E \to V \! - \! P
\text{ is not a } \Bbb P^1\text{-bundle}
\, \, \Longrightarrow \, (CL).
\endcases
$$

\flushpar
The first part is easier, and the proof
will be given in Corollary 3.4.
Next we prove the second part
in Theorem 3.5.
In this case \Rh$=2$, and both of extremal rays
define flipping contractions (Kawamata [Ka4]).
Thus by virtue of Kawamata's characterization [loc.cit]
we determine the structure
of $E$ as described in (MW).
For the last part, we will give
in Theorem 3.7 the structure
theorem of the cone \NE, which gives
a one-to-one correspondence
between the set of extremal rays of \NE
and the set of rational curves of length 1
in $E$ modulo numerical equivalence in $U$.
Then the last part
follows directly from Theorem 3.7
(Corollary 3.8).

In \S 4, we shall prove the relative
freeness of $|-K_X|$ (Theorem 4.1).
This proof heavily relies on Kawamata's
base-point-free technique ({\it e.g.\/} [Ka4]).
In the course of the proof, we need a sort of
property which the direct image sheaf
$g_*\Cal O_U(-K_U)$ should satisfy (Lemma 4.3).

So far we prepared most of the tools
for classifying $E\text{'s}$, and
we then go to the classification from \S 5 on.

First in \S 5, we deal with the case
$l_{E_i}(R)=2$ for some $i$,
especially with the Type (A) in (0.5.1).
For this case the answer is simple:

$E$ is irreducible and is isomorphic to
$\Bbb P^2$.
The normal bundle $N_{E/X}$ is
isomorphic to $\Omega_{\Bbb P^2}^1(1)$.
$P$ is a smooth point of $Y$.
(Theorem 5.1.)

First by the results of the previous
sections, $E$ is irreducible and
$\widetilde{E} \simeq \Bbb P^2$ (Lemma 5.4).
To prove the smoothness of $E$,
we consider a blow-up $\overline{U} \supset \overline{E}$
of $U \supset E$. Then
$-K_{\overline{U}}$ is still
ample over $V$ (Lemma 5.5), and
another contraction $\overline{U} \to U^+$
is found to be a divisorial contraction
with $U^+ \simeq V \times \Bbb P^1$
(Ando [A]).
Consequently $\overline{E}$ is smooth [loc.cit],
and so is $E$
(Proposition 5.6).
The latter half of this section is
devoted to determining $N_{E/X}$.
To do this we consider
$\big| -\dfrac{1}{2}K_U \big|$, and
cut out $N_{E/X}$ by its member.
Then we find from Mori [Mo2] and
Andreatta-Wi\'sniewski [AW1]
that $N_{E/X}$ is a uniform bundle,
and then by Van de Ven's theorem [V],
we conclude
$N_{E/X} \simeq \Omega_{\Bbb P^2}^1(1)$
(Theorem 5.8).

\S 6 is a preparation of later sections.
In this section, we recall Koll\'ar-Miyaoka-Mori's
technique of glueing chain of
rational curves [KoMiMo2,3].
(See Theorem 6.2).
It was
Andreatta-Wi\'sniewski [AW2]
who first applied this method to classify
4-dimensional divisorial contractions,
and they are very much useful also for
our case.
We are much indebted to their idea.

In \S 7, we prove \Rh$\leq 2$ (Theorem 7.1).
To do this, we consider a pair,
say $\{l_1, l_2\}$, of rational curves
of length 1 in $E$ which intersects
with each other, and belong to
distinct extremal rays in \NE (Lemma 7.3).
Then it can be shown in Lemma 7.7 that
$l_1$, $l_2$ are two intersecting
rulings of the same $E_i$ which is
isomorphic to $\Bbb P^1 \times \Bbb P^1$,
unless $l_1 + l_2$ deforms outside $E$.
This part essentially needs
Koll\'ar-Miyaoka-Mori's glueing technique
given in the previous section.
So if we assume \Rh$\geq 3$, then
$E$ admits another irreducible component
which is also isomorphic to
$\Bbb P^1 \times \Bbb P^1$. This
contradicts the condition (CL)
of Theorem 3.1.

In \S 8, we classify
$E$ with $l_E(R)=1$ and \Rh$=2$,
namely, of Type (B) in (0.5.1).
The result is given in Theorem 8.1.
In this case, our $g$ is
factored locally by a birational
contraction. Moreover,
if $g_U|_{U-E} : U \! - \! E \to V \! - \! P$
is a $\Bbb P^1\text{-bundle}$,
then $E$ is of Mukai-Wi\'sniewski type,
as in Theorem 3.5. If not,
then all the local contractions are
of divisorial type.
So in this case we can apply
the classifications of 4-dimensional divisorial
contractions
due to Andreatta-Wi\'sniewski[AW2]
(non-equi-dimensional case)
and Ando[A] (equi-dimensional case)
to get our classification.

The rest sections are devoted to
the study of $E$ of $l_E(R)=1$ and
\Rh$=1$, namely, of Type (C) in (0.5.1).
This case is a little more complicated.

Before going to the classification,
we first prove in \S 9
that the deformation locus of any
rational curve of length 1
is purely codimension 1 (Theorem 9.1)
also for the case \Rh$=1$,
as well as Theorem 3.7
for the case \Rh$\geq 2$ and
not of Mukai-Wi\'sniewski type.
The proof is,
however, quite different from
that of Theorem 3.7.
To be more precise, we follow the argument
of Kawamata [Ka4], as follows:

We have to exclude the case
that $E$ is a union of $\Bbb P^2\text{'s}$,
and $g_U|_{U - E} : U \! - \! E \to V \! - \! P$
is a $\Bbb P^1\text{-bundle}$
(Assumption 9.2).
If such case happened, then,
by the argument of ([Ka4] (2.4)),
together with our Lemma 4.3,
the normal bundle of those $\Bbb P^2\text{'s}$
in $X$ would all be isomorphic to
$\Cal O_{\Bbb P^2}(-1)^{\oplus 2}$
(Proposition 9.6). From this
we easily obtain a contradiction
to the assumption \Rh$=1$.
Kawamata's method used here consists of
the theory of vector bundles
on $\Bbb P^2$ ([V], [GM], [OSS]),
and the formal function theory.

In \S 10 we classify $E\text{'s}$
of Type (C). The hardest part,
in the course of the classification,
is to exclude the possibility that
some $\widetilde{E}_i$ is a
geometrically ruled surface
(Proposition 10.3).
Then we find that $E$ is either
an irreducible $\Bbb P^2$,
a union of two $\Bbb P^2\text{'s}$,
or a union of some rational cones (Proposition 10.6).
In the last case,
it will be proved furthermore that
the number of irreducible components is
at most 2, by considering
the flopping contraction from
a certain blow-up of $U$
(Proposition 10.7).
This is a relative and 4-dimensional
analogue of
(double) projections developped by
Fano-Iskovskikh[Is]-Takeuchi[T].
Thus we get the classification as given
in Theorem 10.1.

Finally in \S 11, we give some examples.

\definition{Notation 0.14}\ \
Throughout this paper, we shall use the symbols
$\Sigma_m$, $S_m$ and
$\Cal O_{S_m}(1)$ in the following sense:

We denote by $\Sigma_m :=
\Bbb P (\Cal O_{\Bbb P^1} \oplus \Cal O_{\Bbb P^1}(m)) \to \Bbb P^1$
the $m\text{-th}$ Hirzebruch surface.
Its any fiber $l$ is called
{\it a ruling \/} of $\Sigma_m$.
If $m \geq 1$, then there is a unique
section $M$ with $(M^2)<0$, called
{\it the minimal section \/}, or
{\it the negative section \/}.

We denote by $S_m$
the image of $\Phi_{|M + ml|} : \Sigma_m \to \Bbb P^{m+1}$.
{\it A ruling \/} on $S_m$ is the image of
a ruling of $\Sigma_m$ by $\Phi$.

Furthermore, the symbol $\Cal O_{S_m}(1)$
always means the line bundle
$\Cal O_{\Bbb P^{m+1}}(1) \otimes \Cal O_{S_m}$
on $S_m$.
\enddefinition

\definition {Acknowledgement}

The author would like to express
his sincere gratitude for Professor
Y.Kawamata, who encouraged him for a long time
and gave him much instruction.

He is grateful to  Professor S.Mori who
suggested this problem to him
with various meaningful advices,
especially told him Beltrametti, Mukai, Reid's
example (Example 11.1), which was
the starting point of this work.

He is also grateful to Professor S.Mukai
who not only invited him to Nagoya University
in August, 1993 with giving him a special
hospitality and entertainment,
but also pointed out to him a mistake
of the first version of this paper,
including the lack of Mukai-Wi\'sniewski type,
and told him a systematic way of constructing
many examples.

He would like to express
his hearty gratitude for Professors M.Andreatta and
J.Wi\'sniewski who kindly communicated to him
their work of 4-dimensional divisorial contractions.
Especially, Professor J.Wi\'sniewski
first suggested to him a possibility of
the existence of Mukai-Wi\'sniewski type,
with much encouragement.

He is grateful to Professor N.I.Shepherd-Barron
who first taught him an example
which is neither $\Bbb P^2$ nor a quadric surface
(Example 11.6),
and this suggestion was very much helpful to
completing the final step of this paper.

He is grateful to Professors F.Hidaka and K.Oguiso
who helped him of constructing Example 11.9,
which gives a counter-example for
irreducible special fibers
of extremal contractions to be
rationally connected.

He is grateful to Professor E.Sato
who kindly invited him to
Kyushu University in March, 1995 and gave him
an opportunity of talking on this work.

He also would like to thank Professors
Y.Miyaoka, T.Fujita, T.Ando,
J.Koll\'ar, A.Corti, K.Matsuki, V.Alexeev
and M.Reid for their lots of helpful comments
and encouragements.
\enddefinition

\newpage

\head
\S 1.\  The universal family of extremal rational curves of length 1
(In the case $l_E(R)=1$)
\endhead

\definition {Assumption 1.0}\ \
Throughout this section, we assume

\flushpar
(1.0.1) \qquad \qquad \qquad \qquad \qquad \qquad \quad
$l_E(R)=1$.
\enddefinition

In this section, we recall the construction
of the universal family
of extremal rational curves of length 1
following Ionescu [Io] after Mori [Mo1],
and to prove Proposition 1.4 below, which is
simple but contains the key idea of this paper.
We are inspired these arguments by Kawamata [Ka4].

\definition {Construction 1.1}\ \ (Ionescu [Io](0.4))

Let $l_0$ be an irreducible rational
curve of length 1 which is contained in $E$.
Let
$$
\alpha : \Bbb P^1 \to l_0 \subset X
$$
be the composition of the normalization
of $l_0$ and the closed immersion.
We consider the deformation of $\alpha$.

\flushpar
(1.1.0) \, \, \, Let $H$ be any irreducible component
containing the point $[\alpha]$ of the
Hilbert scheme $\text{Hom}(\Bbb P^1 \, , \, X)$,

\flushpar
and $\, \Cal U_H \to H$
the corresponding universal family.
There is a natural inclusion
$\, \Cal U_H \subset \Bbb P^1 \times X \times H$.
Let $\Bbb P^1 \times X \times H$ $\to X \times H$
be the projection, and
$\text{Im}(\, \Cal U_H)$ $\subset X \times H$
the image of $\, \Cal U_H$ in $X \times H$.
Consider the projection
$\rho : \text{Im}(\, \Cal U_H) \to H$.
We take the maximal Zariski open set $H_0$
of $H$ over which $\rho$ is flat, and let
$$
\eta : H_0 \to \text{Hilb}_X
$$
be the induced morphism.
Let $T_H \subset \text{Hilb}_X$ be
the closure
of the image of $\eta$,
and
$$
p_H : S_H \subset X \times T_H \to T_H
$$
the induced universal family over $T_H$.
Note that $T_H$ is irreducible by construction.
Let $q_H : S_H \to X$ be the natural projection. Then
by construction,

\flushpar
(1.1.1) \qquad \qquad \qquad \qquad \qquad \,
$\bigl(-K_X \, . \, \, q_H(p_H^{-1}(t)) \bigr)=1$

\flushpar
for all $t \in T_H$.

Let
$$
\overline{\eta} :
H \dashrightarrow T_H \subset \text{Hilb}_X
$$
be the rational map induced from $\eta$.
Then

\flushpar
(1.1.2) \, \, \,
The original $l_0$ is parametrized by points
of $T_H$ which correspond to $[\alpha] \in H$
by this $\overline{\eta}$.
Denote such set of points by $T_H(l_0)$.
\enddefinition

\definition{Definition 1.2}\ \
In the above, let
$V$ be any analytic open subset
of $Y$ passing through $P:=g(E)$, and
let $U:= g^{-1}(V)$.
Obviously
$p_H (q_H^{-1}(U)) \supset T_H(l_0)$ (1.1.2).
Let $T_H(U \supset l_0)$ be the union of all
connected components of $p_H (q_H^{-1}(U))$
which contain at least one point of $T_H(l_0)$.
Then we define the
{\it whole deformation locus \/}
$L_U(l_0)$ {\it of \/} $l_0$ {\it in \/} $U$ as

\flushpar
(1.2.1) \qquad \qquad \qquad \quad \, \,
$L_U(l_0) := \bigcup\limits_H
q_H \bigl( \, p_H^{-1} (T_H(U \supset l_0)) \, \bigr)$,

\flushpar
where the union is taken over all irreducible components
$H$ of $\text{Hom} (\Bbb P^1 \, , \, X)$ containing $[\alpha]$,
as in (1.1.0). Remark that

\flushpar
(1.2.2) \qquad \qquad \, \, \, \,
$L_X(l_0) = \bigcup\limits_H q_H(S_H)$, \, and \,
$L_U(l_0) \subset L_X(l_0) \cap U$.
\enddefinition

The following inequality is due to Ionescu, which is essentially
based on Mori's estimation of the dimension of Hilbert schemes
([Mo1], Proposition 3).

\proclaim {Lemma 1.3}\ \ (Ionescu [Io] (0.4) (2),(3) )

\flushpar
If $H$ is chosen to be of maximal dimension among those in (1.1.0),
then

\flushpar
\qquad \qquad \qquad \qquad \qquad \qquad \qquad \qquad
$\dim T_H \geq 2$. \qquad \qed
\endproclaim

The following is inspired by Kawamata [Ka4]. Actually, it is a
slight generalization of ([Ka4], (2.2)).

\proclaim {Proposition 1.4}\ \
Let $E_1$ be a 2-dimensional
irreducible component of $E$,
let $\nu : \widetilde{E}_1 \to E_1$ be
the normalization of $E_1$, and
$\mu : E_1' \to \widetilde{E}_1$ the
minimal resolution of $\widetilde{E}_1$.
Assume that $E_1$ is covered by
rational curves of length 1 in $E_1$.
Then $E_1'$ is isomorphic either to $\Bbb P^2$
or a geometrically ruled surface.
\endproclaim

\demo{Proof}\ \
Let $p_H : S_H \to T_H$, and
$q_H : S_H \to X$ be as in Construction 1.1.
Then by assumption, for a suitable $H$,
there is a 1-dimensional irreducible closed subset
$T_1$ of $T_H$ such that
$$
q_H(p_H^{-1}(T_1))=E_1.
$$
Let us denote for simplicity
$T_H = T$, $S_H = S$, $p_H = p$, and $q_H = q$.
Let $S_1:= p^{-1}(T_1) \subset S$:

\flushpar
(1.4.0) \qquad \qquad \qquad \qquad \qquad \qquad \, \, \,
$q(S_1)=E_1$.

\flushpar
Moreover, let $\widetilde{S}_1$, $\widetilde{T}_1$
be the normalization of $S_1$, $T_1$,
respectively, and
$\gamma : \widetilde{S}_1 \to \widetilde{T}_1$
the induced morphism
from $p|_{S_1} : S_1 \to T_1$. (See the diagram (1.4.3) below.)
Then we claim that

\flushpar
(1.4.1)\, \, \, $\gamma$ is a $\Bbb P^1\text{-bundle}$.
In particular $\widetilde{S}_1$ is smooth \, \, (cf. [Io] or [Ka4]).

In fact,
a general fiber
of $\gamma$ is isomorphic to $\Bbb P^1$,
since $F$ has at most isolated singularities.
On the other hand,
any fiber of $\gamma$
is irreducible and generically reduced, since so is $p$.
Furthermore it has no embedded points,
since $\widetilde{S}_1$ is Cohen-Macaulay
and $\widetilde{T}_1$ is smooth.
Hence any fiber of $\gamma$ is isomorphic to
$\Bbb P^1$, namely,
$\gamma$ is a $\Bbb P^1\text{-bundle}$, and
we have (1.4.1).

Since $q|_{S_1} : S_1 \to E_1$ is surjective (1.4.0),
we have the induced surjective morphism
$\beta : \widetilde{S}_1 \to \widetilde{E}_1$.
Then there is a suitable birational
morphism $\delta : S_1' \to \widetilde{S}_1$
such that $S_1'$ dominates
$E_1'$ :
$$
\beta' : S_1' \to E_1'.
$$
Since a general fiber
$l \simeq \Bbb P^1$ of
$\gamma \circ \delta : S_1' \to \widetilde{T}_1$
satisfies
$(-K_{S_1'} \, . \, l)=2$ (1.4.1),
and is sent birationally to $l':=\beta'(l)$,
we have

\flushpar
(1.4.2) \qquad \qquad \qquad \qquad \qquad \quad \, \, \, \,
$(-K_{E _1'} \, . \, l') \geq 2$.

$$\CD
l' @. \subset @. E_1' @>\mu>> \widetilde{E}_1 @>\nu>> E_1
@. \subset @. X \\
@AAA @.@AA{\beta'}A @AA{\beta}A @A{q|_{S_1}}AA @. @AA{q}A \\
l @. \subset @. S_1' @>>\delta> \widetilde{S}_1 @>>> S_1 @. \subset @. S \\
@. @. @. @VV{\gamma}V @VVV @. @VV{p}V \\
@. @. @. \widetilde{T}_1 @>>> T_1 @. \subset @. T
\endCD
\tag 1.4.3
$$

Now assume that $E_1'$ is isomorphic
neither to $\Bbb P^2$ nor a geometrically ruled surface,
to get a contradiction.
Then each extremal ray of
$\overline{NE}(E_1')$
is spanned by a $(-1)\text{-curve}$ ([Mo2] Chapter 2).
Thus by (1.4.2), there are two $(-1)\text{-curves}$
$C_1$, $C_2$ (maybe equal) and
a pseudo-effective 1-cycle $C'$
such that

\flushpar
(1.4.4) \qquad \qquad \qquad \qquad \qquad \quad
$l' \equiv C_1 + C_2 + C'$.

\flushpar
Moreover, $\mu(C_i) \not\equiv 0$ $(i=1,2)$,
since $\mu$ is a minimal resolution.
Thus by the ampleness of $-\nu^* K_X$ on $\widetilde{E}_1$,
we have
$$
\cases
(-\mu^* \nu^* K_X \, . \, C_i)
= (-\nu^* K_X \, . \, \mu(C_i))
\geq 1 \quad (i=1, 2), \, \, \, \text{ and} \\
(-\mu^* \nu^* K_X \, . \, C')
= (-\nu^* K_X \, . \, \mu(C')) \geq 0.
\endcases
\tag 1.4.5
$$
\flushpar
By intersecting $-\mu^* \nu^* K_X$
with the numerical equivalence (1.4.4),

\flushpar
(1.4.6) \qquad \qquad \qquad \qquad \qquad \, \, \,
$(- \mu^* \nu^* K_X \, . \, l') \geq 2$.

\flushpar
On the other hand, by construction (1.4.3),
$\nu \circ \mu(l')$ is parametrized
by $T_1$;
$$
\nu \circ \mu(l') = q(p^{-1}(t)) \quad (\exists t \in T_1)
$$
Thus

\flushpar
(1.4.7) \qquad \qquad \qquad \, \, \, \,
$(-\mu^* \nu^* K_X \, . \, l')
= (-K_X \, . \, \nu \circ \mu (l')) =1$

\flushpar
(1.1.1). This contradicts (1.4.6),
and hence the proposition.
\quad \qed
\enddemo

\definition {Definition 1.5}\ \
Let $p_H : S_H \to T_H$ and $q_H : S_H \to X$
be as in Construction 1.1.
Then for any $x \in X$,
we define
$$
T_{H,x} := p_H(q_H^{-1}(x)).
$$
Note that \, \,
$p_H : q_H^{-1}(x) \overset \sim\to\longrightarrow T_{H,x}$.
\enddefinition

This gives a subfamily of $p_H$
parametrizing extremal rational curves
of length 1 and passing through the point $x$ .

\proclaim {Lemma 1.6}
\nolinebreak
\ \ (Wi\'sniewski [Wi\'s], Claim in (1.1))

$$
\dim T_{H,x} \leq 1
$$
for any $ x \in E $. \ \qed
\endproclaim

The following proposition is for the case $\dim T_{H,x} =1$.

\proclaim{Proposition 1.7} \ \
Assume $l_{E_1}(R)=1$
and $\dim T_{H,x}=1$
for some $H$ and for some
$x \in \text{Reg } E_1 - \bigcup\limits_{j \not= 1} E_j$.
(See Definition 1.5.)
Then $\widetilde{E}_1 \simeq \Bbb P^2$
and $\nu^* \Cal O_{E_1}(-K_X) \simeq \Cal O_{\Bbb P^2} (1)$.
\endproclaim

\demo{Proof} \ \
Choose a 1-dimensional
irreducible component $T_1$ of $T_{H,x}$,
and let $\widetilde{T}_1$ be its normalization.
Then, as in the proof of Proposition 1.4,
we have a $\Bbb P^1$-bundle
$\gamma : \widetilde{S}_1 \to \widetilde{T}_1$
over $\widetilde{T}_1$ and a surjective morphism
$\beta : \widetilde{S}_1 \to \widetilde{E}_1$.
(See also the diagram (1.4.3).)
Let $\widetilde{x} \in \widetilde{E}_1$
be such that $\nu (\widetilde{x}) = x$.
Then $\beta^{-1}(\widetilde{x})$ gives a section of $\gamma$.
Since $\widetilde{x}$ is a smooth point of $\widetilde{E}_1$,
we necessarily have
$\widetilde{E}_1 \simeq \Bbb P^2$.
Moreover, if we write
$$
\nu^* \Cal {O}_{E_1}(-K_X) \simeq \Cal {O}_{\Bbb{P}^2}(d) \qquad
(d \in \Bbb{Z}_{>0}),
$$
then $d=1$ by $l_{E_1}(R)=1$. \quad \qed
\enddemo

\head \S 2.\ The universal family of extremal rational curves of
length 2 \endhead

In this section,
we shall first recall the construction
of the universal family of
extremal rational curves of length 2
on $X$ containing a general
fiber of $g$ ([Mo2]). We shall prove that
the total space of such family
is a birational modification of $X$ (Proposition 2.2).
This is a 4-dimensional analogue of [Mo2], (3.24).

As an application, it will be proved
that a 2-dimensional fiber $E$ of $g$ never
admits a 1-dimensional irreducible component
(Proposition 2.4),
and that each irreducible component
of $E$ is a rational surface
(Proposition 2.6).
Moreover, in Theorem 2.8,
we coarsely classify
the normalization of
each irreducible component of $E$.

\definition{Notation 2.1}\ \
Let $f \simeq \Bbb P^1$ be a general fiber of $g$,
let $\text{Hilb}_{X,[f]}$
be the unique irreducible component
of the Hilbert scheme $\text{Hilb}_X$ of $X$
containing the point $[f]$, and let
$W := \left( \text{Hilb}_{X,[f]} \right)_{\text{red}}$.
Let $\phi : Z \to W $ be the induced
universal family over $W$, and
$\pi : Z \to X$ the natural projection.
Note that $\pi$ is surjective,
since $\phi$ parametrizes
a general fiber of $g$.

Let $B$ be the whole set of points of $Y$
which give 2-dimensional fibers
of $g$, and let
$Y^{\circ}:=Y-B$, $X^{\circ}:=g^{-1}(Y^{\circ})$.
Let $E = X_b$ be any 2-dimensional fiber. Then
a {\it limit conic \/} in $E$ as in Definition 0.3
is nothing but the curve parametrized by a point
of $\phi(\pi^{-1}(E))$.
\enddefinition

\proclaim{Proposition 2.2} \ \
Under Notation 2.1 above,

(1)\, \, \, $g$ is a conic bundle on $Y^{\circ}$.
In particular, $Y^{\circ}$ is smooth,
and $g |_{X^{\circ}} : X^{\circ} \to Y^{\circ}$ is flat.

(2)\, \, \, $\dim Z=4$ and
 $\pi : Z \to X$ is a birational morphism onto $X$,
which induces an isomorphism
$\pi^{-1}(X^{\circ}) \overset \sim\to\longrightarrow X^{\circ}$.

(3)\, \, \, $\phi |_{\pi ^{-1} (X^{\circ})}=g|_{X^{\circ}}$
under the identification
$\pi |_{\pi ^{-1}(X^{\circ})}
: \pi^{-1}(X^{\circ}) \overset \sim\to\longrightarrow X^{\circ}$.
\endproclaim

$$
\CD
Z @. \overset\pi\to\longrightarrow @. X @. \supset @. X^{\circ} \\
@V{\phi}VV @. @VV{g}V @. @VV{g|_{X^{\circ}}}V \\
W @. \dashrightarrow @. Y @. \supset @. Y^{\circ}
\endCD
$$

\demo{Proof} \ \
By the similar argument to Ando ([A],Theorem 3.1),
which is a generalization of
Mori ([Mo2], (3.25)),
$g$ is a conic bundle on $Y^{\circ}$.
Thus (1) holds.

For (2),(3), we shall follow the argument of ([Mo2], (3.24)).
Let $f$ be as in Notation 2.1.
First
by $N_{f/X} \simeq \Cal O_{\Bbb P^1}^{\oplus 3}$,
$H^1(N_{f/X})=0$.
So $\text{Hilb}_X$ is smooth at $[f]$ and
$\dim W=\dim_{[f]} \text{Hilb}_X=h^0(N_{f/X})=3$.
Since each fiber of $\phi$ is of dimension 1,
we have

\flushpar
(2.2.1) \qquad \qquad \qquad \qquad \qquad \qquad \, \, \, \, \,
$\dim Z=4$.

\flushpar
By (1), we regard
$g|_{X^{\circ}} : X^{\circ} \to Y^{\circ}$
as a flat family of closed subschemes of $X$
in a trivial way.
Then by the universal property of $\text{Hilb}_X$,
we have morphisms
$j:Y^{\circ} \to W$ and $i:X^{\circ} \to Z$
such that the following is a fibered product diagram:

\vskip 12.5mm
\centerline{$\square$}
\vskip -21mm
$$
\CD
X^{\circ} @>i>> Z \\
@V{g|_{X^{\circ}}}VV @VV{\phi}V \\
Y^{\circ} @>>j> W
\endCD
\tag 2.2.2
$$

\flushpar
Then the composite $\pi \circ i : X^{\circ} \to X$
obviously coincides with
the natural inclusion $X^{\circ} \subset X$.
Hence by (2.2.1),
$i$ is an open immersion
and $\pi$ is a birational morphism,
and we get (2).
Since $\phi$ is flat,
$j$ is also an open immersion
(the diagram (2.2.2)),
and hence (3).
\quad \qed
\enddemo

\proclaim{Remark 2.3} \ \

\flushpar
(1) \, \, \, Let $C_w := \pi (\phi^{-1}(w))$
for each $w \in W$. Then $[C_w] \in R$ and
$(-K_X \, . \, C_w)=2$.

\flushpar
(2) \, \, \, $\pi$ has connected fibers, and
$\text{Exc }\pi$ is purely codimension 1 in $Z$,
by a form of Zariski Main Theorem.

\flushpar
(3) \, \, \, $F := \phi (\pi^{-1}(E))$ is connected.
In particular, any two limit conics in $E$ deform
inside $E$ to each other.
\quad \qed
\endproclaim

\proclaim{Proposition 2.4} \ \
Under the Notation 2.1, let $E$ be a 2-dimensional
fiber of $g$. Then $E$ is purely 2-dimensional.
\endproclaim

\demo{Proof} \ \
Assume that $E$ has
a 1-dimensional irreducible component, to
get a contradiction.
Then by the connectedness of $E$,
$E$ has irreducible components
$l_0$ and $E_1$ such that

\flushpar
(2.4.1) \qquad \qquad \quad \, \, \, \,
$\dim l_0 =1$, \, \, \, $\dim E_1 =2$, \, \, \, and \, \, \,
$l_0 \cap E_1 \not= \emptyset$.

\flushpar
Let $F := \phi(\pi^{-1}(E))$.
First by Proposition 2.2,
$l_0$ is contained in a limit conic
(Definition 0.3 or Notation 2.1):

\flushpar
(2.4.2) \qquad \qquad \qquad \qquad \quad \,
$l_0 \subset \pi(\phi^{-1}(w_0))$ \, \, \,
$\bigl( \exists w_0 \in F \bigr)$.

\flushpar
By Remark 2.3 (3),

\flushpar
(2.4.3) \, \, \,
Any two limit conics deform inside $E$
to each other.

\flushpar
Since $l_0$ forms a whole irreducible component
of $E$, it follows from (2.4.2) and (2.4.3) that
any limit conic $C$ contains $l_0$,
and is of type (0.3.2) in Definition 0.3.
That is:

\flushpar
(2.4.4) \quad \, \, \, \,
$\pi(\phi^{-1}(w)) = l_0 \cup l_w$ \, \, \,
$\bigl( \exists \, l_w :$ rational curve of length 1,
$\, \not= l_0 \bigr)$.

\flushpar
In particular,

\flushpar
(2.4.5) \, \, \,
$l_0$ is a rational curve of length 1.

\flushpar
Then by Lemma 1.3, $l_0$ deforms outside $E$.
Let $V$ be a sufficiently small
analytic neighborhood of $P=g(E)$
in $Y$ and let $U := g^{-1}(V)$.
Then $g|_{U \! - \! E} : U \! - \! E \to V \! - \! P$
is a conic bundle (Proposition 2.2), and
any deformation of $l_0$ in $U$,
other than $l_0$ itself, is
an irreducible component of
a degenerate fiber of $g|_{U \! - \! E}$.
Thus by the result of
Beauville [Bea] (see also
Sarkisov [Sa] or Ando [A]),

\flushpar
(2.4.6) \, \, \,
The whole deformation locus $L := L_U(l_0)$
of $l_0$ inside $U$ (Definition 1.2)
is purely codimension 1.

\flushpar
Take any irreducible component $L_1$ of $L$.
Then $L_1$ is a prime divisor
with $L_1 \cap E = l_0$,
and in particular

\flushpar
(2.4.7) \, \, \, $L_1 \cap E_1 = l_0 \cap E_1$,
which is a non-empty finite set.

\flushpar
Since $U$ is a smooth 4-fold, and since
$\dim L_1 =3$, $\dim E_1 =2$,
(2.4.7) is impossible.
Hence the result.
\quad \qed
\enddemo

\proclaim{Lemma 2.5}\ \
Under the Notation 2.1,
$\pi(\text{Exc }\pi)$ is exactly the
union of all 2-dimensional fibers of $g$:
$$
\pi(\text{Exc }\pi) = g^{-1}(B).
$$
\endproclaim

\demo{Proof}\ \
By Proposition 2.2,
$\pi(\text{Exc }\pi) \subset g^{-1}(B)$.
Assume that $\pi(\text{Exc }\pi) \subsetneq g^{-1}(B)$,
to get a contradiction.
Then there exist a 2-dimensional fiber $E$ of $g$
and an irreducible component
$E_i$ of $E$ such that
$$
\dim (\pi(\text{Exc }\pi) \cap E_i) \leq 1.
$$
Note that

\flushpar
(2.5.0) \qquad \qquad \qquad \qquad \qquad \qquad \quad
$\dim E_i = 2$

\flushpar
by Proposition 2.4.

\flushpar
(2.5.1) \, \, \,
Let $x \in E_i - \pi(\text{Exc }\pi)$ be a
general point.

\flushpar
Then $\pi^{-1}(x)$
is a point, and

\flushpar
(2.5.2) \, \, \,
$C^x := \pi \, \phi^{-1} \, \phi \, \pi^{-1}(x)$
is the unique limit conic which passes
through $x$.

\flushpar
Let $S$ be the intersection
of two general very ample
divisors in $X$ both of which pass through $x$ such that

\flushpar
(2.5.3) \qquad \qquad \qquad \, \, \, \,
$\dim (S \cap E)=0$ \, and \, $(S - x) \cap C^x = \emptyset$.

\flushpar
Let $\overline{S}$ be the proper transform of
$S$ in $Z$, and let

\flushpar
(2.5.4) \qquad \qquad \qquad \qquad \qquad \quad \, \, \,
$D := \pi \, \phi^{-1} \, \phi (\overline{S})$.

\flushpar
Since a general fiber of
$\phi|_{\overline{S}} : \overline{S} \to \phi(\overline{S})$
is isomorphic to $\Bbb P^1$,
$\phi^{-1} \, \phi (\overline{S})$
is a prime divisor in $Z$, and hence
$D$ is a prime divisor in $X$.
Moreover, since $D$ is disjoint from a
general fiber of $g$ (2.5.4),
and since $g$ is the contraction of an extremal ray of $X$,
it follows that

\flushpar
(2.5.5) \qquad \qquad \qquad \qquad \qquad \qquad \quad \, \, \, \,
$D \underset{g}\to {\equiv} 0$.

\flushpar
Note that $D \cap E \not= \emptyset$, since $D \cap E \ni x$.
Thus if $D \not\supset E$, then
we can find an irreducible curve $C$ in $E$
such that $D \not\supset C$ and
$D \cap C \not= \emptyset$,
in particular, $(D \, . \, C)>0$,
which contradicts (2.5.5). Hence we must have

\flushpar
(2.5.6) \qquad \qquad \qquad \qquad \qquad \qquad \quad \, \, \,
$D \supset E$.

\flushpar
{}From this and the definition of $D$ (2.5.4),

\flushpar
(2.5.7) \, \, \,
$E$ is covered by some limit conics
all of which pass through
at least one point of $S \cap E$.

\flushpar
In particular,
$S \cap E \supsetneq \{x\}$ by (2.5.2).
Since $S \cap E$ is a finite set of points (2.5.3),
it follows from (2.5.7), together with (2.5.0), that
there exist
$y \in S \cap E$ and
a 1-parameter family $\{C_t\}_{t \in T}$
of limit conics such that

\flushpar
(2.5.8) \qquad \qquad \qquad \, \, \, \, \,
$E_i \subset \bigcup\limits_{t \in T} C_t$ \, \, and \, \,
$C_t \ni y$ \, \, \, ($\forall t \in T$).

\flushpar
Since $x \in E_i$ (2.5.1), this implies that

\flushpar
(2.5.9) \, \, \,
There is a limit conic which passes through
both $x$ and $y$.

\flushpar
On the other hand,
$y \not\in C^x$ (2.5.3).
This and (2.5.9) contradict (2.5.2).
Hence $\pi(\text{Exc }\pi) \supset E$.
\quad \qed
\enddemo

\proclaim{Proposition 2.6}\ \
Let $E$ be any 2-dimensional
fiber of $g$. Then

\flushpar
(1) \, \, \, Each irreducible component
of $E$ is a rational surface.

\flushpar
(2) \, \, \, There exists a finite set of points
$\{y_1, \dots , y_r\} \subset E$
such that $E$ is covered by some
limit conics all of which pass through
at least one of $y_1, \dots , y_r$.
\endproclaim

\demo{Proof}\ \
Let $\psi : Z' \to Z$ be a resolution of
the normalization of $Z$, and let
$\Phi := \pi \circ \psi : Z' \to X$.
Since $\pi(\text{Exc }\pi) = E$ (Lemma 2.5 (1)),
we have

\flushpar
(2.6.1) \qquad \qquad \qquad \qquad \qquad \qquad
$\Phi(\text{Exc }\Phi) = E$.

\flushpar
Moreover, by the same reason as in
Lemma 2.5 (2),

\flushpar
(2.6.2) \, \, \,
$\text{Exc }\Phi$ is purely codimension 1
in $Z'$.

\flushpar
In particular,
for any irreducible component $E_i$ of $E$
and for a general point $x$ of $E_i$,
$\Phi^{-1}(x)$ is purely 1-dimensional.
Hence if we take
the intersection $S$ of two general very ample
divisors in $X$, then
we may assume that

\flushpar
(2.6.3) \, \, \,
$\dim(S \cap E) = 0$, \, and \,
$\Phi^{-1}(S \cap E)$ is purely 1-dimensional.

\flushpar
Since $R^1 \Phi_* \Cal O_{Z'} =0$,

\flushpar
(2.6.4) \, \, \,
$\Phi^{-1}(S \cap E)$ is a union of finitely many
$\Bbb P^1\text{'s}$.

\flushpar
({\it cf. \/} [D].)
Moreover, by Bertini's Theorem,

\flushpar
(2.6.5) \, \, \,
$\Phi^{-1}(S)$ is an irreducible surface in $Z'$.

\flushpar
By considering the image by $\psi : Z' \to Z$,
it follows from (2.6.4) and (2.6.5) that

\flushpar
(2.6.6) \, \, \, $\pi^{-1}(S \cap E)$ is a union
of finitely many rational curves, and

\flushpar
(2.6.7) \, \, \,
$\pi^{-1}(S)$ is an irreducible surface in $Z$
containing $\pi^{-1}(S \cap E)$.

Consider $\phi|_{\phi^{-1} \, \phi \, \pi^{-1} (S \cap E)} :
\phi^{-1} \, \phi \, \pi^{-1} (S \cap E)
\to \phi \, \pi^{-1} (S \cap E)$.
Since any fiber of $\phi$ is
a union of rational curves,
it follows from (2.6.6) that

\flushpar
(2.6.8) \, \, \,
$\phi^{-1} \, \phi \, \pi^{-1} (S \cap E)$ is
a union of rational surfaces, and so is
$\pi \, \phi^{-1} \, \phi \, \pi^{-1} (S \cap E)$.

On the other hand, by (2.6.7),
$D := \pi \, \phi^{-1} \, \phi \, \pi^{-1} (S)$
is a prime divisor of $X$. Since $D$ is
disjoint from a general fiber of $g$,

\flushpar
(2.6.9) \qquad \qquad \qquad \qquad \qquad \qquad \quad \, \, \,
$D \supset E$,

\flushpar
by the same reason as in
(2.5.6). Hence

\flushpar
(2.6.10) \qquad \qquad \qquad \quad \,
$E = D \cap E = \pi \, \phi^{-1} \, \phi \, \pi^{-1} (S \cap E)$.

\flushpar
(2.6.8) and (2.6.10) prove (1).
(2.6.3) and (2.6.10) prove (2).
\quad \qed
\enddemo

\proclaim{Lemma 2.7}\ \
Let $E$ be a 2-dimensional fiber of $g$,
let $E = \bigcup E_i$ be the irreducible decomposition,
and $\nu_i : \widetilde{E}_i \to E_i$
the normalization of $E_i$.
Let $l$ be a rational curve of length 1 in $E$.

\flushpar
(1) \, \, \, If $l_{E_i}(R)=2$, then
$\bigl( \widetilde{E}_i \, , \,
\nu_i^* \Cal O_{E_i}(-K_X) \bigr)
\simeq \bigl( \Bbb P^2 \, , \, \Cal O_{\Bbb P^2}(2) \bigr)$.

\flushpar
(2) \, \, \, If $L_X(l) \subset E$, then
$L_X(l)$ is a union of some $E_i\text{'s}$
each of which satisfies
$\bigl( \widetilde{E}_i \, , \,
\nu_i^* \Cal O_{E_i}(-K_X) \bigr)
\simeq \bigl( \Bbb P^2 \, , \, \Cal O_{\Bbb P^2}(1) \bigr)$.

\flushpar
(3) \, \, \, If $L_X(l) \not\subset E$, then $L_X(l)$ is
purely codimension 1. Each irreducible component
of $L_X(l)$ contains $E$, whenever it meets $E$.

\flushpar
(4) \, \, \, If $L_X(l) \not\subset E$ for at least one rational
curve $l$ of length 1, then the same
holds for any such $l$.
\endproclaim

\demo{Proof} \ \

\flushpar
{\it (1) \/} \ \
Let $\pi : Z \to X$ and
$\phi : Z \to W$ be as in Notation 2.1.
If we choose a smooth point $x$ of $E_i$
sufficiently general so that
$x \not\in \bigcup\limits_{j \not= i} E_j$,
we have $\dim \, \pi^{-1}(x)=1$
by Lemma 2.5 (2).
Take an irreducible component
$W_1$ of $\phi(\pi^{-1}(x))$,
and let $Z_1:=\phi^{-1}(W_1)$. Moreover
let $\widetilde{Z}_1$, $\widetilde{W}_1$
be the normalization of $Z_1$, $W_1$,
respectively, and
$\gamma : \widetilde{Z}_1 \to \widetilde{W}_1$
the induced morphism. Then by $l_{E_i}(R)=2$,
each fiber of $\gamma$ is
irreducible and generically reduced.
Thus $\gamma$ is a $\Bbb P^1\text{-bundle}$, by the
same argument as in (1.4.1).
Since $x$ is chosen to be a smooth point of $E_i$,
exactly the same argument
as in the proof of Proposition 1.7 proves
$\widetilde{E}_i \simeq \Bbb P^2$.
Then
$\nu^* \Cal O_{E_i}(-K_X) \simeq \Cal O_{\Bbb P^2}(2)$
by $l_{E_i}(R)=2$.

\flushpar
{\it (2) \/} \ \
Let $H$, $T=T_H$, $S=S_H$, $p=p_H$ and $q=q_H$
be as in Construction 1.1 such that
$\dim H$ is maximal among those in (1.1.0).
$\dim T \geq 2$ by Lemma 1.3.
By assumption, $q(S) \subset L_X(l) \subset E$.
Then for any point
$x \in q(S)$,
$\dim p(q^{-1}(x)) = \dim q^{-1}(x) \geq 1$,
since $\dim S = \dim T + 1 \geq 3$
and $\dim q(S) \leq 2$.
Thus the
assumption of Proposition 1.7 is satisfied,
and it follows that
$$
\bigl( \widetilde{E}_i \, , \,
\nu_i^* \Cal O_{E_i}(-K_X) \bigr) \simeq
\bigl( \Bbb P^2 \, , \, \Cal O_{\Bbb P^2}(1) \bigr)
$$
for any $E_i$ such that $E_i \subset L_X(l)$.

\flushpar
{\it (3) \/} \ \
By assumption, there is an irreducible component
$l'$ of a 1-dimensional fiber $X_b$ of $g$
such that $l$ and $l'$ are parametrized
by $T_H$ for a same $H$ (Construction 1.1).
Let $L_1 := q_H (S_H)$ for such an $H$.
Since $g$ is a conic bundle
near $X_b$ (Proposition 2.2),
it follows that $L_1$ is
a prime divisor [Bea].
In particular, $L_1$ forms an irreducible component of
$L_X(l)$ (1.2.2).
Then by the same reason as in (2.5.6),

\flushpar
(2.7.1) \qquad \qquad \qquad \qquad \qquad \qquad \quad \, \,
$L_1 \supset E$.

Next assume that $L_X(l)$ has another
irreducible component $L_2$ with
$L_2 \cap E \not= \emptyset$.
By (2.7.1),
$L_2 \not\subset E$. Thus
$L_2$ contains at least one
irreducible curve, say $l''$,
which is contained in a 1-dimensional
fiber of $g$.
Then again by [Bea], $\dim L_2 =3$.

Combining this and (2.7.1),
we get (3).

\flushpar
{\it (4) \/} \ \
Let $l$ and $l'$ be two rational curves
of length 1 in $E$
and assume that $L_X(l) \not\subset E$.
We would like to deduce
$$
L_X(l') \not\subset E.
$$
In fact, if $L_X(l') \subset E$, then by (2),
we have
$$
\cases
l' \subset E_i, \\
\bigl( \widetilde{E}_i \, , \,
\nu_i^* \Cal O_{E_i}(-K_X) \bigr) \simeq
\bigl( \Bbb P^2 \, , \, \Cal O_{\Bbb P^2}(1) \bigr),
\text{ and} \\
\nu_i^* l' \text{ is a line in } \widetilde{E}_i
\endcases
\tag 2.7.2
$$
for some $i$.
On the other hand, $L_X(l) \supset E \supset E_i$
by (3). Since $E_i$ is an irreducible component
of a fiber of the contraction $g$
associated to the extremal ray $R$ and $[l] \in R$,
it follows that $E_i$ is actually an union
of some rational curves which are all deformations
of $l$ inside $X$. Since
$\bigl( \widetilde{E}_i \, , \,
\nu_i^* \Cal O_{E_i}(-K_X) \bigr) \simeq
\bigl( \Bbb P^2 \, , \, \Cal O_{\Bbb P^2}(1) \bigr)$
(2.7.2) and $(-K_X \, . \, l) =1$,
$\nu_i^* l$ is also a line in $\widetilde{E}_i$.
Again by (2.7.2), we conclude
that $l$ and $l'$ deform inside $X$
to each other,
in particular, $L_X(l) = L_X(l')$.
Since $L_X(l) \not\subset E$ and
$L_X(l') \subset E$, we get
a contradiction,
and (4) is proved.
\quad \qed
\enddemo

\proclaim{Theorem 2.8}\ \
For each $E_i$,
$\bigl( \widetilde{E}_i \, , \, \nu_i^* \Cal O_{E_i}(-K_X) \bigr)$
is isomorphic either to one of the followings:

\flushpar
(a) \, \, \,
$\bigl( \Bbb P^2 \, , \, \Cal O_{\Bbb P^2}(1) \bigr)$,

\flushpar
(b) \, \, \,
$\bigl( \Bbb P^2 \, , \, \Cal O_{\Bbb P^2}(2) \bigr)$,

\flushpar
(c) \, \, \,
$\bigl( \Bbb P^1 \times \Bbb P^1 \, , \,
\Cal O_{\Bbb P^1 \times \Bbb P^1}(1,m) \bigr)$ \, \,
$(m \geq 1)$,

\flushpar
(d) \, \, \,
$\bigl( \Sigma_m \, , \, \text{ an ample section} \bigr)$ \, \,
$(m \geq 1)$, or

\flushpar
(e) \, \, \,
$\bigl( S_m \, , \, \Cal O_{S_m}(1) \bigr)$ \, \, $(m \geq 2)$.

In particular, $\nu_i^* \Cal O_{E_i}(-K_X)$ is globally generated
for each $i$.
\endproclaim

\demo{Proof}\ \
We shall divide our proof into cases.

\flushpar
{\it Case (1) \/}\ \
Either $l_E(R)=2$, or
$l_E(R)=1$ and
$L_X(l) \subset E$
for each rational curve $l$ of length 1 in $E$.

In this case, if $l_{E_i}(R)=2$,
then by Lemma 2.7 (1),
$\widetilde{E}_i$ is of type (b).
If $l_{E_i}(R)=1$, then
by Lemma 2.7 (2),
$\widetilde{E}_i$ is of type (a).

\flushpar
{\it Case (2) \/}\ \
There is at least one rational curve $l$
of length 1 in $E$
such that $L_X(l) \not\subset E$.

In this case, by Lemma 2.7 (3),
$L_X(l) \supset E$, and in particular
each $E_i$ is covered by
rational curves of length 1 in $E_i$.
Hence by Proposition 1.4,
$\widetilde{E}_i$ is isomorphic either to
$\Bbb P^2$, a geometrically ruled surface,
or a cone obtained by contracting the negative section
of a geometrically ruled surface.
Moreover, each $E_i$ is a rational surface
by Proposition 2.6. Hence
$\widetilde{E}_i$
is of type either (a), (c), (d), or (e).
\quad \qed
\enddemo

\head \S 3.\ Rational 2-chain connectedness of $E$
\endhead

The aim of this section is to prove the following
theorem. We are inspired
this by Campana [C] and
Koll\'ar-Miyaoka-Mori [KoMiMo1].

\proclaim{Theorem 3.1}\ \
Either one of the followings holds;

\flushpar
(CL) \, \, \, (Connected by limit conics)

\flushpar
For any two points $x$, $y$ of $E$,
there exists a limit conic in $E$ which passes through
both $x$ and $y$, or

\flushpar
(MW) \, \, \, (Mukai-Wi\'sniewski type)

\flushpar
$E=E_1 \cup E_2$, $E_1 \simeq E_2 \simeq \Bbb P^2$,
$\dim (E_1 \cap E_2) =0$, and
$N_{E_1/X} \simeq N_{E_2/X} \simeq \Cal O_{\Bbb P^2}(-1)^{\oplus 2}$.
$g_U|_{U-E} : U \! - \! E \to V \! - \! P$
is a $\Bbb P^1\text{-bundle}$.
In particular \Rh $=2$.
In this case, for any points $x \in E_1$ and $y \in E_2$,
there exists a limit conic in $E$ which passes through
both $x$ and $y$.
\endproclaim

(In the case (MW), we will prove later
in Proposition 8.2 that $E_1 \cap E_2$
is actually a single point.)

\definition{Notation 3.2}\ \
Let $E_i$ be any irreducible component of $E$,
and $x$ a sufficiently general point
of $E_i$. Let $S_x \ni x$ be the connected component
of the intersection of two general
$g_U\text{-very}$ ample divisors in $U$
both of which pass through $x$, so that
$$
\cases
S_x \text{ is a smooth surface proper over } V, \text{ and} \\
S_x \cap E = \{x\} \text{ intersecting transversally.}
\endcases
\tag 3.2.1
$$
Let $\pi : Z \to X$ and $\phi : Z \to W$
be as in Notation 2.1, and let
$$
\cases
Z_U := \pi^{-1}(U), \, \, \, W_U := \phi(\pi^{-1}(U)), \\
\pi_U := \pi|_{Z_U} : Z_U \to U, \, \, \text{ and } \, \,
\phi_U := \phi|_{Z_U} : Z_U \to W_U.
\endcases
\tag 3.2.2
$$
Moreover let

\flushpar
(3.2.3) \qquad \qquad \qquad \qquad \quad \, \, \,
$D_x := \pi_U \, \phi_U^{-1} \, \phi_U \, \pi_U^{-1} (S_x)$.

\flushpar
This is disjoint from a general fiber
of $g_U : U \to V$.

We remark that
$\pi_U^{-1}(S_x)$ is an irreducible
surface in $Z_U$
by Bertini's theorem, as in (2.6.5).
Hence

\flushpar
(3.2.4) \, \, \,
$D_x$ (3.2.3) is a prime divisor in $U$.
\enddefinition

\proclaim{Proposition 3.3}\ \
Under the Notation 3.2,
assume that

\flushpar
\qquad \qquad \qquad \qquad \qquad \qquad \qquad \quad \, \,
$D_x \underset{g_U}\to {\equiv} 0 \, \text{ in } \, U$.

\flushpar
Then $x$ and any point of $E$
are contained in a same limit conic.
\endproclaim

\demo{Proof}\ \
If $D_x \not\supset E$, then
there is an irreducible curve
$C$ in $E$ such that $D_x \not\supset C$ and
$D_x \cap C \not= \emptyset$, in particular
$(D_x \, . \, C) >0$, which contradicts
the assumption $D_x \underset{g_U}\to {\equiv} 0$.
Hence

\flushpar
(3.3.1) \qquad \qquad \qquad \qquad \qquad \qquad \quad \, \,
$D_x \supset E$.

\flushpar
By the definition of $D_x$ (3.2.3),
$D_x$ is a union of curves
which are parametrized by $\phi$
and intersect with $S_x$. Since
$S_x \cap E = \{x\}$ (3.2.1),
$D_x \cap E$ is a union of
limit conics which pass through $x$.
Thus (3.3.1) implies that
$x$ and any point of $E$ are
contained in a same limit conic.
\quad \qed
\enddemo

\proclaim{Corollary 3.4}\ \
Assume

\flushpar
\qquad \qquad \qquad \qquad \qquad \qquad \quad \, \, \, \,
\Rh$=1$.

\flushpar
Then any two points of $E$ are contained in
a same limit conic in $E$.
\endproclaim

\demo{Proof}\ \
Let $x$ be as in Notation 3.2. Then
the assumption \Rh$=1$ implies
$D_x \underset{g_U}\to {\equiv} 0$,
since $D_x$ is disjoint from
a general fiber of $g_U : U \to V$.
Hence by Proposition 3.3,

\flushpar
(3.4.1) \, \, \,
$x$ and any point of $E$ are joined by a limit conic.

\flushpar
Since $x$ is chosen to be a general point of an
arbitrary irreducible component $E_i$ of $E$,
and since the property of being a limit conic is
preserved under degenerations,
the result follows.
\quad \qed
\enddemo

Next we shall prove Theorem 3.1
also in the case \Rh$\geq 2$.

\proclaim{Theorem 3.5}\ \
Assume that

\flushpar
\qquad \quad \, \, \, \,
\Rh$\geq 2$ \, and \,
$g_U|_{U - E} : U \! - \! E \to V \! - \! P$
is a $\Bbb P^1\text{-bundle}$.

\flushpar
Then $E$ satisfies the following:

\flushpar
(1) \, \, \,
$E=E_1 \cup E_2$, $E_1 \simeq E_2 \simeq \Bbb P^2$,
$\dim (E_1 \cap E_2) =0$, and

\flushpar
\qquad \qquad \qquad \qquad \qquad \quad \, \, \,
$N_{E_1/X} \simeq N_{E_2/X} \simeq \Cal O_{\Bbb P^2}(-1)^{\oplus 2}$.

\flushpar
(2) \, \, \, For any point $x$ of $E_1$ and any point
$y$ of $E_2$, there exists a limit conic which
passes through both $x$ and $y$.

\flushpar
(3) \, \, \, (Structure of the cone \NE $-I$)

\flushpar
\qquad \qquad \qquad \qquad \qquad \qquad \quad \, \, \, \,
\Rh$=2$.

\flushpar
Let $\{R_1, R_2\}$ be the set of extremal rays
of $\overline{NE}(U/V)$. Then $R_i$ defines
the flipping contraction [Ka4] which contracts
$E_i$ to a point $(i=1,2)$.
\endproclaim

\demo{Proof}\ \
First by the assumption, any rational curve
$l$ of length 1 in $E$
(if exists) satisfies $L_X(l) \subset E$.
Hence by Lemma 2.7 (1) and (2),

\flushpar
(3.5.0) \qquad \qquad \qquad \qquad \qquad \qquad
$\widetilde{E}_j \simeq \Bbb P^2$ \, \, \, $(\forall j)$,

\flushpar
regardless of $l_{E_j}(R)=1$ or $2$.

Consider the cone of curves $\overline{NE}(U/V)$.
This is spanned by finitely many
extremal rays, since $-K_U$ is $g_U\text{-ample}$.
In particular

\flushpar
(3.5.1) \, \, \, $\overline{NE}(U/V)$ has at least two
extremal rays.

\flushpar
Choose any extremal ray of $\overline{NE}(U/V)$ and
let $h : U \to U_h$ be the associated contraction.
Then we claim that

\flushpar
(3.5.2) \, \, \, $h$ never contracts
a general fiber of $g_U : U \to V$.

\flushpar
In fact, if $h$ contracts a general fiber of $g_U$,
then obviously $h$ contracts all limit conics in $E$.
Hence by Proposition 2.6 (2), $h(E)$ must be a point,
and $h$ coincides with $g_U$, a contradiction to our
assumption $\rho(U/V) \geq 2$.
Hence (3.5.2).

Since $g_U|_{U - E} : U \! - \! E \to V \! - \! P$
is a $\Bbb P^1\text{-bundle}$ by the assumption,
(3.5.2) implies that

\flushpar
(3.5.3) \, \, \, $h$ is a flipping contraction with
$\text{Exc }h \subset E$.

\flushpar
Then by Kawamata [Ka4],

\flushpar
(3.5.4) \, \, \, $h$ contracts an unique irreducible
component $E_i$ of $E$ such that

\flushpar
\qquad \qquad \qquad \qquad \qquad \, \, \, \, \,
$E_i \simeq \Bbb P^2$ \, and \,
$N_{E_i/X} \simeq \Cal O_{\Bbb P^2}(-1)^{\oplus 2}$.

\flushpar
{}From (3.5.1) and (3.5.4), there are at least two
irreducible components, say $E_1$ and $E_2$,
of $E$ such that

\flushpar
(3.5.5) \qquad \qquad
$E_1 \simeq E_2 \simeq \Bbb P^2$, \, and \,
$N_{E_1/X} \simeq N_{E_2/X} \simeq \Cal O_{\Bbb P^2}(-1)^{\oplus 2}$.

\flushpar
Let $h_i$ $(i=1,2)$ be the flipping contraction
which contracts $E_i$ to a point.

Assume that $E$ has at least three irreducible components,
to get a contradiction.

Let $E_3$ be any other one. Since
$\widetilde{E}_3 \simeq \Bbb P^2$ (3.5.0),
$h_i$ never contracts any curve in $E_j$
$(i \in \{1, 2\}, j \in \{1, 2, 3\} - \{i\})$.
Thus

\flushpar
(3.5.6) \qquad \qquad \quad \, \, \, \, \,
$\dim (E_i \cap E_j) \leq 0$ \, \, \,
$(\forall i, j \in \{1, 2, 3\}, \, i \not= j)$.

\flushpar
This contradicts the fact that any two limit conics
in $E$ deform inside $E$ to each other
(Remark 2.3 (3)).

Hence

\flushpar
(3.5.7) \qquad \qquad \qquad \qquad \qquad \qquad \,
$E = E_1 \cup E_2$.

\flushpar
Moreover by the same reason as in (3.5.6), we have

\flushpar
(3.5.8) \qquad \qquad \qquad \qquad \qquad \quad \,
$\dim (E_1 \cap E_2) =0$.

\flushpar
(3.5.5), (3.5.7) and (3.5.8) prove the Theorem.
\quad \qed
\enddemo

Next we investigate the structure
of the cone \NE in the rest case,
{\it i.e.\/} the case \Rh$\geq 2$ and
$g_U|_{U - E} : U \! - \! E \to V \! - \! P$
is not a $\Bbb P^1\text{-bundle}$.

\definition{Notation 3.6}\ \
Till the end of this section, for any curve $l$
which is contained in a fiber of $g_U$,
the symbol $[l]$ means the numerical equivalence
class of $l$ in $U$, not in $X$:
$[l] \in \overline{NE}(U/V)$.
\enddefinition

\proclaim{Theorem 3.7}\ \ (Structure of the cone \NE $-II$)

\flushpar
Assume that

\flushpar
\qquad \quad
\Rh$\geq 2$ \, and \,
$g_U|_{U - E} : U \! - \! E \to V \! - \! P$
is not a $\Bbb P^1\text{-bundle}$.

\flushpar
Then the followings hold:

\flushpar
(1)\, \, \, For each extremal ray $R_k$ of
$\overline{NE}(U/V)$,
there exists a rational curve $l_k$ of length 1
in $E$ such that
$$
R_k = \Bbb R_{\geq 0}[l_k].
$$

\flushpar
Moreover, $R_k$ defines a divisorial contraction
$\varphi_k : U \to U_k'$
which contracts an irreducible divisor $D_k$ to a surface, and
$(D_k \, . \, l_k)=-1$. Set theoretically, $D_k = L_U(l_k)$.

\flushpar
(2)\, \, \, Conversely, for each rational curve
$l$ of length 1 in $E$, $\Bbb R_{\geq 0}[l]$ forms
an extremal ray of
$\overline{NE}(U/V)$. If we regard $L_U(l)$
as a reduced divisor, then $-K_U + L_U(l)$ is a
supporting divisor of $\Bbb R_{\geq 0}[l]$.
\endproclaim

\demo{Proof}\ \
First by the assumption and Lemma 2.7 (4),

\flushpar
(3.7.0) \, \, \,
$L_U(l) \not\subset E$ for {\it any\/} rational curve
$l$ of length 1 in $E$.

\flushpar
{\it (1) \/}\ \
Take an arbitrary extremal ray $R_k$ of $\overline{NE}(U/V)$
and let $\varphi_k$ be the associated contraction.
First we claim:

\flushpar
(3.7.1) \, \, \, $\varphi_k$ is a divisorial contraction.

In fact, if $\varphi_k$ is of fiber type,
then by Proposition 2.6 (2),
it must coincide with $g_U$,
which is absurd.
Thus $\varphi_k$ is a birational contraction.
Assume that $\varphi_k$ is a flipping type
contraction.
Then by [Ka4],
$\varphi_k$ contracts an irreducible component
$E_i \simeq \Bbb P^2$ of $E$ to a point,
with a line $l$ of $E_i$ being of length 1.
Since $\varphi_k$ is an isomorphism
outside $E_i$, $l$ never
deforms outside $E$, which
contradicts (3.7.0).
Hence (3.7.1).

Moreover, since any divisor of $U$
is sent by $g_U$ either to
$V$ itself or to a divisor on $V$,
it follows that

\flushpar
(3.7.2) \, \, \, $\varphi_k$ contracts
an irreducible divisor $D_k$ to a surface.

Next, consider
$\varphi_k|_{D_k \! - \! D_k \cap E} : D_k \! - \! D_k \cap E \to
\varphi_k(D_k \! - \! D_k \cap E)$.
Since
$g|_{U-E} : U \! - \! E \to V \! - \! P$
is a conic bundle
(Proposition 2.2) and $\varphi_k$ never
contracts a whole fiber of $g_U$ as above,

\flushpar
(3.7.3) \, \, \,
$\varphi_k|_{D_k \! - \! D_k \cap E}$
is a $\Bbb P^1\text{-bundle}$
whose fiber $l_k$ satisfies
$$
N_{l_k/U} \simeq
\Cal O_{\Bbb P^1}^{\oplus 2} \oplus \Cal O_{\Bbb P^1}(-1), \,
\text{ and } \, (-K_U \, . \, l_k)=1.
$$
Moreover for such $l_k$, consider
the exact sequence:
$$
0 \longrightarrow N_{l_k/D_k} \longrightarrow N_{l_k/U}
\longrightarrow \Cal O_{D_k}(D_k)\otimes \Cal O_{l_k} \longrightarrow 0
$$

\flushpar
Since $\text{deg } N_{l_k/D_k} = 0$
and $\text{deg } N_{l_k/U} = -1$ (3.7.3),

\flushpar
(3.7.4) \qquad \qquad \qquad \qquad \qquad \qquad
$(D_k \, . \, l_k)=-1$.

By considering a degeneration
of $l_k \text{'s}$ into $E$, it follows that
$R_k$ is spanned by
a rational curve of length 1 which is
actually contained in $E$.
Hence we get (1).

\flushpar
{\it (2) \/}\ \
Assume that there exists a rational curve
$l$ of length 1 in $E$
such that $\Bbb R_{\geq 0}[l]$
is not an extremal ray of $\overline{NE}(U/V)$,
and shall derive a contradiction.

By (3.7.0), $L_U(l)$
is purely codimension 1. So we regard this
as a reduced divisor on $U$ and denote by $L$;
$$
L:= L_U(l).
$$
Then we claim:

\flushpar
(3.7.5)\, \, \, $-K_U + L$ is $g_U\text{-ample}$.

In fact by (1), any extremal ray of
$\overline{NE}(U/V)$ is written as
$$
R_k = \Bbb R_{\geq 0} [l_k]
$$
for a certain rational curve $l_k$ of length 1 in $E$.
Then again by (3.7.0), $l_k$ deforms to
a rational curve $l_k'$ which is
an irreducible component of
a degenerate fiber of the conic bundle
$g|_{U \! - \! E} : U \! - \! E \to V \! - \! P$.
If $L \supset l_k'$, then obviously
$l \underset{g_U}\to {\equiv} l_k'$,
which contradicts the assumption
that $\Bbb R_{\geq 0}[l]$ is not an extremal ray
of $\overline{NE}(U/V)$. Thus
$L \not\supset l_k'$, in particular
$$
(L \, . \, l_k) \geq 0.
$$
On the other hand, $(-K_U \, . \, l_k) >0$ for each $k$,
since $-K_U$ is $g_U\text{-ample}$.
Thus
$(-K_U+L \, . \, l_k) >0$ for each $k$.
This means that $-K_U + L$ is $g_U\text{-ample}$,
since $\overline{NE}(U/V) = \sum \Bbb R_{\geq 0}[l_k]$,
and (3.7.5) is proved.

On the other hand,

\flushpar
(3.7.6) \qquad \qquad \qquad \qquad \qquad \qquad \, \, \,
$(L \, . \, l)=-1$.

\flushpar
In fact,
let $l'$ be a deformation of $l$ which
is not contained in $E$ (3.7.0),
and is contained in an unique
irreducible component, say $L_1$, of $L$.
Then $l'$ is a fiber of the $\Bbb P^1\text{-bundle}$
$g|_{L_1 \! - \! L_1 \cap E} : L_1 \! - \! L_1 \cap E \to
g(L_1 \! - \! L_1 \cap E)$. Thus from
$$
0 \longrightarrow N_{l'/L_1} \longrightarrow N_{l'/U}
\longrightarrow \Cal O_{L_1}(L_1)\otimes \Cal O_{l'} \longrightarrow 0
$$
we get (3.7.6).

Since $(-K_U \, . \, l')=1$,
(3.7.6) says that $(-K_U + L \, . \, l')=0$, which contradicts
the above claim, and hence (2).
\quad \qed
\enddemo

Now we come to the proof
of Theorem 3.1 in the case
\Rh$\geq 2$ and
$g_U|_{U - E} : U \! - \! E \to V \! - \! P$
is not a $\Bbb P^1\text{-bundle}$.

\proclaim{Corollary 3.8}\ \
Assume that

\flushpar
\qquad \quad
\Rh$\geq 2$ \, and \,
$g_U|_{U - E} : U \! - \! E \to V \! - \! P$
is not a $\Bbb P^1\text{-bundle}$.

\flushpar
Then any two points of $E$ are
contained in a same limit conic.
\endproclaim

\demo{Proof}\ \
Let $E_i$ be any irreducible component of $E$, and
$x$ a general smooth point of $E_i$.
Let $S_x \ni x$ be as in Notation 3.2,
in particular be satisfying
the condition (3.2.1).
Let $\Delta^{\circ} \subset V-P$ be the discriminant
locus of the conic bundle
$g|_{U \! - \! E} : U \! - \! E \to V \! - \! P$
(Proposition 2.2),
$\Delta := \Delta^{\circ} \cup \{P\}$, and
$\Delta = \bigcup\limits_j \Delta_j$
the irreducible decomposition.
Note that

\flushpar
(3.8.0) \qquad \qquad \qquad \qquad \qquad \qquad \quad \, \, \, \,
$\Delta^{\circ} \not= \emptyset$

\flushpar
by the assumption.
Then

\flushpar
(3.8.1) \qquad \qquad \qquad \qquad \qquad \quad \,
$\dim \Delta_j =2$ \quad $(\forall j)$

\flushpar
[Bea]. Hence by the definition of $S_x$,

\flushpar
(3.8.2) \qquad \qquad \qquad \qquad \qquad \quad
$g(S_x) \not\supset \Delta_j$ \quad $(\forall j)$.

\flushpar
Let $D_x$ be as in (3.2.3). Recall that

\flushpar
(3.8.3) \, \, \, $D_x$ is a prime divisor in $U$

\flushpar
(3.2.4).
Take any extremal ray $R_k$ of $\overline{NE}(U/V)$.
This is spanned by a rational curve $l$ of length 1 in $U$
(Theorem 3.7).
Consider the whole deformation locus $L_U(l)$
of $l$ in $U$.
Then by Theorem 3.7,
$g(L_U(l))$ coincides with some $\Delta_j$.
Thus by (3.8.2), we can choose
a deformation of $l$ which is disjoint from $D_x$,
in particular $(D_x \, . \, l)=0$.
Since $\overline{NE}(U/V)$ is spanned by extremal rays,
it follows that

\flushpar
(3.8.4) \qquad \qquad \qquad \qquad \qquad \qquad
$D_x \underset{g_U}\to {\equiv} 0$ \, \, \, in \, \, \, $U$.

\flushpar
Thus the assumption of Proposition 3.3
is satisfied, and hence

\flushpar
(3.8.5) \, \, \,
$x$ and any point of $E$ are joined by a limit conic.

\flushpar
Since $x$ is chosen to be a general point of an
arbitrary irreducible component $E_i$ of $E$,
and since the property of being a limit conic is
preserved under degenerations,
the result follows.
\quad \qed
\enddemo

\proclaim{3.9}\ \ Proof of Theorem 3.1.
\endproclaim

Corollary 3.4, Theorem 3.5 and
Corollary 3.8 immediately imply Theorem 3.1.
\quad \qed

\head \S 4.\ Relative base-point-freeness of $|-K_X|$
\endhead

In this section,
we shall prove the following theorem,
whose proof heavily relies on
Kawamata's base-point-free technique,
especially ([Ka4], (2.3)):

\proclaim{Theorem 4.1} \ \
The linear system $|-K_X|$ is $g$-free.
Namely, the natural homomorphism
$$
\rho : g^*g_* \Cal O_X(-K_X) \to \Cal O_X(-K_X)
$$
is surjective.
\endproclaim

\definition{Definition 4.2}\ \
Let $B \subset Y$ be the whole set of points
which give 2-dimensional fibers of $g$,
as in the Notation 2.1.
Since $g$ is a conic bundle
over $Y \! - \! B$
(Proposition 2.2),
the above $\rho$ is surjective
on $g^{-1}(Y \! - \! B)$.
Thus we concentrate on
an neighborhood $U = g^{-1}(V)$ of
any 2-dimensional fiber $E$ of $g$.
We shall prove that
$\rho|_U$ is surjective.

For any torsion-free $\Cal O_V\text{-module}$
$\Cal G$, and for any section $v \in \Cal G$,
we define $Z_{\Cal G}(v)$,
{\it the zero locus of \/} $v$ {\it in \/} $\Cal G$,
as the closed analytic subspace of $V$ defined by the
image ideal of the $\Cal O_V\text{-homomorphism}$:
$$
{}^av : \Cal G^{\vee} \to \Cal O_V, \quad
\Cal G^{\vee} \ni \varphi \mapsto \varphi(v) \in \Cal O_V.
$$
Moreover let
$$
B_{\Cal G} :=
\{Q \in V | \, \Cal G_Q \text{ is not a free }
\Cal O_{V,Q}\text{-module} \}.
$$
Note that

\flushpar
(4.2.1) \qquad \quad
$\dim B_{\Cal G} \leq 1$, \, and \,
$\{ Q \in V | \text{Tor}_1^{\Cal O_{V,Q}}
(\Cal G_Q, \Bbb C(Q))\not= 0 \} \subset B$,

\flushpar
where $\Bbb C(Q) = \Cal O_{V,Q}/\frak m_Q$
is the residue field of $\Cal O_{V,Q}$.
If $Q \in V - B_{\Cal G}$,
then $Q \in Z_{\Cal G}(v)$
if and only if $v(Q)=0$ in $\Bbb C(Q)$.
\enddefinition

\proclaim{Lemma 4.3}\ \ (Communicated by
Y.Kawamata and N.Takahashi)

For any irreducible member $D \in |-K_U|$, there exists
another $D' \in |-K_U|$ such that

\flushpar
(1) \, \, \, $D$ and $D'$ give the same element of
$g_{U *}\Cal O_U(-K_U) \otimes \Bbb C(P)$, and

\flushpar
(2) \, \, \, $D'$ never contains a whole fiber of
$g_U|_{U \! - \! E} : U \! - \! E \to V \! - \! P$,
as analytic subspaces of $U$.
\endproclaim

\demo{Proof}\ \
Let $D \in |-K_U|$ be any irreducible member.
If the condition (2) of this lemma fails for
$D'=D$, then from $D$
we shall find another $D' \in |-K_U|$ which satisfies
both (1) and (2).

Let $\Cal F_3 := (g_U)_* \Cal O_U(-K_U)$.
Since $g_U|_{U \! - \! E} : U \! - \! E \to V \! - \! P$
is a conic bundle (Proposition 2.2),
it follows that
$\dim H^0(U_Q \, , \Cal O_U(-K_U) \otimes \Cal O_{U_Q})=3$
for any $Q \in V \! - \! P$,
and thus $\Cal F_3$ is a reflexive sheaf of rank 3
on $V$ which is locally free outside $P$.
Let $s \in \Cal F_3$ be the section
corresponding to $D$. Then

\flushpar
(4.3.0)\, \, \, For any $Q \in V \! - \! P$,
$Q \in Z_{\Cal F_3}(s)$ is equivalent to
$D \supset U_Q$ as analytic subspaces of $U$.

\flushpar
In particular, $Z_{\Cal F_3}(s)$
is at most 1-dimensional, since
$g_U|_D : D \to V$ is generically finite,
and since $D$ is irreducible.

If $\dim Z_{\Cal F_3}(s) =0$, then by (4.3.0) there is
nothing to prove, so assume $\dim Z_{\Cal F_3}(s)$ $=1$.
By (4.3.0) again,

\flushpar
(4.3.1)\ \ \  For the lemma, it is enough
to to show the existence of
$t \in \frak m_P\Cal F_3$
such that
$Z_{\Cal F_3}(s + t) \subset \{P\}$.

Take first of all a saturated filtration
$$
\cases
0 \subset \Cal F_1 \subset \Cal F_2 \subset \Cal F_3 \quad \text{with }
\text{rk } \Cal F_i = i \quad \text{such that} \\
s \not\in \Cal F_2, \, \, \, \text{ and} \\
Z_{\Cal F_3}(s) \cap B \subset \{P\},
\endcases
\tag 4.3.2
$$
where $B := B_{\Cal F_3/\Cal F_1}$, as defined in Definition 4.2.
Note that

\flushpar
(4.3.3) \qquad
$\dim B \leq 1$ \, and \,
$\{ Q \in V | \text{Tor}_1^{\Cal O_{V,Q}}
((\Cal F_3/\Cal F_1)_Q, \Bbb C(Q))\not= 0 \}
\subset B$

\flushpar
(4.2.1). Then we claim:
\enddemo

\proclaim{claim 1}\ \
There exists
$t \in \frak m_P(\Cal F_2/\Cal F_1)$ such that
for any $\varepsilon \in \Bbb C$ with $0<|\varepsilon|<1$,
$$
\dim Z_{\Cal F_3/\Cal F_1}(s+ \varepsilon t) \leq 1.
$$
\endproclaim

\demo{Proof}\ \
If $\dim Z_{\Cal F_3/\Cal F_1}(s) \leq 1$,
then the claim trivially holds by letting $t=0$.
Assume
$\dim Z_{\Cal F_3/\Cal F_1}(s) =2$.
First we take
$t \in \frak m_P(\Cal F_2/\Cal F_1)$ so that

\flushpar
(4.3.4) \qquad \qquad \qquad \qquad
$\dim \bigl( Z_{\Cal F_2/\Cal F_1}(t)
\cap Z_{\Cal F_3/\Cal F_1}(s) \bigr) =1$.

\flushpar
Let $B':= B_{\Cal F_3/\Cal F_2}$.
Then
$$
\dim B' \leq 1 \, \text{ and } \,
\{ Q \in V | \text{Tor}_1^{\Cal O_{V,Q}}
((\Cal F_3/\Cal F_2)_Q , \Bbb C(Q)) \not= 0 \} \subset B'
$$
(4.2.1). In particular,
$$
Z_{\Cal F_2/\Cal F_1}(t) \subset
Z_{\Cal F_3/\Cal F_1}(t) \subset Z_{\Cal F_2/\Cal F_1}(t) \cup B'.
$$
{}From this and (4.3.4),

\flushpar
(4.3.5) \quad \qquad \qquad \qquad \quad
$\dim \bigl( Z_{\Cal F_3/\Cal F_1}(t)
\cap Z_{\Cal F_3/\Cal F_1}(s) \bigr) =1$.

Second, for any $\varepsilon \in \Bbb C$,

\flushpar
(4.3.6) \qquad \qquad \quad \, \,
$Z_{\Cal F_3/\Cal F_1}(s+ \varepsilon t)
\subset Z_{\Cal F_3/\Cal F_2}(s+ \varepsilon t)
= Z_{\Cal F_3/\Cal F_2}(s)$

\flushpar
by $t \equiv 0$ in $\Cal F_3/\Cal F_2$. Let $\{V_k\}_k$ be the
set of all 2-dimensional irreducible subspaces $V_k$ of $V$
passing through $P$ such that
$$
V_k \not\subset Z_{\Cal F_3/\Cal F_1}(s) \text{ and }
V_k \subset Z_{\Cal F_3/\Cal F_2}(s).
$$
In particular,

\flushpar
(4.3.7) \quad $Z_{\Cal F_3/\Cal F_2}(s) =
\bigcup\limits_k V_k \cup Z_{\Cal F_3/\Cal F_1}(s) \cup
(\text{one dimensional subspaces of } V)$.

\flushpar
Let
$$
r_k := \text{Inf} \bigl\{ r \in \Bbb R_{>0} | \, V_k \subset
Z_{\Cal F_3/\Cal F_1}(s+ \varepsilon t) \text{ for some }
\varepsilon \in \Bbb C -\{0\} \text{ with } |\varepsilon|=r \bigr\}.
$$

\flushpar
Then obviously $r_k>0$ for each $k$,
and for any $\varepsilon \in \Bbb C$
with $0< |\varepsilon| < \operatornamewithlimits{Min}_{k} r_k$,
we have

\flushpar
(4.3.8) \qquad \qquad \qquad \qquad
$Z_{\Cal F_3/\Cal F_1}(s+ \varepsilon t) \not\supset V_k$ \, \, \,
for each $k$.

\flushpar
Thus
$$
\align
Z_{\Cal F_3/\Cal F_1}(s+ \varepsilon t)
= & Z_{\Cal F_3/\Cal F_1}(s+ \varepsilon t) \cap
Z_{\Cal F_3/\Cal F_2}(s)
\qquad \qquad \qquad \qquad \, \,
(\text{by } (4.3.6)) \\
= & \bigcup\limits_k \bigl( Z_{\Cal F_3/\Cal F_1}
(s+ \varepsilon t) \cap V_k \bigr)
\, \, \cup \, \bigl( Z_{\Cal F_3/\Cal F_1}(s+ \varepsilon t) \cap
Z_{\Cal F_3/\Cal F_1}(s) \bigr) \\
& \, \, \cup \bigl( Z_{\Cal F_3/\Cal F_1}(s+ \varepsilon t) \cap
( \text{one dimensional subspace of }V ) \bigr) \\
& \qquad \qquad \qquad \qquad \qquad \qquad
\qquad \qquad \qquad \qquad
(\text{by } (4.3.7)) \\
= & \bigcup\limits_k \bigl(
Z_{\Cal F_3/\Cal F_1}(s+ \varepsilon t) \cap V_k \bigr)
\, \, \cup \, \bigl( (Z_{\Cal F_3/\Cal F_1}(t)
\cap Z_{\Cal F_3/\Cal F_1}(s) \bigr) \\
& \, \, \cup \bigl( \text{one dimensional
subspace of } V \bigr),
\endalign
$$
which is of dimension at most 1 (4.3.5),(4.3.8).
Hence we get claim 1.
\quad \qed
\enddemo

Next we prove the following:

\proclaim{claim 2}\ \
Let $t \in \frak m_P(\Cal F_2/\Cal F_1)$ be as given in
claim 1 above. Let $t' \in \frak m_P\Cal F_2$ be an
arbitrary lifting of $t$. Then for any $\varepsilon \in \Bbb C$
with $0<|\varepsilon|<1$,
there is a suitable
$u_{\varepsilon} \in \frak m_P\Cal F_1$
such that for any $\eta \in \Bbb C$ with
$0<|\eta|<1$,
$$
Z_{\Cal F_3}(s+ \varepsilon t'+ \eta u_{\varepsilon}) \subset B,
$$
where $B = B_{\Cal F_3/\Cal F_1}$ as in (4.3.2).
\endproclaim

\demo{Proof}\ \
Since
$$
Z_{\Cal F_3}(s+\varepsilon t')
\subset Z_{\Cal F_3/\Cal F_1}(s+\varepsilon t),
$$
we have
$\dim Z_{\Cal F_3}(s+\varepsilon t') \leq 1$ by claim 1.
If $Z_{\Cal F_3}(s+\varepsilon t') \subset B$,
then there is nothing to prove.
Assume $Z_{\Cal F_3}(s+\varepsilon t') \not\subset B$.
First we can take
$u=u_{\varepsilon} \in \frak m_P \Cal F_1$,
depending on the value $\varepsilon$, so that

\flushpar
(4.3.9) \qquad \qquad \qquad \qquad \quad
$Z_{\Cal F_1}(u) \cap Z_{\Cal F_3}(s+\varepsilon t') = \{P\}$.

\flushpar
(Recall that $\Cal F_1$ is a reflexive sheaf of rank 1 on $V$.)
Since $Z_{\Cal F_3}(u) \subset Z_{\Cal F_1}(u) \cup B$ (4.3.3),
we have

\flushpar
(4.3.10) \qquad \qquad \qquad \qquad \quad
$Z_{\Cal F_3}(u) \cap Z_{\Cal F_3}(s+\varepsilon t') \subset B$

\flushpar
by (4.3.9). Second, for any $\eta \in \Bbb C$,

\flushpar
(4.3.11) \qquad \, \, \,
$Z_{\Cal F_3}(s+\varepsilon t'+ \eta u) \subset
Z_{\Cal F_3/\Cal F_1}(s+\varepsilon t'+ \eta u)
= Z_{\Cal F_3/\Cal F_1}(s+\varepsilon t)$

\flushpar
by $u \equiv 0$ in $\Cal F_3/\Cal F_1$.
Let $\{W_k\}_k$ be the set of all
1-dimensional irreducible subspaces of $V$ passing through
$P$ such that
$$
W_k \not\subset Z_{\Cal F_3}(s+\varepsilon t') \text{ and }
W_k \subset Z_{\Cal F_3/\Cal F_1}(s+\varepsilon t),
$$
and in particular

\flushpar
(4.3.12) \qquad \qquad \qquad \, \, \, \,
$Z_{\Cal F_3/\Cal F_1}(s+\varepsilon t)
= \bigcup\limits_k W_k \, \cup Z_{\Cal F_3}(s+\varepsilon t')$.

\flushpar
Let
$$
q_k := \text{Inf} \bigl\{ q \in \Bbb R_{>0} | \, W_k \subset
Z_{\Cal F_3}(s+\varepsilon t'+ \eta u) \text{ for some }
\eta \in \Bbb C -\{0\} \text{ with } |\eta|=q \bigr\}.
$$

Then obviously $q_k>0$ for each $k$,
and for any $\eta \in \Bbb C$
with $0< \eta < \operatornamewithlimits{Min}_{k} q_k$,
we have

\flushpar
(4.3.13) \qquad \qquad \quad \, \, \,
$Z_{\Cal F_3}(s+\varepsilon t'+ \eta u) \cap W_k = \{P\}$ \, \, \,
for each $k$.

\flushpar
Thus
$$
\align
Z_{\Cal F_3}(s+\varepsilon t'+ \eta u)
= & Z_{\Cal F_3}(s+\varepsilon t'+ \eta u) \cap
Z_{\Cal F_3/\Cal F_1}(s+\varepsilon t)
\qquad \qquad
(\text{by } (4.3.11) ) \\
= & \bigcup\limits_k \bigl( Z_{\Cal F_3}
(s+\varepsilon t'+ \eta u) \cap W_k \bigr) \\
& \qquad \qquad \qquad \qquad \, \,
\cup \, \bigl( Z_{\Cal F_3}(s+\varepsilon t'+ \eta u) \cap
Z_{\Cal F_3}(s+\varepsilon t') \bigr) \\
& \qquad \qquad \qquad \qquad \qquad \qquad \qquad
\qquad \qquad \, \, \, \, \,
(\text{by } (4.3.12)) \\
= & \bigcup\limits_k \bigl(
Z_{\Cal F_3}(s+\varepsilon t'+ \eta u) \cap W_k \bigr)
\, \, \cup \, \bigl( (Z_{\Cal F_3}(s+\varepsilon t')
\cap Z_{\Cal F_3}(u) \bigr) \\
\subset & \, \, B, \qquad \qquad \qquad \qquad
\qquad \qquad \quad \, \, \, \,
(\text{by } (4.3.10), (4.3.13))
\endalign
$$
and the claim 2 is proved. \quad \qed

Finally, if we choose a sufficiently small
$\varepsilon >0$, and then choose a sufficiently small
$\eta >0$,
$Z_{\Cal F_3}(s+ \varepsilon t'
+ \eta u_{\varepsilon}) \cap B \subset \{P\}$
by (4.3.2). From this together with claim 2 above,
we actually have
$Z_{\Cal F_3}(s+ \varepsilon t'
+ \eta u_{\varepsilon}) \subset \{P\}$.
Hence we get (4.3.1), and the proof of
Lemma 4.3 is completed.
\quad \qed
\enddemo

{}From now on,
assume that the cokernel $\text{Bs } |-K_U|$
of $\rho|_U$ is nonempty so as to derive a
contradiction.

\proclaim{Lemma 4.4} \ \
A general member $D$ of $|-K_U|$
has at most canonical singularities.
\endproclaim

\demo{Proof} \ \
We will follow the argument of
the proof of the base-point-free theorem
originated by Kawamata.

We take $\varphi : U' \to U$,
a resolution of the base locus
$\text{Bs }|-K_U|$ of $|-K_U|$ such that
$$
\cases
\varphi ^* |-K_U|=|D'|+ \sum r_iG_i,\\
K_{U'}=\varphi ^* K_U+\sum a_iG_i,\\
-\varphi ^* K_U-\sum \delta _i G_i \quad
\text{is $g \circ \varphi$-ample,}
\endcases
\tag 4.4.1
$$
where $\sum\limits_{i=1}^N G_i$
is a simple normal crossing divisor
on $U'$, $D'$ is the proper transform of $D$
such that

\flushpar
(4.4.2) \qquad \qquad \qquad \qquad \qquad \qquad \, \, \,
$\text{Bs }|D'|=\emptyset$,

\flushpar
$r_i, a_i \in \Bbb Z_{ \geq 0}$
with $(r_i, a_i) \not= (0,0)$,
and $\delta_i \in \Bbb Q_{>0}$ with $0 < \delta_i < 1$.
Note that

\flushpar
(4.4.3) \qquad \qquad \qquad \qquad \qquad \, \,
$\varphi(\bigcup\limits_i G_i) = \text{Bs }|-K_U|$.

\flushpar
Let
$$
c:=\min \dfrac{a_i+1- \delta_i}{r_i}.
$$
By openness of the condition of ampleness,
we may shrink $\delta _i$'s
if necessary so that the minimum $c$ is
attained exactly for a single $i$,
say $i=1$.
Let
$$
\cases
A:= \sum\limits_{i \geq 2} \left(-cr_i+a_i- \delta_i \right) G_i,
\, \, \text{ and} \\
B:= G_1.
\endcases
$$

Here we note
$-cr_i+a_i- \delta_i \geq -1$,
with equality if and only if $i=1$.
In particular,

$$
\cases
\ulcorner A \urcorner \geq 0, \\
\text{Supp } \ulcorner A \urcorner \subset
\text{Exc } \varphi, \text{ and} \\
A-B= \sum\limits_{i \geq 1}
\left( -cr_i+a_i- \delta_i \right)G_i.
\endcases
\tag 4.4.4
$$
Then
\enddemo

\proclaim{claim} \ \
$ a_i\geq r_i$ for all $i$.
\endproclaim

\demo{Proof} \ \
If $a_i+1 \leq r_i$ for some $i$,
then $c< \dfrac{a_i+1}{r_i} \leq 1$ {\it i.e. \/}
$2-c >1$. Moreover since $D'$ is nef (4.4.2) and $c>0$,
$$
\align
C :=&- \varphi^*K_U-K_{U'}+(A-B) \\
  =&cD'-(2-c)\varphi^*K_U-\sum_{i \geq1} \delta_iG_i
\endalign
$$
is $g_U \circ \varphi$-ample. Thus
$$
R^1(g_U \circ \varphi)_*\Cal O_{U'}
 \left( -\varphi^*K_U+\ulcorner A \urcorner -B \right) =0
$$
by the Kawamata-Viehweg vanishing theorem
(e.g.[KaMaMa]), and hence the following
restriction homomorphism is surjective:
$$
s : (g_U \circ \varphi)_*\Cal O_{U'}
\left(-\varphi^*K_U+ \ulcorner A \urcorner \right)
\longrightarrow
H^0 \bigl( B, \,
\Cal O_B (-\varphi^*K_U + \ulcorner A \urcorner) \bigr)
$$
The left-hand side is naturally
equal to $g_{U,*} \Cal O_U(-K_U)$ by
$\varphi_* \Cal O_{U'}(\ulcorner A \urcorner) \simeq \Cal O_U$
(4.4.4).
On the other hand,
since $\nu_i^* \Cal O_{E_i}(-K_U)$
is globally generated for each
$i$ (Theorem 2.8), and since
$\varphi(B) \subset \text{Bs }|-K_U| \subset E$
(4.4.3), the right-hand side of
$s$ does not vanish. Hence
$\varphi(B) \not\subset \text{Bs }|-K_U|$,
which contradicts (4.4.3),
and the claim is proved.
\qed
\enddemo

\demo{Proof of Lemma 4.4 continued} \ \
Now we have
$a_i-r_i \geq 0$.
Let $D'$ be a general member of $|D'|$
by abuse of notation. $D'$ is smooth by (4.4.2)
and Bertini's theorem.
Then by (4.4.1) together with the adjunction,
we have
$$
K_{D'}=\varphi^*K_D+\sum_{i \geq 1}(a_i-r_i)(G_i|_{D'}),
$$
which implies that
$D$ has at most canonical singularities.
\qed
\enddemo

Let $D$ be a general member of $|-K_U|$
as in Lemma 4.4, and
$h:D \to \widetilde{V}$ the Stein
factorization of
$g_U|_D:D \to V$.
Then $h$ is a projective bimeromorphic
morphism.
Let $F$ be its
exceptional locus.

\proclaim{Lemma 4.5}\ \
Let
$F=\sum\limits_i F_{1,i}+\sum\limits_j F_{2,j}+\sum\limits_k F_{2,k}'$
be the irreducible decomposition of $F$, where
$$
\cases
\dim F_{1,i} =1, \\
\bigl( \dim F_{2,j}, \dim h(F_{2,j}) \bigr) = (2, 0), \text{ and} \\
\bigl( \dim F_{2,k}', \dim h(F_{2,k}') \bigr) = (2, 1).
\endcases
$$
Let
$F_1 := \sum\limits_i F_{1,i}$,
$F_2 := \sum\limits_j F_{2,j}$, and
$F_2':= \sum\limits_k F_{2,k}'$.
Then
$$
\text{Bs }|-K_U| =F_2.
$$
\endproclaim

(Note that each $F_{2,j}$ coincides with an
irreducible component of $E$.)

\demo{Proof} \ \
Let $L :=D|_D$,
then by the exact sequence
$$
0 \longrightarrow \Cal O_U \longrightarrow \Cal O_U(D)
\longrightarrow \Cal O_D(L) \longrightarrow 0
$$
together with $R^1 g_{U \, *} \Cal O_U =0$ [KaMaMa],
we have

\flushpar
(4.5.0) \qquad \qquad \qquad \qquad \qquad \quad \, \, \,
$\text{Bs }|D|=\text{Bs }|L|$.

\flushpar
Let us choose
$x \in (F_1 \cup F_2') \cap E -F_2$ arbitrarily.
By Lemma 4.3, if two 1-dimensional
irreducible components, say $l_1$ and $l_2$,
of a fiber of $h$ meets at $x$,
then $l_1 \cup l_2$ never deforms outside $E$ in $D$.
So we can find an effective Cartier divisor $M$ of $D$
with $M \ni x$ such that
$$
\cases
& \! \! \! \! \! \! (4.5.1) \quad M \cap (\text{the fiber
of } h \text{ passing through } x) = \{x\},
\text{ and} \\
& \! \! \! \! \! \! (4.5.2) \quad (M \, . \, l) =1
\text{ for any irreducible component }
l \text { of any one-dimensional} \\
&\qquad \quad \,
\text{fiber of } h \text{ with } M \cap l \not= \phi.
\endcases
$$
Hence $L-M$ is $h$-nef,
so again by [loc.cit],
$R^1h_* \Cal O_D(L-M)=0$.
In particular
$$
h_* \Cal O_D(L) \to h_*\Cal O_M(L)
$$
is surjective.
Thus
$x \not\in \text{Bs }|L| = \text{Bs }|-K_U|$,
and we have
$\text{Bs }|-K_U| \subset F_2$.

On the other hand,
$F_2$ is a union of some irreducible components of $E$
(as remarked just before this proof),
in particular,
$F_2$ does not depend on the choice of
a general $D \in |-K_U|$.
Thus we necessarily have
$F_2 \subset \text{Bs }|-K_U|$.
Hence
\linebreak
$\text{Bs }|-K_U|=F_2$.
\qed
\enddemo

\proclaim{4.6}\ \ Proof of Theorem 4.1.
\endproclaim

For a general
$h$-very ample divisor $H$ on $D$,
we have $H = H_1 + H_2$,
$H_1 \cap (F_1 \cup F_2') = H_2 \cap F_2 = \emptyset$
by shrinking $\widetilde{V}$
if necessary.
Then $|H_2|$ gives a projective bimeromorphic morphism
$h_1 : D \to \widetilde{V}'$
whose exceptional locus is just $F_2$.
Let $F^{\circ}$ be a connected component
of $F_2$, $D^{\circ}$ a neighborhood
of $F^{\circ}$ in $D$, $V^{\circ} := h_1(D^{\circ})$,
and consider
$$
h^{\circ} := h_1|_{D^{\circ}} : D^{\circ} \to V^{\circ}.
$$

\flushpar
(4.6.0) \, \, \,
$h^{\circ}$ is a projective bimeromorphic morphism
which contracts a connected exceptional divisor $F^{\circ}$
to a point,
and
each irreducible component
of $F^{\circ}$ coincides with some $E_i$.

\flushpar
Moreover, let us denote
$L|_{D^{\circ}}$ again by $L$ for simplicity.
Then by (4.5.0) and Bertini's theorem,

\flushpar
(4.6.1) \qquad \qquad \qquad \qquad \qquad \quad \,
$\text{Sing }D^{\circ} \subset \text{Bs }|L|$.

Let $h^{\circ}(F^{\circ}) =: P^{\circ} \in V^{\circ}$.
Let $L_0 \in |L|$ be a general member,
$r \in \Bbb Z$ with $r \gg 0$,
and $L'$ a general hyperplane section
of $(V^{\circ}, P^{\circ})$.
Let

\flushpar
(4.6.2) \qquad \qquad \qquad \qquad \qquad
$L_r :=L_0 + r \, h^{\circ *} L' \in |L|$.

We take a resolution
$\psi : D' \to D^{\circ}$
of singularities of $D^{\circ}$ such that
$$
\cases
\psi^* L_r= \sum r_i' G_i',\\
K_{D'} = \psi^* K_{D^{\circ}} + \sum a_i' G_i'
\, \, \sim \, \, \sum a_i' G_i', \, \text{ and} \\
\psi^* L - \sum \delta_i' G_i'
\quad \text{is $h^{\circ} \circ \psi$-ample,}
\endcases
\tag 4.6.3
$$
where $\sum G_i'$ is
a simple normal crossing divisor on $D'$,
$r_i'$, $a_i' \in \Bbb Z$ with $r_i' \geq 0$,
and $\delta_i' \in \Bbb Q_{>0}$
with $0 < \delta_i' < 1$.
Note that

\flushpar
(4.6.4) \qquad \qquad \qquad \qquad
$\psi(\bigcup\limits_i G_i') \subset
\text{Bs }|L| \cup \text{Sing }D^{\circ} = F^{\circ}$

\flushpar
by (4.6.1) and Lemma 4.5. Moreover by Lemma 4.4,

\flushpar
(4.6.5) \qquad \qquad \qquad \qquad \qquad \qquad \qquad
$a_i' \geq 0$

\flushpar
for each $i$. Then, for any $G_i'$,
we have $r_i' \geq r$
by virtue of the definition of $L_r$ (4.6.2),
and $\psi(G_i') \subset F^{\circ}$ (4.6.4).
So if we define
$$
c':=\min \dfrac{a_i'+1- \delta_i'}{r_i'},
$$
then we have

\flushpar
(4.6.6) \qquad \qquad \qquad \qquad \qquad \qquad \, \, \, \,
$0 < c' \ll 1$,

\flushpar
since $a_i'$ does not depend on
the choice of the number $r$.
We may assume that
the minimum $c'$ is attained
exactly for a single $i$, say $i=1$,
by shrinking $\delta_i'$'s
if necessary. Let
$$
\cases
A':=\sum\limits_{i \geq 2}(-c'r_i'+a_i'- \delta_i')G_i', \text{ and} \\
B':=G_1'.
\endcases
$$
Then
$$
\ulcorner A' \urcorner \geq 0\, ,
\,\, A'-B'=\sum_{i \geq 1}(-c'r_i'+a_i'- \delta_i')G_i'
$$
and
$$
\align
C':=\psi^* L-K_{D'}+(A'-B')
&=(1-c')\psi^*L-\sum \delta_i' G_i' \\
&=(1-c')\left( \psi^*L - \sum \dfrac{\delta_i'}{1-c'}G_i' \right)
\endalign
$$
is $h^{\circ} \circ \psi\text{-ample}$ (4.6.3), (4.6.6).
Thus
$R^1(h^{\circ} \circ \psi)_*
\Cal O_{D'} (\psi^*L + \ulcorner A' \urcorner - B')=0$
[loc.cit], and the restriction homomorphism
$$
(h^{\circ} \circ \psi)_* \Cal O_{D'}
(\psi^*L + \ulcorner A' \urcorner)
= {h^{\circ}}_* \Cal O_{D^{\circ}}(L)
\longrightarrow
H^0 \bigl( B', \,
\Cal O_{B'}(\psi^* L+ \ulcorner A' \urcorner) \bigr)
$$
is surjective.
Moreover, since
$\psi(B') \subset E$ (4.6.4),
and since each $\nu_i^* \Cal O_{E_i}(L)$ is globally
generated (Theorem 2.8),
$$
H^0 \bigl( B', \,
\Cal O_{B'}(\psi^* L+ \ulcorner A' \urcorner) \bigr) \not= 0,
$$
as in the proof of Lemma 4.4.
Thus we obtain
$$
\psi(B') \not\subset \text{Bs }|L|.
$$
This contradicts
$\psi(B') \subset \text{Bs }|L|$
(4.6.4),
and the proof of Theorem 4.1 is completed.
\quad \qed

\head \S 5.\ The case $l_{E_i}(R)=2$ for some $i$
\endhead

\definition{Assumption 5.0}\ \
In this section,
we assume that

\flushpar
(5.0.1) \, \, \, $E$ contains an irreducible component, say
$E_1$, with $l_{E_1}(R)=2$.
\enddefinition

The following is the main result of
this section:

\proclaim{Theorem 5.1}\ \
Under the Assumption 5.0,

\flushpar
(I) \, \, \, $E$ is irreducible; $E=E_1$,
and is isomorphic to $\Bbb P^2$.
$\Cal O_E(-K_X) \simeq$ $\Cal O_{\Bbb P^2}(2)$.
$Y$ is smooth at $P$.
In particular, \Rh$=1$.

\flushpar
(II) \, \, \,
The normal bundle $N_{E/X}$ of $E$ in $X$
is isomorphic to $\Omega_{\Bbb P^2}^1(1)$.

\flushpar
(III)\, \, \, (local elementary transformation)

Let $x \in E$ be
an arbitrary point. Then there exists a smooth surface
$S_x \subset U$ proper over $V$ such that
$g|_{S_x - x} : S_x \! - \! x \to g(S_x) \! - \! P$
is an isomorphism, and that $S_x \cap E=\{x\}$
intersecting transversally.
Let $\varphi : \overline{U} \to U$
be the blow-up with center $S_x$.
Then $-K_{\overline{U}}$ is
$g_U \circ \varphi\text{-ample}$,
and $\varphi^{-1}(E) \simeq \Sigma_1$.
Let $\varphi^+ : \overline{U} \to U^+$
be the contraction
associated to the extremal ray
of $\overline{NE}(\overline{U}/V)$
other than $\varphi$.
Then $U^+ \simeq V \times \Bbb P^1$.
Let $g^+: U^+ \to V$ be the first projection.
Then
$$
g_U \circ \varphi = g^+ \circ \varphi^+,
$$
and
$\varphi^+$ is the blow-up with center $S^+$
which is a smooth surface proper over $V$ such that
$S^+ \supset g^{+ \, -1}(P) \simeq \Bbb P^1$
with the normal bundle
$N_{g^{+ \, -1}(P)/S^+} \simeq \Cal O_{\Bbb P^1}(-1)$.
\endproclaim

First we shall prove:

\proclaim{Lemma 5.2}\ \
$g|_{U \! - \! E} : U \! - \! E \to V \! - \! P$
is a $\Bbb P^1\text{-bundle}$.
\endproclaim

\demo{Proof}\ \
Assume to the contrary. Then
there is a rational curve $l$ of length 1
which deforms outside $E$.
Then by Lemma 2.7 (3),
$L_X(l) \supset E$,
which contradicts the assumption (5.0.1).
\quad \qed
\enddemo

\proclaim{Lemma 5.3} \qquad \, \,
\Rh $=1$. \quad (See Definition 0.4.)
\endproclaim

\demo{Proof}\ \
If \Rh $\geq 2$,
then by Lemma 5.2 and
Theorem 3.5, $E$ is a union of
two $\Bbb P^2\text{'s}$
whose lines are of length 1.
This contradicts our
assumption (5.0.1).
\quad \qed
\enddemo

\proclaim{Lemma 5.4}\ \
Under the Assumption 5.0, $E$ is irreducible: $E=E_1$.
\endproclaim

\demo{Proof}\ \
By Lemma 5.3 and Theorem 3.1,
any two points of $E$ are
joined by a limit conic.
Since $l_{E_1}(R)=2$,
$E$ is necessarily irreducible.
\quad \qed
\enddemo

\proclaim{Lemma 5.5}\ \
Let $x$ be a general smooth point of $E$,
and let $S_x$ be as in Notation 3.2.
Let $\varphi : \overline{U} \to U$
be the blow-up of $U$ with center $S_x$.
Then $-K_{\overline{U}}$ is
$g_U \circ \varphi\text{-ample}$.
\endproclaim

\demo{Proof}\ \
First by (3.2.1),

\flushpar
(5.5.0) \, \, \, For every one-dimensional fiber $C$
of $\varphi$, $C \simeq \Bbb P^1$ and
$(-K_{\overline{U}} \, . \, C)=1$.

\flushpar
Let
$\overline{E} := \varphi^{-1}(E)$,
$\varphi_{\overline{E}}
:= \varphi|_{\overline{E}} : \overline{E} \to E$,
and
$L := \text{Exc } \varphi$.
Moreover let
$\nu : \widetilde{E} \to E$,
$\overline{\nu} : (\overline{E})\sptilde \to \overline{E}$
be the normalization of $E$, $\overline{E}$,
respectively, and
$(\varphi_{\overline{E}})\sptilde :
(\overline{E})\sptilde \to \widetilde{E}$
the induced morphism.

$$\CD
\Sigma_1 @. \simeq @. (\overline{E})\sptilde
@>\overline{\nu}>> \overline{E} @. \subset @. \overline{U} \\
@. @. @V{(\varphi_{\overline{E}})\sptilde}VV
@V{\varphi_{\overline{E}}}VV @. @VV{\varphi}V \\
\Bbb P^2 @. \simeq @. \widetilde{E} @>>{\nu}> E @. \subset @. U \\
@. @. @. @VVV @. @VV{g_U}V \\
@. @. @. P @. \in @. V
\endCD
\tag 5.5.1
$$

By Lemma 2.7 (1) and Lemma 5.4,

\flushpar
(5.5.2) \qquad \qquad \qquad \qquad
$\widetilde{E} \simeq \Bbb P^2$, \, \, and \, \,
$\nu^*\Cal O_E(-K_X) \simeq \Cal O_{\Bbb P^2}(2)$.

By (3.2.1) again,
$\varphi_{\overline{E}}$
is just the blow-up of $E$ with center $x$.
Since $x$ is a smooth point of $E$,
$(\varphi_{\overline{E}})\sptilde$ coincides with
$\varphi_{\overline{E}}$
in an analytic neighborhood of $x$, and
it follows that $(\varphi_{\overline{E}})\sptilde$
is the blow-up of
$\widetilde{E} \simeq \Bbb P^2$ with center $x$.
In particular

\flushpar
(5.5.3) \qquad \qquad \qquad \qquad \qquad \qquad \quad
$(\overline{E})\sptilde \simeq \Sigma_1$.

\flushpar
(5.5.4) \, \, \,
Let $\widetilde{M}$, $\widetilde{l}$
be the minimal section and a ruling of
$(\overline{E})\sptilde \simeq \Sigma_1$, respectively,
and let $M := \overline{\nu}(\widetilde{M})$,
$l := \overline{\nu}(\widetilde{l})$.
Note that $\varphi(M)=\{x\}.$

\flushpar
We claim that

\flushpar
(5.5.5) \, \, \, $\overline{\nu}^*
\Cal O_{\overline{E}}(-K_{\overline{U}})$ is ample.

In fact,
since $\varphi$ contracts $M$ (5.5.4),

\flushpar
(5.5.6) \qquad \qquad \qquad \, \, \, \, \,
$(\overline{\nu}^*(-K_{\overline{U}}) \, . \, \widetilde{M})
= (-K_{\overline{U}} \, . \, M) =1 >0$

\flushpar
(5.5.0). On the other hand, since
$(\varphi_{\overline{E}})\sptilde(\widetilde{l})$
is a line in
$\widetilde{E} \simeq \Bbb P^2$,
and since $\nu^*\Cal O_E(-K_U)$
$\simeq \Cal O_{\Bbb P^2}(2)$
(5.5.2), we have
$$
(-K_U \, . \, \varphi(l)) = (-K_U \, . \,
(\nu \circ (\varphi_{\overline{E}})\sptilde)(\widetilde{l})) = 2.
$$
Moreover, since
$\varphi(l)$ intersects with
the center $S_x$ of $\varphi$
transversally at the single point $\{x\}$,
we have
$$
(L \, . \, l)=1.
$$
Hence
$$
\align
(5.5.7) \qquad \qquad \qquad \quad \, \,
(\overline{\nu}^*(-K_{\overline{U}}) \, . \, \widetilde{l})
&= (-K_{\overline{U}} \, . \, l) \\
&= (\varphi^*(-K_U)-L \, . \, l) \\
&= (-K_U \, . \, \varphi(l))
- (L \, . \, l) \qquad \qquad \qquad \qquad \qquad \, \, \\
&= 2 - 1 >0.
\endalign
$$

\flushpar
(5.5.6) and (5.5.7) prove (5.5.5).

Finally,
recall that
$g_U|_{U - E} : U \! - \! E \to V \! - \! P$
is a $\Bbb P^1\text{-bundle}$ (Lemma 5.2).
Moreover since $l_E(R)=2$,
any limit conic in $E$ is
generically reduced, and thus
$S_x - \{x\}$ is a subsection of
$g_U|_{U - E}$, namely,
$g|_{S_x - \{x\}} : S_x \! - \! \{x\} \to g(S_x) \! - \! P$
is an isomorphism.
Hence

\flushpar
(5.5.8) \, \, \,
$g_U \circ \varphi$
is a conic bundle over $V \! - \! P$.

\flushpar
In particular,
$-K_{\overline{U}}$ is
$g_U \circ \varphi\text{-ample}$
on $\overline{U} - \overline{E}$.
By combining this and (5.5.5), we conclude that
$-K_{\overline{U}}$ is
$g_U \circ \varphi\text{-ample}$.
\quad \qed
\enddemo

\proclaim{Proposition 5.6}\ \
Under the Assumption 5.0,
$E$ is normal: $E \simeq \Bbb P^2$,
and $V$ is smooth.
Moreover,
(III) of Theorem 5.1 holds.
\endproclaim

\demo{Proof}\ \
We will follow the notations
in the proof of Lemma 5.5.
Since $\rho(\overline{U}/V) =2$,
it follows from Lemma 5.5 that

\flushpar
(5.6.1) \, \, \, $\overline{NE}(\overline{U}/V)$
admits exactly two extremal rays, one of which
defines $\varphi$.

\flushpar
Let $\varphi^+ : \overline{U} \to U^+$ be
the contraction associated to the other extremal ray,
and $g^+ : U^+ \to V$ the structure morphism.
By construction, $\varphi^+$ is a divisorial contraction
such that $g^+(\text{Exc }(\varphi^+))=g(S_x)$ (5.5.8).

Let $E^+ := \varphi^+(\overline{E})$ and
$\varphi_{\overline{E}}^+ := \varphi^+|_{\overline{E}}
: \overline{E} \to E^+$.
Let
$\nu^+ : (E^+)\sptilde \to E^+$
be the normalization of $E^+$, and
$(\varphi_{\overline{E}}^+)\sptilde :
(\overline{E})\sptilde \to (E^+)\sptilde$
the induced morphism.

$$\CD
\Sigma_1 @. \simeq @. (\overline{E})\sptilde
@>\overline{\nu}>> \overline{E} @. \subset @. \overline{U} \\
@. @. @V{(\varphi_{\overline{E}}^+)\sptilde}VV
@V{\varphi_{\overline{E}}^+}VV @. @VV{\varphi^+}V \\
@. @. (E^+)\sptilde @>>{\nu^+}> E^+ @. \subset @. U^+ \\
@. @. @. @VVV @. @VV{g^+}V \\
@. @. @. P @. \in @. V
\endCD
\tag 5.6.2
$$

By the upper-semi-continuity
of the fiber dimension of $\varphi^+$,
$\varphi^+$ contracts at least one curve
in $\overline{E}$.
Let $M$, $l$ be as in (5.5.4),
then $\varphi$ contracts $M$. Thus
$\varphi^+$ never contracts $M$,
and $\varphi^+$ necessarily contracts
$l\text{'s}$.
In particular,

\flushpar
(5.6.3) \, \, \, $E^+$ is an irreducible rational curve, and
$\varphi^+$ has no 2-dimensional fibers.
Thus by the result of Ando ([A] Theorem 2.3),

\flushpar
(5.6.4) \, \, \, Both $U^+$ and
$S^+ := \varphi^+(\text{Exc }(\varphi^+))$
are smooth, and $\varphi^+$ is the blow-up
with center $S^+$.
In particular, all fibers of
$\varphi_{\overline{E}}^+ : \overline{E} \to E^+$
are isomorphic to $\Bbb P^1$.

\flushpar
Let
$\text{Sing }E^+ = \{q_1, \dots , q_r\}$.
Then by (5.6.4)

\flushpar
(5.6.5) \qquad \qquad \qquad \qquad \qquad
$\text{Sing } \overline{E}
= \coprod\limits_{k=1}^r \varphi^{+ \, -1}(q_k)$.

\flushpar
Since $M$ dominates $E^+$ through $\varphi^+$,
and since $\overline{E}$ is smooth along $M$,
it follows that
$\text{Sing }E^+ = \emptyset$.
Hence $\overline{E}$ is also smooth,
by (5.6.5). Namely, we have

\flushpar
(5.6.6) \qquad \qquad \qquad \qquad \quad \, \, \,
$E^+ \simeq \Bbb P^1$,
\, \, and \, \, $\overline{E} \simeq \Sigma_1$

\flushpar
(5.5.3). Note in particular that every fiber of $g^+ : U^+ \to V$
is isomorphic to $\Bbb P^1$.
Since $U^+$ is smooth (5.6.4),

\flushpar
(5.6.7) \, \, \, $(V,P)$ is a smooth germ, and
$g^+$ is a trivial $\Bbb P^1\text{-bundle}$:
$U^+ \simeq V \times \Bbb P^1$

\flushpar
([loc.cit], Theorem 3.1).
Finally, since $\varphi_{\overline{E}}$ is an isomorphism
outside $M = \varphi^{-1}(x)$,
(5.6.6) implies that
$E$ is smooth outside $x$. Moreover $E$ is smooth
also at $x$, by the assumption of Lemma 5.5. Hence
we get the smoothness of the whole $E$, {\it i.e. \/}

\flushpar
(5.6.8) \qquad \qquad \qquad \qquad \qquad \qquad \quad
$E \simeq \Bbb P^2$.

At the same time, the above proof
just gives the elementary
transformation as described in (III) of Theorem 5.1.
\quad \qed
\enddemo

\demo{Determination of the normal bundle $N_{E/X}$}\ \
({\it cf. \/} Andreatta-Wi\'sniewski [AW1])
\enddemo

{}From now on we shall prove (II) of Theorem 5.1.
First we shall prove:

\proclaim{Lemma 5.7}\ \
Under the Assumption 5.0 and
Proposition 5.6,

\flushpar
(1) \, \, \, There is a divisor $L$ in $U$ such that
$2L \sim -K_U$. Moreover $\text{Bs }|L| = \emptyset$.

\flushpar
(2) \, \, \, Take a member $L$ of $|L|$
such that $L \not\supset E$. Then
$L$ is a smooth 3-fold, and
$g|_L : L \to V$
is a divisorial contraction which contracts
a divisor to a curve.
\endproclaim

\demo{Proof}\ \
Since $-K_U$ is $g_U\text{-ample}$,
$H^i(U, \Cal O_U)=0$ for $i \geq 1$ [KaMaMa].
On the other hand, since
$E \simeq \Bbb P^2$ (Proposition 5.6),
$H^i(E, \Cal O_E)=0$ for $i \geq 1$.
Furthermore since $U$ can be retracted to $E$,
$H^2(U, \Bbb Z) \simeq H^2(E, \Bbb Z)$.
Thus

\flushpar
(5.7.1) \, \, \, The natural $\text{Pic }U \to \text{Pic }E$
is an isomorphism.

\flushpar
(5.7.2) \, \, \, Let $L$ be a divisor
on $U$ which corresponds to $\Cal O_{\Bbb P^2}(1)$
on $E$ under (5.7.1).

\flushpar
Since $-K_U$ corresponds to $\Cal O_{\Bbb P^2}(2)$
(Lemma 2.7 (1) with Proposition 5.6),
we have

\flushpar
(5.7.3) \qquad \qquad \qquad \qquad \qquad \qquad \, \, \,
$2L \sim -K_U$.

\flushpar
Then by a version of the non-vanishing theorem
([AW1] Theorem 3.1),

\flushpar
(5.7.4) \qquad \qquad \qquad \qquad \qquad \qquad \, \, \,
$\text{Bs }|L| \not\supset E$.

\flushpar
Choose any $L \in |L|$ such that $L \not\supset E$.
Then $C:=L \cap E$ is a line in $\Bbb P^2$ (5.7.2).
Then

\flushpar
(5.7.5) \, \, \, $L$ is smooth.

\flushpar
Indeed, let $x$ be any point of $C$ and
$B_x$ an analytic neighborhood
of $x$ in $U$. Then in $B_x$,
$E$ is the intersection of two smooth hypersurfaces.
In order to obtain $C$, it is enough to cut out
by one more hypersurface $L$. Since $C$ is still smooth,
$L$ is necessarily smooth in $B_x$.
Thus $L$ is smooth along $C$.
By shrinking $V$ if necessary,
we can make $L$ everywhere smooth, and we have (5.7.5).

Consider $g|_L : L \to V$.
This has connected fibers by $(L \, . \, f)=1$ (5.7.3),
and thus is a birational morphism onto $V$.
Moreover since

\flushpar
(5.7.6) \qquad \qquad \qquad \quad \, \, \, \,
$-2K_L \sim -2K_U-2L |_L \sim -K_U|_L$

\flushpar
(5.7.3), $-K_L$ is $g|_L\text{-ample}$.
Thus by [Mo2],
$g|_L$ is a divisorial contraction, and
$|-K_L|$ is $g|_L\text{-free}$.
Since $-K_L \sim L|_L$,

\flushpar
(5.7.7) \qquad \qquad \qquad \qquad \qquad \qquad
$\text{Bs }\bigm|L|_L\bigm| = \emptyset$.

\flushpar
Finally, by the exact sequence
$$
0 \longrightarrow \Cal O_U \longrightarrow \Cal O_U(L)
\longrightarrow \Cal O_L(L) \longrightarrow 0
$$
with $R^1(g_U)_*\Cal O_U=0$, we conclude that
$$
\text{Bs }|L| = \text{Bs }\bigm|L|_L\bigm| = \emptyset
$$
(5.7.7). \quad \qed
\enddemo

\proclaim{Theorem 5.8}\ \
Under the Assumption 5.0, let
$N_{E/X}$ be the normal bundle
of $E \simeq \Bbb P^2$ (Proposition 5.6)
in $X$. Then
$$
N_{E/X} \simeq \Omega_{\Bbb P^2}^1(1).
$$
\endproclaim

\demo{Proof}\ \
Consider the linear system
$$
\Lambda := \bigl\{ L|_E \in |\Cal O_{\Bbb P^2}(1)| \, \bigm|
\, L \in |L| \bigr\}.
$$
Then by Lemma 5.7 $\Lambda$ is free, and thus
it must coincide with the whole $|\Cal O_{\Bbb P^2}(1)|$.
(Recall that any proper sub-linear system of
$|\Cal O_{\Bbb P^2}(1)|$ has a base point.)
Namely,

\flushpar
(5.8.1) \, \, \,
For any line $C$
in $E \simeq \Bbb P^2$,
there exists $L_C \in |L|$ such that $L_C |_E =C$.

\flushpar
In particular, $L_C$ is smooth,
and $g|_{L_C} : L_C \to V$ is a divisorial
contraction which contracts a divisor to a curve
(Lemma 5.7 (2)).
Note that $C$ is a fiber of $g|_{L_C}$.
Then
$$
N_{E/X} \otimes \Cal O_C \simeq N_{C/{L_C}}
\simeq \Cal O_{\Bbb P^1} \oplus \Cal O_{\Bbb P^1}(-1)
$$
[Mo2]. Since the line $C$ is arbitrarily chosen,
it follows that $N_{E/X}$ is a uniform bundle.
Since $c_1(N_{E/X})=-1$, by Van de Ven's theorem [V],
$N_{E/X}$ is isomorphic either to
$$
\Cal O_{\Bbb P^2} \oplus \Cal O_{\Bbb P^2}(-1)
\, \, \text{ or } \, \,
\Omega_{\Bbb P^2}^1(1).
$$
In the former case, we have
$H^0(N_{E/X}) \not= 0$ and  $H^1(N_{E/X}) = 0$,
which implies that $E$ deforms inside $X$,
a contradiction. Hence
$$
N_{E/X} \simeq \Omega_{\Bbb P^2}^1(1). \quad \qed
$$
\enddemo

\head \S 6.\ Glueing chain of rational curves
\endhead

\definition{Definition 6.0}\ \
Let $B$ be a 1-dimensional reduced proper scheme,
and $B=\bigcup_{i \in I} B_i$ the irreducible decomposition.
Then $B$ is called a {\it tree \/}
of $\Bbb P^1\text{'s}$ if
the following conditions are satisfied:
$$
\cases
(6.0.1)\ \ B_i \simeq \Bbb P^1 \text{ for each } i,\\
(6.0.2)\ \ \text{for any } i,j (i \not= j), \,
B_i \cap B_j \text{ is either }
\emptyset \text{ or one point intersecting transver-}\\
\qquad \quad \, \, \,
\text{sally, and} \\
(6.0.3)\ \ \text{there is no subset } \{i_1, \dots , i_r \}
\subset I \, \, (r \geq 2)
\text{ such that for each } k=1, \dots , \\
\qquad \quad \, \, \,
r-1, \,
B_{i_k} \cap B_{i_{k+1}} \not= \emptyset \text{ and }
B_{i_r} \cap B_{i_1} \not= \emptyset.
\endcases
$$

Let $X$ be an algebraic variety.
A 1-dimensional closed subset
$C$ of $X$
is called a {\it chain of rational curves \/} if
there is a tree $B=\cup B_i$ of $\Bbb P^1\text{'s}$
and a morphism $\varphi_0 : B \to X$
whose image is equal to $C$ such that
$\varphi_0|_{B_i}$ is birational onto its image for each $i$.
\enddefinition

In this section, we shall recall Koll\'ar-Miyaoka-Mori's
method of glueing chain of rational curves
[KoMiMo2,3].
To be more precise, we shall give a sufficient condition
when we can glue a given chain $C$ of rational curves
on $X$ to an irreducible curve $C'$
by a deformation inside $X$.
We will give it only for the simplest case,
namely, the case
that $C'$ is rational, and that
there is no prescribed fixing point
under deformations,
which is suitable for our aim.

Andreatta-Wi\'sniewski [AW2]
essentially used this method
to classify 4-dimensional
divisorial contractions
which contracts divisors
to surfaces.
We are much inspired by their idea.

\definition{Notation 6.1}\ \
Let $X$ be a smooth projective variety and $C \subset X$
a chain of rational curves on $X$.
By Definition 6.0, there is
a tree $B$ of $\Bbb P^1\text{'s}$
and a morphism $\varphi_0 : B \to X$
which is birational onto $C$.

Let $S \to \Bbb A^1$ be
a proper surjective morphism
from a smooth surface $S$
whose fiber at $0$ is isomorphic to $B$, and
another fiber all isomorphic to $\Bbb P^1$.
(Such a morphism can easily be constructed
by a successive blowing-up
from $\Bbb P^1 \times \Bbb A^1$.)

We consider the connected component of the
relative Hilbert scheme
$\text{Hom}_{\Bbb A^1}(S,X \times \Bbb A^1)$
over $\Bbb A^1$ containing $[\varphi_0]$,
with the structure morphism
$$
\lambda : \text{Hom}_{\Bbb A^1}
(S,X \times \Bbb A^1)_{[\varphi_0]} \to \Bbb A^1
$$
([KoMiMo2],(1.1)).
\enddefinition

\proclaim{Theorem 6.2}\,
(Koll\'ar-Miyaoka-Mori. See [KoMiMo2](1.2))\ \

Under the Notation 6.1,

\flushpar
(1)\, \, \, $\dim \text{Hom}_{\Bbb A^1}
(S,X \times \Bbb A^1)_{[\varphi_0]}
\geq \dim X + (-K_X \, . \, C) + 1$.

\flushpar
(2)\, \, \, $\dim \text{Hom}(B, X)_{[\varphi_0]}
\geq \dim X + (-K_X \, . \, C)$.

\flushpar
(3)\, \, \, If the equality holds in (2), then
$C$, as a 1-cycle of $X$, deforms
to an irreducible
rational curve inside $X$.
Namely, there is an
analytic open subset $A \ni 0$ of $\Bbb A^1$, and
a morphism
$\varphi : S \times_{\Bbb A^1} A \to X$
such that
$\varphi|_{S_0} = \varphi_0$.
\endproclaim

\demo{Proof}\ \ (1) is a direct consequence of
[KoMiMo2](1.2). Since $\text{Hom}(B, X)_{[\varphi_0]}$
is the fiber of $\lambda$ at $0$, (2)
follows from (1). Finally, the assumption of (3)
implies that
$\text{Hom}_{\Bbb A^1}
(S,X \times \Bbb A^1)_{[\varphi_0]}$
dominates $\Bbb A^1$ through $\lambda$.
Thus (3) holds.
\quad \qed
\enddemo

\head \S 7.\ Proof of \Rh $\leq 2$
\endhead

The main result of this section is:

\proclaim{Theorem 7.1}\qquad \qquad \qquad \, \,
\Rh$\leq 2$.
\endproclaim

\definition{Assumption 7.2}\ \
To prove Theorem 7.1, we may assume

\flushpar
(7.2.1) \qquad \qquad \qquad \qquad \qquad \qquad \quad
$l_E(R)=1$,

\flushpar
according to Theorem 5.1.
{}From now on we assume

\flushpar
(7.2.2) \qquad \qquad \qquad \qquad \qquad \qquad \, \, \,
$\rho (U/V) \geq 3$

\flushpar
to get a contradiction, till the end of this section.
In particular, we may assume that

\flushpar
(7.2.3) \, \, \,
$g_U|_{U - E} : U \! - \! E \to V \! - \! P$
is not a $\Bbb P^1\text{-bundle}$
(Theorem 3.5), and the assumption
of Theorem 3.7 is satisfied.

In this section, for any curve $l$
which is contained in a fiber of $g_U$,
the symbol $[l]$ means the numerical equivalence
class of $l$ in $U$, not in $X$:
$[l] \in \overline{NE}(U/V)$,
as in the Notation 3.6.
\enddefinition

\proclaim{Lemma 7.3}\ \
Assume $\rho (U/V) \geq 3$, and let $f$ be
a general fiber of $g_U$. Then there exist two
rational curves $l$, $l'$ of length 1 in $E$ which satisfy
$$
l \cap l' \not= \emptyset, \, \, \,
\Bbb R_{\geq 0}[l] \not= \Bbb R_{\geq 0}[l'], \, \text{ and } \,
[f] \not\in \Bbb R_{\geq 0}[l] + \Bbb R_{\geq 0}[l'].
$$
\endproclaim

\demo{Proof}\ \
By assumption $l_E(R)=1$, there exists a couple
$\{l_1, l_2\}$ of rational curves in $E$ with
$l_1 \cap l_2 \not= \emptyset$ and

\flushpar
(7.3.1) \qquad \qquad \qquad \qquad \quad \qquad \qquad \! \!
$f \equiv l_1 + l_2$.

\flushpar
Since $\rho (U/V) \geq 3$ by assumption,
$\overline{NE}(U/V)$ has at least three extremal rays. Thus
by Theorem 3.7 and the connectedness of $E$,
there is another rational curve $l_3$ of length 1 in $E$
such that $l_i \cap l_3 \not= \emptyset$ for $i=1$ or $2$, and

\flushpar
(7.3.2) \qquad \qquad \qquad \qquad \quad \, \,
$\Bbb R_{\geq 0}[l_3] \not= \Bbb R_{\geq 0}[l_1], \,
\Bbb R_{\geq 0}[l_2]$.

\flushpar
Without loss of generality, we may assume

\flushpar
(7.3.3) \qquad \qquad \qquad \qquad \qquad \qquad \, \, \, \,
$l_1 \cap l_3 \not= \emptyset$.

\flushpar
Then
$$
\align
[l_1] + [l_3] &= [f] - [l_2] + [l_3] \qquad (\text{by } (7.3.1)) \\
&\not= [f]. \qquad \qquad \qquad \, \, \, \, (\text{by } (7.3.2))
\endalign
$$
This, together with (7.3.2) and (7.3.3),
proves our lemma.
\quad \qed
\enddemo

\proclaim{Lemma 7.4}\ \
Let $\{l, l'\}$ be as in Lemma 7.3.
Then these are contained in the
same irreducible component, say $E_i$, of $E$.
\endproclaim

\demo{Proof}\ \
Let $B=B_1 \cup B_2$ be the reduced scheme
which is isomorphic to the union of
two distinct lines in $\Bbb P^2$.
Let $C:= l \cup l'$, and consider
any morphism
$\alpha : B \to C \subset X$
such that $\alpha|_{B_1}$, $\alpha|_{B_2}$
gives the normalization
of $l$, $l'$, respectively.

We shall divide our proof into cases:

\flushpar
{\it Case (1) \/} \qquad \qquad \qquad \qquad \quad
$\dim \text{Hom}_{[\alpha]}(B \, , \, X) \geq 7$.

In this case, we shall derive a contradiction
by assuming that $l$ and $l'$ are not contained
in the same irreducible component of $E$.
By $\dim \text{Aut } B = 4$, it follows that
$C$ deforms, as a 1-cycle of $U$,
with an at least 3-dimensional
parameter space. This is possible only when
$$
\cases
l \subset E_i, l' \subset E_j \, \, \, \,
(\exists i \not= \exists j), \\
\widetilde{E}_i \simeq \widetilde{E}_j \simeq \Bbb P^2,
\text{ and} \\
\dim (E_i \cap E_j) =1,
\endcases
$$
in view of Theorem 2.8, where
$\widetilde{E}_i$ is the normalization of $E_i$, etc.
Let $C$ be any irreducible
curve in $E_i \cap E_j$. Then
$$
[l], [l'] \in \Bbb R_{\geq 0}[C],
$$
which contradicts
$\Bbb R_{\geq 0}[l] \not= \Bbb R_{\geq 0}[l']$.

\flushpar
{\it Case (2) \/} \qquad \qquad \qquad \qquad \quad
$\dim \text{Hom}_{[\alpha]}(B \, , \, X) \leq 6$.

In this case, Theorem 6.2 (3) says that
$l \cup l'$ deforms to an irreducible rational
curve of length 2 inside $X$.
If such a deformation goes outside $E$, then it must
coincide with the whole fiber of $g_U$, a contradiction
to $f \not\equiv l + l'$. Thus
$l \cup l'$ deforms inside $E$ to an irreducible
rational curve of length 2, in particular,
$l$ and $l'$ must be contained in the same
irreducible component of $E$.
\quad \qed
\enddemo

\proclaim{Lemma 7.5}\ \ Let $\{l, l'\}$ and $E_i$ be
as in Lemma 7.3 and 7.4, and
let $\varphi$, $\varphi'$ be the contraction
associated to the extremal rays
$\Bbb R_{\geq 0}[l]$ and $\Bbb R_{\geq 0}[l']$,
respectively. Then
$$
\dim \varphi(E_i) = \dim \varphi'(E_i) = 1.
$$
\endproclaim

\demo{Proof}\ \
By assumption,
$$
\cases
\dim \varphi(l) = \dim \varphi'(l') =0, \text{ and} \\
\dim \varphi(l') = \dim \varphi'(l) =1. \\
\endcases
\tag 7.5.1
$$
In particular
$$
\dim \varphi(E_i) \geq 1, \, \, \dim \varphi'(E_i) \geq 1.
$$
First assume $\dim \varphi(E_i) = \dim \varphi'(E_i) = 2$,
and let $\varphi(E_i)\sptilde$,
$\varphi'(E_i)\sptilde$
be the normalization of
$\varphi(E_i)$, $\varphi'(E_i)$,
respectively. Then
$\widetilde{E}_i \to \varphi(E_i)\sptilde$
and
$\widetilde{E}_i \to \varphi'(E_i)\sptilde$
are birational morphisms
which are both not isomorphisms
and contracts distinct curves
(7.5.1).
This is impossible, by Theorem 2.8.

Thus let us assume
$\dim \varphi(E_i)=2$ and
$\dim \varphi'(E_i)=1$, say,
and consider $D := \text{Exc } \varphi = L_U(l)$,
which is a prime divisor (Theorem 3.7 (1)).
In this case, obviously $D \not\supset E_i$.
In particular $(D \, . \, l') >0$,
since $l \cap l' \not= \emptyset$.
On the other hand, take a deformation
$\overline{l'}$ of $l'$ which is not contained in $E$
(Theorem 3.7). Then since
$f \not\equiv l + \overline{l'} \equiv l + l'$
by assumption,
$D \cap \overline{l'} = \emptyset$
and thus $(D \, . \, \overline{l'})=0$,
a contradiction. Hence Lemma 7.5.
\quad \qed
\enddemo

\proclaim{Lemma 7.6}\ \ Under Lemma 7.3, 7.4 and 7.5,
$\varphi$ and $\varphi'$ have no 2-dimensional fibers
which intersect with $E_i$.
\endproclaim

\demo{Proof}\ \
Assume that $\varphi$, say, has a 2-dimensional fiber
$F$ with $F \cap E_i \not= \emptyset$.
Such an $F$ is classified by Andreatta-Wi\'sniewski
[AW2]:

\flushpar
(7.6.1) \, \, \, $F$ is irreducible and is isomorphic
either to $\Bbb P^2$ or a singular quadric surface
\qquad \qquad
$S_2$ in $\Bbb P^3$, and
$\Cal O_F(-K_X) \simeq
\Cal O_{\Bbb P^2}(1), \Cal O_{S_2}(1)$, respectively.

\flushpar
Let

\flushpar
(7.6.2) \qquad \qquad \qquad \qquad \qquad \qquad \, \, \,
$l := F \cap E_i$.

\flushpar
Then $l$ is the fiber of
$\varphi|_{E_i}$ at the point
$\varphi(F)$. By Lemma 7.5, $l$
is purely 1-dimensional in $F$.
{}From this and (7.6.1),
$l$ is connected.
In particular, $\varphi|_{E_i}$
has connected fibers near $\varphi(F)$.
Thus
$l$ is an irreducible rational curve
of length 1, since so is a general fiber
of $\varphi|_{E_i}$ [A].
Again by (7.6.1),

\flushpar
(7.6.3) \, \, \,
$l$ is a line ({\it resp. \/} a ruling) of $F$
when $F \simeq \Bbb P^2$ ({\it resp. \/} $F \simeq S_2$).

Let $l''$ be a general line
({\it resp. \/} a general ruling)
of $F$, and let $l'$ be
the unique fiber of $\varphi'|_{E_i}$
passing through the point $Q := l \cap l''$.
In particular $l' \cap l'' \ni Q$.
Let $C:=l' \cup l''$. Then $C$ deforms,
as a 1-cycle of $U$,
with a 2-dimensional ({\it resp. \/}1-dimensional)
parameter space. So denote by
$B=B_1 \cup B_2$ the reduced scheme
isomorphic to two intersecting lines
in $\Bbb P^2$, and by $\alpha : B \to C \subset X$
a morphism such that $\alpha|_{B_1}$, $\alpha|_{B_2}$
are isomorphisms with $\alpha(B_1 \cap B_2)=Q$.
Then
$$
\dim \text{Hom}_{[\alpha]} (B \, , \, X) \leq 6 \, \,
(\text{{\it resp. \/}} \leq 5).
$$

Thus by Theorem 6.2,

\flushpar
(7.6.4) \, \, \, The case
$F \simeq S_2$ cannot happen,

\flushpar
and for the case $F \simeq \Bbb P^2$,
$C$ deforms to an irreducible rational curve
inside $X$. Moreover, such a deformation
never goes outside $E$, since
$l'+l'' \equiv l'+l \not\equiv f$.
On the other hand,
$l' \subset E_i, \not\subset F$, and
$l'' \not\subset E_i, \subset F$, a contradiction.
Hence the Lemma 7.6.
\quad \qed
\enddemo

\proclaim{Lemma 7.7}\ \
Under Lemma 7.3 and 7.4,
$E_i \simeq \Bbb P^1 \times \Bbb P^1$.
$l$ and $l'$ are two intersecting rulings.
\endproclaim

\demo{Proof}\ \
Let $C:=\varphi(E_i)$ and $C':=\varphi'(E_i)$
in view of Lemma 7.5:
$$
\varphi|_{E_i} : E_i \to C, \, \, \, \,
\varphi'|_{E_i} : E_i \to C'.
$$
Then
$C$ and $C'$ are both rational curves,
since $\varphi(l')=C$ and
$\varphi'(l)=C'$.
By [A] and Lemma 7.6,

\flushpar
(7.7.1) \, \, \,
$\varphi|_{\text{Exc } \varphi}$
and $\varphi'|_{\text{Exc }\varphi'}$
are both $\Bbb P^1\text{-bundles}$
near $E_i$.

\flushpar
We shall prove that
$C \simeq C' \simeq \Bbb P^1$.
Assume $C$, say, is singular:
$\text{Sing }C = \{ P_1, \dots , P_r\}$.
Then by (7.7.1)

\flushpar
(7.7.2) \qquad \qquad \qquad \qquad \qquad
$\text{Sing }E_i = \coprod\limits_{k=1}^r \varphi^{-1}(P_k)$.

Since $\varphi$ and $\varphi'$ are
the contractions of two distinct
extremal rays, $\varphi'(\varphi^{-1}(P_k))=C'$
for each $k$.
On the other hand, again by (7.7.1),
$\varphi'{}^{-1}(C' \! - \! \text{Sing }C')$
is smooth, a contradiction. Hence $C$ must be
smooth, and similarly for $C'$.

Again by (7.7.1), $E_i$ is smooth.
Since $E_i$ has two distinct
$\Bbb P^1\text{-bundle}$ structures:
$$
\varphi|_{E_i} : E_i \to C \simeq \Bbb P^1 \, \, \text{ and } \, \,
\varphi'|_{E_i} : E_i \to C' \simeq \Bbb P^1,
$$
it follows that
$E_i \simeq \Bbb P^1 \times \Bbb P^1$.
\quad \qed
\enddemo

\proclaim{Lemma 7.8}\ \
Let $\{l, l'\}$ and $E_i$ be as in Lemma 7.7.
Then there exists another
irreducible component
$E_j \simeq \Bbb P^1 \times \Bbb P^1$
of $E$, together with a couple of
intersecting rulings $\{m, m'\}$,
such that $f \equiv l + m$.
\endproclaim

\demo{Proof}\ \
Let $[f]$ be the numerical class of
a general fiber of $g_U$, as in Lemma 7.3.

First, even if we replace the original
$\{l, l'\}$ by another couple of intersecting
rulings of $E_i \simeq \Bbb P^1 \times \Bbb P^1$,
the assumption of Lemma 7.3 is preserved.
So we shall make the following
assumption:

\flushpar
(7.8.0) \, \, \,
By abuse of notation, we shall also denote
by $l$ ({\it resp. \/} $l'$) an
{\it arbitrary \/} ruling of
$E_i \simeq \Bbb P^1 \times \Bbb P^1$
which is linearly equivalent to the original
$l$ ({\it resp. \/} $l'$)
given in Lemma 7.3.
Recall that $\Bbb R_{\geq 0}[l] \not= \Bbb R_{\geq 0}[l']$.

Since each $l$ as in (7.8.0) is isomorphic to $\Bbb P^1$,
$l$ deforms, as a 1-cycle of $U$,
with an at least 2-dimensional
parameter space (Lemma 1.3),
while in $E_i$, $l$ deforms
only with 1-dimensional parameters.
Thus those $l\text{'s}$ which are
not contained in any other $E_j$
are actually contained in some limit conics:

\flushpar
(7.8.1)\, \, \,
There is a connected 1-cycle $f_l \subset E$
of length 2 which contains $l$ as its
irreducible component, and $[f_l]=[f]$.

Let us denote the other irreducible component
of $f_l$ by $m$:

\flushpar
(7.8.2) \qquad \qquad \qquad \qquad \qquad \qquad \, \, \,
$f_l =: l + m$.

Note that since
$[f_l] \not\in \Bbb R_{\geq 0}[l]  + \Bbb R_{\geq 0}[l']$
by assumption,

\flushpar
(7.8.3) \qquad \qquad \qquad \qquad \qquad
$[m] \not\in \Bbb R_{\geq 0}[l]  + \Bbb R_{\geq 0}[l']$.

\flushpar
Since every irreducible curve in $E_i$ is
spanned by $[l]$ and $[l']$ (Lemma 7.7),

\flushpar
(7.8.4) \qquad \qquad \qquad \qquad \qquad \qquad \quad \,
$E_i \not\supset m$.

Let $Q \in l \cap m$, and let $l' \subset E_i$
be as in (7.8.0) which passes
through $Q$. Then
$$
\align
[m] + [l'] &= [l] + [m] - [l] + [l'] \\
&= [f] - [l] + [l'] \qquad \, \, \, \, \, (7.8.2) \\
&\not= [f]. \qquad \qquad \qquad \quad (7.8.0)
\endalign
$$
Hence the pair $\{m, l'\}$ satisfies
the assumption of Lemma 7.3, and it follows
from Lemma 7.7 that there exists an irreducible
component $E_j \simeq \Bbb P^1 \times \Bbb P^1$
of $E$ such that $m$, $l'$ are
intersecting rulings of $E_j$.
By (7.8.4), $E_i \not= E_j$,
and Lemma 7.8 is proved.
\quad \qed
\enddemo

Now we come to the proof of Theorem 7.1,
namely, we assume \Rh$\geq 3$ to get a contradiction.

\proclaim{7.9}\ \ Proof of Theorem 7.1.
\endproclaim

Assume \Rh$\geq 3$. Then
(CL) of Theorem 3.1 holds.

Let $E_i \simeq \Bbb P^1 \times \Bbb P^1$ and
$E_j \simeq \Bbb P^1 \times \Bbb P^1$
be as in Lemma 7.8, and
take a general $x \in E_i$,
and a general $y \in E_j$.
Since
$\Cal O_{E_i}(-K_X) \simeq
\Cal O_{E_j}(-K_X) \simeq
\Cal O_{\Bbb P^1 \times \Bbb P^1}(1,1)$,
(CL) of Theorem 3.1 is impossible.
Hence \Rh $\leq 2$.
\quad \qed

\head \S 8.\ The classification of $E$ in the case
$l_E(R)=1$ and $\text{\Rh}=2$
\endhead

In this section, we shall give
a classification of
$E$ in the case
$l_E(R)=1$ and $\text{\Rh}=2$.

\proclaim{Theorem 8.1}\ \
Assume that

\flushpar
\qquad \qquad \qquad \qquad \qquad \, \,
$l_E(R)=1$ and \Rh$=2$.

\flushpar
Then $E$ is one of the followings:

\flushpar
(0) \, \, \, (Mukai-Wi\'sniewski type)

\flushpar
$E=E_1 \cup E_2$, $E_1 \simeq E_2 \simeq \Bbb P^2$,
$E_1 \cap E_2$ is a point, and
$N_{E_1/X} \simeq N_{E_2/X} \simeq \Cal O_{\Bbb P^2}(-1)^{\oplus 2}$.
In particular

\flushpar
\qquad \qquad \qquad \qquad \qquad \quad \, \,
$\Cal O_{E_i}(-K_X) \simeq \Cal O_{\Bbb P^2}(1)$ \, \, $(i=1,2)$.

\flushpar
(1) \, \, \, $E$ is irreducible and is isomorphic to
$\Bbb P^1 \times \Bbb P^1$,

\flushpar
\qquad \qquad \qquad \qquad \qquad \qquad \, \, \,
$\Cal O_E(-K_X) \simeq \Cal O_{\Bbb P^1 \times \Bbb P^1}(1,1)$.

\flushpar
(2) \, \, \, $E$ is irreducible and is isomorphic to
$\Sigma_1$,

\flushpar
\qquad \qquad \qquad \qquad \qquad \qquad \, \, \,
$\Cal O_E(-K_X) \simeq \Cal O_{\Sigma_1}(M + 2l)$

\flushpar
(Notation 0.14).

\flushpar
(3) \, \, \, $E=E_1 \cup E_2$, $E_1 \simeq \Bbb P^2$ and
$E_2 \simeq \Bbb P^1 \times \Bbb P^1$.
$E_1 \cap E_2$ is a line of $E_1$
and is a ruling of $E_2$.

\flushpar
\qquad \qquad \qquad \,
$\Cal O_{E_1}(-K_X) \simeq \Cal O_{\Bbb P^2}(1)$, \, and \,
$\Cal O_{E_2}(-K_X) \simeq \Cal O_{\Bbb P^1 \times \Bbb P^1}(1,1)$.

\flushpar
(4) \, \, \, $E=E_1 \cup E_2$, $E_1 \simeq \Bbb P^2$ and
$E_2 \simeq \Sigma_1$.
$E_1 \cap E_2$ is a line of $E_1$
and is the negative section of $E_2$.

\flushpar
\qquad \qquad \qquad \, \,
$\Cal O_{E_1}(-K_X) \simeq \Cal O_{\Bbb P^2}(1)$, \, and \,
$\Cal O_{E_2}(-K_X) \simeq \Cal O_{\Sigma_1}(M + 2l)$.

In (0), $P$ is an ordinary double singular point of $Y$,
and $g_U|_{U-E} : U \! - \! E \to V \! - \! P$
is a $\Bbb P^1\text{-bundle}$.

In (1) $\sim$ (4), $P$ is a smooth point of $Y$,
and $g_U|_{U-E} : U \! - \! E \to V \! - \! P$ is
a conic bundle with an irreducible discriminant divisor.
\endproclaim

First we shall prove:

\proclaim{Proposition 8.2}\ \
Assume that
\Rh$=2$ and
$g_U|_{U-E} : U \! - \! E \to V \! - \! P$
is a $\Bbb P^1\text{-bundle}$.
Then (0) of Theorem 8.1 holds.
\endproclaim

\demo{Proof}\ \
By Theorem 3.5, the rest things we have to check are
$$
\cases
(8.2.1) \, \, \, \# (E_1 \cap E_2) =1, \, \text{ and} \\
(8.2.2) \, \, \, (V,P) \text{ is an ordinary double point.}
\qquad \qquad \qquad \qquad \qquad \qquad \qquad
\endcases
$$

For (8.2.1): \quad
Assume $\# (E_1 \cap E_2) \geq 2$,
to get a contradiction.
Let $l$ be any line of $E_1 \simeq \Bbb P^2$
which passes through at least two points $\{Q_1, Q_2\}$ of
$E_1 \cap E_2$. Then by $g\text{-freeness}$
of $|-K_X|$, there exists a member
$D \in |-K_U|$ such that $D|_{E_1} = l$.
Such a $D$ is necessarily smooth, by exactly the
same argument as in (5.7.5).
Hence $R^1 (g_U|_D)_* \Cal O_D =0$ [KaMaMa], and
$H^1(\Cal O_{D \cap E})=0$.
On the other hand,
$D \cap E$ has two irreducible components
$D \cap E_1$ and $D \cap E_2$,
which intersect at two distinct points $Q_1, Q_2$
with each other, a contradiction.
Thus we have (8.2.1).

For (8.2.2): \quad
Consider a general member $D \in |-K_U|$.
Then $D \cap E = l_1 \cup l_2$
such that $l_i$ is a line in $E_i$ $(i=1,2)$ and
$l_1 \cap l_2 = \emptyset$.
Moreover

(8.2.3) \qquad \qquad \qquad \quad \, \,
$N_{l_i/D} \simeq N_{E_i/X} \otimes \Cal O_{l_i}
\simeq \Cal O_{\Bbb P^1}(-1)^{\oplus 2}$.

\flushpar
Take the Stein factorization
$h : D \to \widetilde{V}$ of $g_U|_D : D \to V$.
We note that since the fiber of $g_U$ at $P$
has two connected components,
the double cover $\widetilde{V} \to V$
is \'etale.
Moreover by Lemma 4.3, $h$ is a small contraction,
and $\text{Exc }h = l_1 \coprod l_2$.
Thus by (8.2.3), it follows that
$h(l_i)$ is an ordinary double point of
$\widetilde{V}$ $(i=1,2)$.
Since $(V,P) \simeq (\widetilde{V},h(l_i))$
through the \'etale double cover $\widetilde{V} \to V$,
we are done.
\quad \qed
\enddemo

\definition{Assumption 8.3}\ \
In the rest of this section, we assume that

\flushpar
(8.3.1) \, \, \,
\Rh$=2$ \, and \,
$g_U|_{U - E} : U \! - \! E \to V \! - \! P$
is not a $\Bbb P^1\text{-bundle}$.

In particular, the assumption of Theorem 3.7
is satisfied, and (CL) of Theorem 3.1 holds.
\enddefinition

\proclaim{Lemma 8.4}\ \
Under the Assumption 8.3,
let $E_1$ be any irreducible component
of $E$, and let $\varphi$, $\varphi'$ be
the contraction morphisms associated to
the two extremal rays of $\overline{NE}(U/V)$.
Then, up to the permutation of $\{ \varphi, \varphi' \}$,
one of the followings holds:

\flushpar
(a) \, \, \, $(\dim \varphi(E_1) \, , \, \dim \varphi'(E_1))=(2,1)$,
\, \, \, $\widetilde{E}_1 \simeq \Sigma_1$, and

\flushpar
\qquad \qquad \qquad \qquad \qquad \qquad
$\nu^* \Cal O_{E_1}(-K_X) \simeq
\Cal O_{\Sigma_1}(\widetilde{M} + 2\widetilde{l})$,

\flushpar
where $\widetilde{M}$ and $\widetilde{l}$
are the minimal section and
a ruling of $\Sigma_1$, respectively,

\flushpar
(b) \, \, \, $(\dim \varphi(E_1) \, , \, \dim \varphi'(E_1))=(1,1)$,
\, \, \, $\widetilde{E}_1 \simeq \Bbb P^1 \times \Bbb P^1$, and

\flushpar
\qquad \qquad \qquad \qquad \qquad \qquad
$\nu^* \Cal O_{E_1}(-K_X) \simeq
\Cal O_{\Bbb P^1 \times \Bbb P^1}(1,1)$, \, \, or

\flushpar
(c) \, \, \, $(\dim \varphi(E_1) \, , \, \dim \varphi'(E_1))=(2,0)$,
$E_1$ is normal, is isomorphic to $\Bbb P^2$ or $S_2$,
and

\flushpar
\qquad \qquad \qquad \qquad \qquad \qquad
$\Cal O_{E_1}(-K_X) \simeq \Cal O_{\Bbb P^2}(1), \, \Cal O_{S_2}(1)$,

\flushpar
respectively.
\endproclaim

\demo{Proof}\ \
First of all, we shall prove

\flushpar
(8.4.1) \qquad \qquad \, \, \,
$(\dim \varphi(E_1) \, , \, \dim \varphi'(E_1))=(2,1)$,
$(1,1)$, or $(2,0)$.

In fact, $(\dim \varphi(E_1) \, , \, \dim \varphi'(E_1))=(2,2)$
is impossible by Theorem 2.8,
as exactly the same argument
as in Lemma 7.5. Next if $\dim \varphi'(E_1)=0$,
then $\varphi$ never contracts any curves, in particular
$\dim \varphi(E_1)=2$. Hence (8.4.1) holds.

\flushpar
{\it Case \/} (a), (b) \qquad \quad \, \,
$(\dim \varphi(E_1) \, , \, \dim \varphi'(E_1))=(2,1)$, \, $(1,1)$.

\flushpar
In these cases $\rho(E_1) \geq 2$, and
it follows from Theorem 2.8 that
$\widetilde{E}_1$ is isomorphic to a Hirzebruch surface:
$$
\widetilde{E}_1 \simeq \Sigma_m.
$$
Moreover since

\flushpar
(8.4.2) \, \, \,
$\nu^*\Cal O_{E_1}(-K_X)$ is ample on $\widetilde{E}_1$,

\flushpar
any section of
$\widetilde{E}_1 \simeq \Sigma_m \to \Bbb P^1$,
except for the minimal one, has length
at least $m+1$. Thus by Theorem 3.1,
it necessarily follows that $m \leq 1$, {\it i.e. \/}

\flushpar
(8.4.3) \qquad \qquad \qquad \qquad \qquad \, \, \,
$\widetilde{E}_1 \simeq \Sigma_1$ \, or \,
$\Bbb P^1 \times \Bbb P^1$.

Assume $\widetilde{E}_1 \simeq \Sigma_1$.

\flushpar
(8.4.4) \, \, \,
Let $\widetilde{M}$, $\widetilde{l}$ be the minimal section
and a ruling of $\widetilde{E}_1 \simeq \Sigma_1$,
respectively, and let
$M := \nu(\widetilde{M})$, $l := \nu(\widetilde{l})$.

\flushpar
Then again by Theorem 3.1
$$
(-K_X \, . \, M) = (-K_X \, . \, l) = 1,
$$
in particular,
$$
-\nu^* K_X \sim \widetilde{M} + 2 \widetilde{l}.
$$
Obviously
$(\dim \varphi(E_1) \, , \, \dim \varphi'(E_1)) = (2,1)$.

Next, assume $\widetilde{E}_1 \simeq \Bbb P^1 \times \Bbb P^1$.

\flushpar
(8.4.5) \, \, \,
Let $\{ \widetilde{l}, \widetilde{l'} \}$
be a couple of intersecting
rulings in $\widetilde{E}_1$,
and let
$l := \nu(\widetilde{l})$, $l' := \nu(\widetilde{l'})$,
such that
$\varphi'$ contracts $\nu(l')$.

\flushpar
Then by Theorem 3.1,
$$
(-K_X \, . \, l)=(-K_X \, . \, l')=1,
$$
in particular
$$
- \nu^* K_X \sim \widetilde{l} + \widetilde{l'}.
$$
We have two possibilities:
$$
(\dim \varphi(E_1) \, , \, \dim \varphi'(E_1))
= (2,1), \, \text{ or } \, (1,1).
$$
Since by Theorem 3.7, $\Bbb R_{\geq 0}[l]$
is the extremal ray of $\overline{NE}(U/V)$
defining $\varphi$,
the former one cannot happen, and we have
$(\dim \varphi(E_1) \, , \, \dim \varphi'(E_1)) = (1,1)$.
Hence (a) and (b).

\flushpar
{\it Case \/} (c) \qquad \qquad \qquad \quad
$(\dim \varphi(E_1) \, , \, \dim \varphi'(E_1))=(2,0)$.

In this case, $E_1$ is normal and is isomorphic
either to $\Bbb P^2$ or $S_2$, as a direct
consequence of Andreatta-Wi\'sniewski [AW2].
\quad \qed
\enddemo

\proclaim{Proposition 8.5}\ \
In Lemma 8.4 (a), $E_1$ is normal: $E_1 \simeq \Sigma_1$.
$Y$ is smooth at $P$.

Moreover if $E$ is reducible,
then $E$ has exactly two irreducible components.
Let $E_2$ be another one.
Then $E_2 \simeq \Bbb P^2$, and
$E_1 \cap E_2$ is the minimal section of
$E_1$ and is a line of $E_2$.
\endproclaim

\demo{Proof}\ \
Let $\varphi : U \to U_{\varphi}$,
$\varphi' : U \to U_{\varphi'}$
be as in Lemma 8.4 such that

\flushpar
(8.5.0) \qquad \qquad \qquad \qquad
$(\dim \varphi(E_1) \, , \, \dim \varphi'(E_1))=(2,1)$.

\flushpar
Let $\nu : \widetilde{E}_1 \to E_1$ be the normalization
of $E_1$, and let
$\widetilde{M}$, $\widetilde{l} \subset \widetilde{E}_1$,
$M$, $l \subset E_1$ be
as in (8.4.4).
Recall that

\flushpar
(8.5.1) \qquad \qquad \quad
$\widetilde{E}_1 \simeq \Sigma_1$, \, \, \, and \, \, \,
$\nu^* \Cal O_{E_1}(-K_X) \simeq
\Cal O_{\Sigma_1}(\widetilde{M} + 2 \widetilde{l})$.

\flushpar
We first claim that

\flushpar
(8.5.2) \, \, \, $\varphi' : U \to U_{\varphi'}$
has no 2-dimensional
fibers which intersect with $E_1$.

Actually, assume that $\varphi'$ has a fiber $F$
with $\dim F = 2$ and $m:= F \cap E_1 \not= \emptyset$.
Then $m$ is a fiber of
$\varphi'|_{E_1} : E_1 \to \varphi'(E_1)$.
In particular, a general fiber of
$\varphi'|_{E_1}$ is disjoint from
$m$. Then by (8.5.1),
a general point of $E_1$ and that of $F$ never
can be joined by a limit conic, which contradicts
Theorem 3.1. Hence $\varphi'$ has
no 2-dimensional fibers, and
(8.5.2) is proved.

Let $h' : U_{\varphi'} \to V$ be the
structure morphism. By (8.5.2) and [A],

\flushpar
(8.5.3) \, \, \, $U_{\varphi'}$ is smooth,
$\varphi'$ is the blow-up of
a smooth codimension 2 center, and
$\varphi'|_{\text{Exc }\varphi'} :
\text{Exc }\varphi' \to \varphi'(\text{Exc }\varphi')$
is a $\Bbb P^1\text{-bundle}$.

We divide our proof into cases.

\flushpar
{\it Case \/} (1) \, \, \, $E$ is irreducible: $E=E_1$.

In this case, since $\dim \varphi'(E)=1$
(8.5.0), $h'$ has no 2-dimensional fibers.
Since $U_{\varphi'}$ is smooth (8.5.3),
$h'$ is a $\Bbb P^1\text{-bundle}$ [loc.cit],
in particular
$V$ is smooth, and
$\varphi'(E) \simeq \Bbb P^1$.
Moreover since $\varphi'|_E : E \to \varphi'(E)$
is a $\Bbb P^1\text{-bundle}$ (8.5.3),
$E$ is smooth: $E \simeq \Sigma_1$.
Thus we get the proposition in the Case (1).

\flushpar
{\it Case \/} (2) \, \, \, $E$ is reducible.

First $\varphi(M)$ and $\varphi'(l)$
are points. Let $E_2$ be any other irreducible
components of $E$. Then by Theorem 3.1 and Lemma 8.4,

\flushpar
(8.5.4) \qquad \qquad \qquad \qquad \qquad \qquad \quad \,
$E_2 \simeq \Bbb P^2$.

\flushpar
Since $\varphi'$ does not contract $E_2$ (8.5.2),
$\varphi$ contracts $E_i$, and hence

\flushpar
(8.5.5) \qquad \qquad \qquad \qquad \qquad \qquad
$E_1 \cap E_2 \supset M$.

\flushpar
Since $M$ is a rational curve of length 1,
$M$ is a line in $E_2$ in $\Bbb P^2$
and in particular
$M \simeq \Bbb P^1$. Hence by (8.5.3),
$E_1$ is smooth:

\flushpar
(8.5.6) \qquad \qquad \qquad \qquad \qquad \qquad \quad
$E_1 \simeq \Sigma_1$.

\flushpar
Next we claim that $E_2$ is unique:

\flushpar
(8.5.7) \qquad \qquad \qquad \qquad \qquad \qquad \, \,
$E = E_1 \cup E_2$.

In fact, take a general $D \in |-K_U|$.
Since $|-K_U|$ is free (Theorem 4.1),
$D$ is smooth, and hence
$R^1 (g|_D)_* \Cal O_D =0$ [KaMaMa].
In particular

\flushpar
(8.5.8) \qquad \qquad \qquad \qquad \qquad \quad \, \, \, \,
$H^1 (\Cal O_{D \cap E}) =0$.

\flushpar
If there are at least two irreducible
components $E_2$, $E_3$ of $E$
which both meet $E_1$,
then by (8.5.5) $D \cap E$ contains
at least three irreducible components
meeting at the point $D \cap M$,
which contradicts (8.5.8).
Hence (8.5.7) is proved.

Consider $\varphi' : U \to U_{\varphi'}$.
Let
$$
\overline{E} := \varphi'(E) = \varphi'(E_2).
$$
Then

\flushpar
(8.5.9) \, \, \, $\varphi'|_{E_2} :
\Bbb P^2 \simeq E_2 \to \overline{E}$
is a finite birational morphism, namely,
is the normalization morphism.

\flushpar
Let $L := \text{Exc } \varphi'$;

\flushpar
(8.5.10) \qquad \qquad \qquad \qquad \quad \, \, \, \,
$-K_U = -{\varphi'}^* K_{U_{\varphi'}} - L$.

\flushpar
Then for a general line $l$ in $E_2 \simeq \Bbb P^2$,
$(L \, . \, l)>0$ and hence

\flushpar
(8.5.11) \qquad \qquad \qquad \,
$(-K_{U_{\varphi'}} \, . \, \varphi'(l))
= (-K_U \, . \, l) + (L \, . \, l) \geq 2$

\flushpar
(8.5.10). In particular, $h' : U_{\varphi'} \to V$
satisfies the following:
$$
\cases
{h'}^{-1}(P) = \overline{E}, \, \text{ which is
irreducible and whose normalziation is } \, \Bbb P^2, \\
-K_{U_{\varphi'}} \, \text{ is } \, h'\text{-ample, and} \\
\varphi'|_{U_{\varphi'} - \overline{E}} :
U_{\varphi'} \! - \! \overline{E} \to V \! - \! P \,
\text{ is a } \, \Bbb P^1\text{-bundle}.
\endcases
\tag 8.5.12
$$
Based on the facts (8.5.3), (8.5.11) and (8.5.12),
we can proceed exactly the same argument
as in the proof of Lemma 5.5 and Proposition 5.6
to deduce that

\flushpar
(8.5.13) \, \, \,
$\overline{E}$ is smooth:
$\overline{E} \simeq \Bbb P^2$,
and $V$ is smooth.

\flushpar
Then by (8.5.9), $\varphi'|_{E_2}$
is an isomorphism. In particular

\flushpar
(8.5.14) \qquad \qquad \qquad \qquad \qquad \quad \, \, \, \,
$E_1 \cap E_2 = M$.

(8.5.4), (8.5.6), (8.5.7), (8.5.13) and (8.5.14)
prove our proposition also in the Case (2).
\quad \qed
\enddemo

\proclaim{Proposition 8.6}\ \
In Lemma 8.4 (b), $E_1$ is normal:
$E_1 \simeq \Bbb P^1 \times \Bbb P^1$.
$Y$ is smooth at $P$.

Moreover if $E$ is reducible,
then $E$ has exactly two irreducible components.
Let $E_2$ be another one.
Then $E_2 \simeq \Bbb P^2$, and
$E_1 \cap E_2$ is a ruling of
$E_1$ and is a line of $E_2$.
\endproclaim

\demo{Proof}\ \
Let $\varphi : U \to U_{\varphi}$,
$\varphi' : U \to U_{\varphi'}$
be as in Lemma 8.4 such that

\flushpar
(8.6.0) \qquad \qquad \qquad \qquad
$(\dim \varphi(E_1) \, , \, \dim \varphi'(E_1))=(1,1)$.

\flushpar
Recall that

\flushpar
(8.6.1) \qquad \qquad \, \,
$\widetilde{E}_1 \simeq \Bbb P^1 \times \Bbb P^1$, \, and \,
$\nu^*\Cal O_{E_1}(-K_X) \simeq
\Cal O_{\Bbb P^1 \times \Bbb P^1} (1,1)$.

\flushpar
We divide our proof into cases.

\flushpar
{\it Case \/} (1) \, \, \, $E$ is irreducible: $E=E_1$.

In this case, consider
$\varphi|_E : E \to \varphi(E)$ and
$\varphi'|_E : E \to \varphi'(E)$.
$\varphi(E)$, $\varphi'(E)$ are both irreducible curves
(8.6.0). Since $\varphi$ and $\varphi'$ has no
2-dimensional fibers,
the smoothness of $E$
follows exactly from the same argument as in Lemma 7.7.
Moreover $U_{\varphi}$ is smooth [loc.cit], and hence
the structure morphism $h : U_{\varphi} \to V$
is a $\Bbb P^1\text{-bundle}$ [loc.cit],
in particular $V$ is smooth.
Thus we are done in the Case (1).

\flushpar
{\it Case \/} (2) \, \, \, $E$ is reducible.

In this case, as in (8.5.4),

\flushpar
(8.6.2) \, \, \,
Any other irreducible component of $E$ is
isomorphic to $\Bbb P^2$, and is a 2-dimensional
fiber of $\varphi$ or $\varphi'$.

\flushpar
We claim that

\flushpar
(8.6.3) \, \, \, Either $\varphi$ or $\varphi'$
has no 2-dimensional fibers.

{}From now on, we assume that
$\varphi$ and $\varphi'$ both have
2-dimensional fibers, to get a contradiction,
till (8.6.12).

Let $E_2$ ({\it resp. \/} $E_3$) be an arbitrary
2-dimensional fiber of $\varphi$
({\it resp. \/} $\varphi'$). Note that

\flushpar
(8.6.4) \qquad \qquad \qquad \qquad \qquad \qquad \,
$E_2 \simeq E_3 \simeq \Bbb P^2$

\flushpar
(8.6.2), and

\flushpar
(8.6.5) \, \, \,
$E_1 \cap E_2$ ({\it resp. \/} $E_1 \cap E_3$) is
the fiber of $\varphi|_{E_1}$ ({\it resp. \/} $\varphi'|_{E_1}$)
at the point $\varphi(E_2)$ ({\it resp. \/} $\varphi(E_3)$).

\flushpar
As in (7.6.3),

\flushpar
(8.6.6) \, \, \,
$l := E_1 \cap E_2$ ({\it resp. \/} $l' := E_1 \cap E_3$)
is a line in $E_2$ ({\it resp. \/} $E_3$).

Choose $Q \in E_1 \cap E_2 \cap E_3 = l \cap l'$
arbitrarily, and consider the linear system
$$
\Lambda_Q := \bigl\{ D \in |-K_U| \big| D \ni Q \bigr\}.
$$
We claim that

\flushpar
(8.6.7) \, \, \,
There exists a point $Q' \in E_1$ such that
$\text{Bs }\Lambda_Q = \{Q, Q'\}$.

In fact, since $\text{Bs }|-K_U|= \emptyset$
(Theorem 4.1), for any line $m$ in
$E_2$ ({\it resp. \/} $E_3$), there exists
$D \in |-K_U|$ such that
$D|_{E_2} = m$ ({\it resp. \/} $D|_{E_3} = m$),
as in (5.8.1).
Hence

\flushpar
(8.6.8) \qquad \qquad \qquad \qquad \qquad \qquad \, \,
$\text{Bs }\Lambda_Q \subset E_1$.

\flushpar
In particular, for a general $D \in \Lambda_Q$,
$$
D \not\supset l, \, \, \, \, D \not\supset l'.
$$
Hence if we let
$$
\nu^* \Lambda_Q :=
\bigl\{ \nu^*(D|_{E_1}) \in | \Cal O_{\Bbb P^1 \times \Bbb P^1}(1,1) |
\big| D \in \Lambda_Q \bigr\},
$$
then

\flushpar
(8.6.9) \, \, \,
$\nu^* \Lambda_Q$ contains an irreducible member.

\flushpar
Since
$$
\nu^*|-K_U| := \bigl\{ \nu^*(D|_{E_1}) \in
| \Cal O_{\Bbb P^1 \times \Bbb P^1}(1,1) |
\big| D \in |-K_U| \bigr\}
$$
is a free linear sub-system of
$| \Cal O_{\Bbb P^1 \times \Bbb P^1}(1,1) |$
(Theorem 4.1), we have
$\dim \nu^*|-K_U| \geq 2$, and hence

\flushpar
(8.6.10) \qquad \qquad \qquad \qquad \qquad \quad \, \, \,
$\dim \nu^* \Lambda_Q \geq 1$.

{}From (8.6.9) and (8.6.10),
$$
\# \text{Bs } \nu^*\Lambda_Q
= (\nu^* (D|_{E_1}))_{\widetilde{E}_1}^2 =2.
$$
This, together with (8.6.8), proves (8.6.7).

Next we claim that

\flushpar
(8.6.11) \, \, \,
A general member $D$ of $\Lambda_Q$ is smooth.

By (8.6.7) and Bertini's theorem,
it is sufficient to prove that
a general $D \in \Lambda_Q$ is smooth at $Q$ and $Q'$.
Since $D|_{E_2}$ is a line in $E_2 \simeq \Bbb P^2$,
$D$ is smooth at $Q$, as the same reason
as in (5.7.5).
Similarly if $Q'$ is contained in
a 2-dimensional fiber of $\varphi$ or $\varphi'$,
then $D$ is smooth at $Q'$.
So assume that $Q'$ is not contained in
any 2-dimensional fiber of $\varphi$ and $\varphi'$.
Then
$$
\varphi|_{E_1} : E_1 \to \varphi(E_1) \, \, \text{ and } \, \,
\varphi'|_{E_1} : E_1 \to \varphi'(E_1)
$$
are both $\Bbb P^1\text{-bundles}$
in a neighborhood of $\varphi(Q')$ and
$\varphi'(Q')$, respectively [A].
Hence $E_1$ is snooth at $Q'$, and
the normalization morphism $\nu : \widetilde{E}_1 \to E_1$
is an isomorphism at $Q'$.
Then by (8.6.9), $D|_{E_1}$ is smooth at $Q'$.
Then by the same reason as above, $D$ is smooth at $Q'$.
Hence (8.6.11).

\flushpar
(8.6.12) \, \, \,
Next we shall prove (8.6.3).

\flushpar
By (8.6.11), $H^1(\Cal O_{D \cap E}) =0$
[KaMaMa].
On the other hand, $D \cap E$
contains three distinct irreducible
rational curves $D \cap E_i$ $(i=1,2,3)$,
meeting at $Q$, a contradiction.
Thus (8.6.3) is proved.

Assume that $\varphi'$, say, has no 2-dimensional fibers.
Then $\varphi'|_{E_1} : E_1 \to \varphi'(E_1)$
is a $\Bbb P^1\text{-bundle}$ [loc.cit],
and $l \simeq \Bbb P^1$ (8.6.6) is a section.
Hence $E_1$ is smooth:

\flushpar
(8.6.13) \qquad \qquad \qquad \qquad \qquad \quad \, \, \, \,
$E_1 \simeq \Bbb P^1 \times \Bbb P^1$.

By combining the above,
$$
\cases
E = \bigcup\limits_{i=1}^n E_i, \\
E_1 \simeq \Bbb P^1 \times \Bbb P^1, \, \,
E_2 \simeq \dots \simeq E_n \simeq \Bbb P^2, \\
E_1 \cap E_i \, \text{ is a ruling in } \, E_1 \,
\text{ and a line in } \, E_i \,
(i=2, \dots n), \, \text{ and} \\
E_1 \cap E_i \sim E_1 \cap E_j \,
\text{ in } \, E_1.
\endcases
\tag 8.6.14
$$

The rest we have to do
is to prove $E = E_1 \cup E_2$, and the
smoothness of $V$.

Assume that $E$ has more than
two irreducible components,
and let $E_2 \simeq E_3 \simeq \Bbb P^2$
be as in (8.6.14). First
$E_1 \cap E_2 = E_1 \cap E_3$
is impossible, exactly by the same reason
as in (8.5.7).
So $E_1 \cap E_2$ and $E_1 \cap E_3$
are two disjoint rulings in $E_1$.
Then $E_2 \cap E_3 \not= \emptyset$
by Theorem 3.1.
Hence

\flushpar
(8.6.15) \qquad \qquad \qquad \, \, \, \, \,
$E_2 \cap E_3 \not= \emptyset$ \, and \,
$E_1 \cap E_2 \cap E_3 = \emptyset$.

\flushpar
Since $E_2 \simeq E_3 \simeq \Bbb P^2$,

\flushpar
(8.6.16) \qquad \qquad \qquad \qquad \qquad \, \, \, \,
$\dim (E_2 \cap E_3) =0$.

\flushpar
Choose $Q \in E_2 \cap E_3$, and general lines
$l_Q \subset E_2$, $l_Q' \subset E_3$
which both pass through $Q$.
Let
$B := l_Q \cup l_Q'$ regarding as a reduced scheme, and
consider
the closed immersion $i : B \hookrightarrow X$.
Since
$l_Q \equiv l_Q'$ (8.6.14),
$l_Q + l_Q' \not\equiv f$
(a general fiber of $g_U$), and hence
$l_Q \cup l_Q'$ never deform outside $E$ (Theorem 3.7).
Thus
$$
\dim \text{Hom }_{[i]} (B, X) =6.
$$
By Theorem 6.2, $l_Q \cup l_Q'$ deforms
to an irreducible curve. This is absurd,
since $l_Q \subset E_2$, $\not\subset E_3$,
and $l_Q' \not\subset E_2$, $\subset E_3$.

Hence $E$ has exactly two irreducible
components, namely,
$E= E_1 \cup E_2$ in (8.6.14).
Assume $\varphi'$, say, satisfies
$\dim \varphi'(E_2)=2$.
Then $\overline{E} := \varphi'(E)$
is an irreducible surface, and hence
the structure morphism $h' : U_{\varphi'} \to V$
satisfies the same condition as in (8.5.12).
Thus as in (8.5.13), we conclude that
$V$ is smooth.
\quad \qed
\enddemo

\proclaim{Proposition 8.7}\ \
Under the same assumption as in Lemma 8.4,
assume that $E$ has at least one
irreducible component $E_1$ which is isomorphic
either to $\Bbb P^2$ or $S_2$ (Lemma 8.4 (c)).
Then $E$ has another irreducible component
which is isomorphic either to $\Sigma_1$ or
$\Bbb P^1 \times \Bbb P^1$.
\endproclaim

\demo{Proof}
Since $\rho(U/V)=2$ while $\rho(E_1)=1$,
$E$ necessarily admits another irreducible component.
Now we assume furthermore that
all the irreducible components of $E$
are isomorphic either to $\Bbb P^2$ or
$S_2$, to get a contradiction.

We put the equivalence relation into
$\{ E_1, \dots , E_n \}$, the whole set of
irreducible components of $E$,
as follows:

\flushpar
(8.7.1) \, \, \,
$E_i \sim E_j$ if and only if
there exists a subset
$\{ E_{i,0}=E_i, E_{i,1}, \dots , E_{i,r}=E_j \}$
of $\{E_1, \dots , E_n\}$
such that for each $k=0, \dots , r-1$,
$E_{i,k} \cap E_{i,k+1}$ contains a rational curve
of length 1.

\flushpar
Then by the same reason as above,
this relation is disconnected, {\it i.e. \/}

\flushpar
(8.7.2) \, \, \, There is a partition
$\{1, \dots , n\} = I \coprod J$ $(I, J \not= \emptyset)$
such that
$E_i \sim E_{i'}$ for each $i, i' \in I$, and
$E_i \not\sim E_j$ for each
$i \in I$ and $j \in J$.

Choose $i \in I$ arbitrarily, and let $l$
be any line ({\it resp. \/} any ruling)
of $E_i$ when $E_i \simeq \Bbb P^2$
({\it resp. \/} $E_i \simeq S_2$).
Consider the deformation locus $L_U(l)$.
By Theorem 3.7, this is a prime divisor in $U$.
Moreover
$$
\cases
L_U(l) \supset \bigcup\limits_{i \in I} E_i, \, \text{ and} \\
L_U(l) \not\supset E_j \, \, (\forall j \in J)
\endcases
\tag 8.7.3
$$
by the above construction.
By the connectedness of $E$, we find
some $j \in J$ such that
$L_U(l) \cap E_j \not= \emptyset$.
Since $U$ is a smooth 4-fold and $\dim L_U(l) =3$,
$\dim E_j = 2$,

\flushpar
(8.7.4) \qquad \qquad \qquad \qquad \qquad \, \, \, \, \,
$\dim (L_U(l) \cap E_j) =1$

\flushpar
for such $j$. Let
$C$ be any
irreducible component of $L_U(l) \cap E_j$,
and let
$l_j$ be
a general line or a general ruling in $E_j$.
Then
$[l], [l_j] \in \Bbb R_{\geq 0}[C]$,
and hence

\flushpar
(8.7.5) \qquad \qquad \qquad \qquad \qquad \qquad \, \, \, \,
$l \equiv l_j$ in $U$.

\flushpar
Since $l_j \not\subset L_U(l)$,

\flushpar
(8.7.6) \qquad \qquad \qquad \qquad \qquad \quad \, \, \, \, \,
$(L_U(l) \, . \, l_j) \geq 0$.

\flushpar
On the other hand,

\flushpar
(8.7.7) \qquad \qquad \qquad \qquad \qquad \quad \, \, \, \, \,
$(L_U(l) \, . \, l) =-1$

\flushpar
(Theorem 3.7). (8.7.6) and (8.7.7) contradict (8.7.5).

Hence $E$ admits
at least one irreducible component which
is isomorphic neither to $\Bbb P^2$ nor $S_2$.
By Lemma 8.4, Proposition 8.5 and 8.6,
such a component must be
isomorphic either to $\Sigma_1$ or
$\Bbb P^1 \times \Bbb P^1$.
\quad \qed
\enddemo

\proclaim{8.8} \ \ Proof of Theorem 8.1.
\endproclaim

Now our Theorem 8.1 is an easy consequence of
Lemma 8.4, Proposition 8.5, 8.6, and 8.7.
\quad \qed

\head \S 9.\ Deformation loci of extremal
rational curves of length 1 (the case \Rh$=1$)
\endhead

\definition{Assumption 9.0}\ \
Throughout this and the next sections,
we assume

\flushpar
(9.0.1) \qquad \qquad \qquad \quad
$l_E(R)=1$ \, and \, \Rh$=1$.
\enddefinition

In this and the next sections, we shall
investigate the structure of $E$
under the above assumption.
In this section, we shall prove
the following:

\proclaim{Theorem 9.1}\ \
Let $g : X \to Y$, $R$, $E$ and $P$
be as above. Assume
$$
l_E(R)=1 \, \text{ and } \, \text{\Rh}=1.
$$
Then for any rational curve $l$
of length 1 in $E$,
the whole deformation locus
$L_X(l)$ of $l$ in $X$ (Definition 1.2)
is purely codimension 1 in $X$.
Each irreducible component of $L_X(l)$ contains $E$,
whenever it meets $E$.
\endproclaim

\definition{Assumption 9.2}\ \
According to Lemma 2.7, it is sufficient
to do the following in order to prove
Theorem 9.1:

\flushpar
(9.2.1) \, \, \,
We assume that
$g|_{U \! - \! E} : U \! - \! E \to V \! - \! P$
is a $\Bbb P^1\text{-bundle}$,
and shall derive a contradiction.
\enddefinition

\proclaim{Lemma 9.3}\ \
Under the Assumption 9.0 and 9.2,
$$
\bigl( E_i \, , \, \Cal O_{E_i}(-K_X) \bigr)
\simeq \bigl( \Bbb P^2 \, , \, \Cal O_{\Bbb P^2}(1) \bigr)
\text{ for any irreducible component } E_i \text{ of } E.
$$
\endproclaim

\demo{Proof} \ \
For any irreducible component
$E_i$ of $E$, we have $l_{E_i}(R)=1$ by Theorem 5.1.
Then
$$
\bigl( \widetilde{E}_i \, , \, \nu_i^* \Cal O_{E_i}(-K_X) \bigr)
\simeq \bigl( \Bbb P^2 \, , \, \Cal O_{\Bbb P^2}(1) \bigr)
\text{ for each } i,
$$
by Assumption 9.2 and Lemma 2.7 (2).

We shall prove that each $E_i$ is normal,
{\it i.e. \/} $E_i$ itself is isomorphic to $\Bbb P^2$.
First of all, consider the rational map associated to
$\big| -K_X|_{E_i} \big|$ :
$$
\Phi : E_i \to \Bbb P^N.
$$
Note that since
$\big| -K_X|_{E_i} \big|$
is free (Theorem 4.1), the above $\Phi$ is actually
a morphism.

We claim that $N= \dim \big| -K_X|_{E_i} \big| =2$.
In fact, there is a natural injective homomorphism
$$
H^0 \bigl( E_i \, , \, \Cal O_{E_i}(-K_X) \bigr)
\hookrightarrow
H^0 \bigl( \widetilde{E}_i \, , \, \nu_i^* \Cal O_{E_i}(-K_X) \bigr)
\simeq
H^0 \bigl( \Bbb P^2 \, , \, \Cal O_{\Bbb P^2}(1) \bigr).
$$
The last term has dimension 3, whereas for the first term,
$$
\dim H^0 \bigl( E_i \, , \, \Cal O_{E_i}(-K_X) \bigr) \geq 3,
$$
since the linear system
$\big| -K_X|_{E_i} \big|$
on a surface $E_i$
is free and ample. Thus the above injection
must be an isomorphism and we deduce
$N=2$.
Then we consider the composition;
$$
\Bbb P^2 \simeq \widetilde{E}_i
\overset \nu_i\to\longrightarrow E_i
\overset \Phi\to\longrightarrow \Bbb P^2
$$
which is exactly
the rational map associated to the linear system
$\big| \nu_i^* (-K_X |_{E_i}) \big|
= \big| \Cal O_{\Bbb P^2}(1) \big|$,
namely, an isomorphism $\Bbb P^2 \to \Bbb P^2$.
Thus the morphism $\Phi : E_i \to \Bbb P^2$
is set theoretically a bijection.
By Zariski Main Theorem, we conclude that
$\Phi$ is actually an isomorphism:
$$
E_i \simeq \Bbb P^2. \quad \qed
$$
\enddemo

{}From now on, we shall show that the normal
bundle $N_{E_i/X}$ is isomorphic to
$\Cal O_{\Bbb P^2}(-1)^{\oplus 2}$ for each $E_i$ (Proposition 9.6),
following the argument of Kawamata ([Ka4]
(2.4)).

\definition{Notation 9.4}\ \
Let $E_i$ be an arbitrary irrreducible component
of $E$, and $l_0$ any line in $E_i$ which is not contained in
any other $E_j$.
Then we can take a member $D$ of $|-K_U|$ such that
$$
\cases
(1)& \, \, D \cap E_i = l_0, \text{ and}
\qquad \qquad \qquad \qquad \\
(2)& \, \, D \text{ is smooth along } l_0.
\qquad \qquad \qquad \qquad
\endcases
\tag 9.4.1
$$
In fact, we can choose first of all
$D \in |-K_U|$ which satisfies (1)
by the $g\text{-freeness}$ of $|-K_U|$
and $-K_U|_{E_i} = \Cal O_{\Bbb P^2}(1)$.
Then the same argument as in (5.7.5)
shows that $D$ is necessarily smooth along $l_0$,
namely, $D$ satisfies also (2).

Let $h_D : D \to \widetilde{V}$ be the Stein factorization
of $g|_D : D \to V$. This is a bimeromorphic morphism
and $h_D(l_0)$ is a point. Let $\text{Exc}(h_D)$ be
its exceptional locus. Each irreducible component
of $\text{Exc}(h_D)$ has dimension either 1 or 2.
\enddefinition

\proclaim{Lemma 9.5}\ \
Under the Notation 9.4, if we choose
a sufficiently general $D \in |-K_U|$ which
satisfies (9.4.1), then
$g|_{D - (D \cap E)} : D \! - \! (D \cap E) \to Y \! - \! P$
is a finite morphism. In particular,
$l_0$ forms a whole
irreducible component of $\text{Exc }(h_D)$.
\endproclaim

\demo{Proof}\ \
Since each fiber of
$g|_{U \! - \! E} : U \! - \! E \to V \! - \! P$
is isomorphic to the reduced $\Bbb P^1$ (Assumption 9.2),
this is an immediate consequence of Lemma 4.3.
\quad \qed
\enddemo

\proclaim{Proposition 9.6}\ \ (Kawamata [Ka4] (2.4))

Under the Assumption 9.2, Lemma 9.3 and Notation 9.4,
the normal bundle $N_{E_i/X}$ is isomorphic to
$\Cal O_{\Bbb P^2}(-1)^{\oplus 2}$ for each $E_i$.
\endproclaim

\demo{Proof}\ \
Let $l_0 \subset E_i$ be any line
with $l_0 \not\subset \bigcup\limits_{j \not= i} E_j$.
Then by Lemma 9.5,
there is a member $D_0 \in |-K_U|$
which satisfies (9.4.1) such that
$g|_{D_0 - D_0 \cap E} : D_0 \! - \! D_0 \cap E \to Y \! - \! P$
is a finite morphism.
By deforming $D_0$ inside $|-K_U|$,
we obtain a one-parameter family
$$
\Cal D \to \varDelta
$$
(where $\varDelta$ is the unit disc)
of members $D_t :=(\Cal D)_t$ of $|-K_U|$
which satisfies the following:
$$
\cases
l_t := D_t \cap E_i \text{ is a line in } E_i, \\
D_t \text{ is smooth along } l_t, \text{ and} \\
g|_{D_t \! - \!(D_t \cap E)} :
D_t \! - \! (D_t \cap E) \to Y \! - \! P
\text{ is a finite morphism.}
\endcases
\tag 9.6.0
$$
Let
$$
l_t':=D_t \cap (\bigcup\limits_{j \not= i}E_j).
$$
Then $g|_U$ induces a surjective morphism $\psi$ defined as:

\vskip 0.9cm
\centerline{\quad \; $\cup$ \hskip 1.7cm \; $\cup$}
\vskip -1.2cm
$$
\CD
{\psi \;:} @.  \Cal D @>>> {V\times \varDelta} \\
{g|_{D_t} :\;} @. D_t @>>> {V \times \{t\}}
\endCD
$$

\flushpar
We take the Stein factorization
$\psi ': \Cal D \to (V \times \varDelta)\sptilde$
of $\psi$. Then

\flushpar
(9.6.1) \, \, \, $l_t$ is an irreducible component
of the exceptional locus of $(\psi')_t : D_t \to
(V \times \Delta)_t\sptilde$
for each $t \in \varDelta$

\flushpar
by (9.6.0).
A general $\psi '$-very ample divisor $H$ on $\Cal D$
is decomposed as
$$
H=H_1+H_2, \quad
H_1 \cap l_t=\emptyset, \quad
H_2 \cap l_t'=\emptyset \quad
(\forall t \in \varDelta),
$$
by shrinking $(V \times \varDelta)\sptilde$
if necessary.
Then $|H_1|$ gives
a projective surjective morphism
$$
\varphi : \Cal D \to \Cal D'
$$
which is,
on each fiber over $\varDelta$,
a flopping contraction of the smooth 3-fold
$D_t$ whose exceptional locus
is just $l_t \simeq \Bbb P^1$ (9.6.0).
Thus there are only
three possibilities:
$$
N_{l_t/D_t}\simeq \Cal O_{\Bbb P^1}(-1)^{\oplus 2},
\Cal O_{\Bbb P^1} \oplus \Cal O_{\Bbb P^1}(-2)\, \text{ or }\,
\Cal O_{\Bbb P^1}(1) \oplus \Cal O_{\Bbb P^1}(-3)
\tag 9.6.2
$$
Reid [R2], whereas the left-hand side is naturally isomorphic to
$N_{E_i/X} \otimes \Cal O_{l_t}$.

{}From now on,
we shall show
that the type of $N_{l_t/D_t}$ in (9.6.2)
is locally constant on $t$.

Let $L$ be an arbitrary divisor on $\Cal D$
such that $(L\, . \,\,l_t)<0$
for all $t \in \varDelta$,
and let $L':=\varphi(L)$.
\enddemo
\proclaim{Claim 1}\ \
dim$\,R^1\varphi_{t*}\Cal O_{\Cal D_t}(L)$
is constant on $t \in \varDelta$,
by shrinking $\varDelta$ if necessary.
\endproclaim

\demo{Proof}\ \
By ([P] Theorem 3),
we have another projective bimeromorphic morphism
(so-called the {\it simultaneous flop \/} of $\varphi$)
$$
\varphi^+ : \Cal D^+ \to \Cal D'
$$
such that the proper transform
$L^+$ of $L$ is $\varphi^+$-ample.
Then by [KaMaMa],
$R^1\varphi_*^+$ $\Cal O_{\Cal D^+}(L^+-D_t^+)=0$,
since $D_t^+ \underset{\varphi^+}\to {\equiv}0$.
Thus by the exact sequsnce;
$$
0 \longrightarrow \Cal O_{\Cal D^+}(L^+-D_t^+) \longrightarrow
\Cal O_{\Cal D^+}(L^+) \longrightarrow \Cal O_{D_t^+}(L^+)
\longrightarrow 0
$$
we have a surjection
$$
\varphi_*^+ \Cal O_{\Cal D^+}(L^+)
\twoheadrightarrow \varphi_{t*}^+ \Cal O_{D_t^+}(L^+).
$$
Let $D_t':=(\Cal D')_t$,
then this surjection is nothing but
$$
\Cal O_{\Cal D'}(L') \twoheadrightarrow \Cal O_{D_t'}(L').
\tag 9.6.3
$$
On the other hand,
by
$$
0 \longrightarrow \Cal O_{\Cal D}(L-D_t)
\longrightarrow \Cal O_{\Cal D}(L)
\longrightarrow \Cal O_{D_t}(L) \longrightarrow 0,
$$
we have the induced exact sequence:
$$
\varphi_*\Cal O_{\Cal D}(L) \to
\varphi_{t*}\Cal O_{D_t}(L) \to
R^1\varphi_*\Cal O_{\Cal D}(L-D_t)
\to R^1\varphi_*\Cal O_{\Cal D}(L) \to R^1
\varphi_{t*}\Cal O_{D_t}(L) \to 0
$$
Here the first arrow is equal to
$$
\Cal O_{\Cal D'}(L') \longrightarrow \Cal O_{D_t'}(L')
$$
and this is surjective (9.6.3).
Thus the following is exact:
$$
0 \longrightarrow R^1\varphi_*\Cal O_{\Cal D}(L-D_t) \longrightarrow
R^1\varphi_*\Cal O_{\Cal D}(L) \longrightarrow
R^1\varphi_{t*}\Cal O_{D_t}(L) \longrightarrow 0.
$$
Hence
$\dim$ $R^1\varphi_{t *} \Cal O_{D_t}(L)$
is independent on $t$.
\qed
\enddemo

\demo{Proof of Proposition 9.6 Continued}\ \
Let $I_t$ be the defining ideal of $l_t$
(with reduced structure)
in $D_t$.
By the formal function theorem,
$$
\bigl( R^1\varphi_{t*} \Cal O_{D_t}(L) \bigr)
_{\varphi_t(l_t)}^{\wedge}
=\varprojlim H^1 \bigl(l_t\,,\,\,
\Cal O_{D_t}/I_t^m \otimes \Cal O_{D_t}(L) \bigr).
$$
Here the inverse limit above is taken with respect to
the sequence of homomorphisms $\{ \rho_m \}_{m \geq 2}$
given by:
$$
\cases
H^1 \bigl(I_t^{m-1}/I_t^m \otimes \Cal O_{D_t}(L) \bigr)
\longrightarrow
H^1 \bigl( \Cal O_{D_t}/I_t^m \otimes \Cal O_{D_t}(L) \bigr) \\
\qquad \qquad \qquad \qquad \qquad
\overset{\rho_m} \to {\longrightarrow}
H^1 \bigl(\Cal O_{D_t}/I_t^{m-1} \otimes \Cal O_{D_t}(L) \bigr)
\longrightarrow 0
\quad (m \geq 3) \\
\text{and} \\
0 \longrightarrow H^1 \bigl(I_t/I_t^2 \otimes \Cal O_{D_t}(L) \bigr)
\longrightarrow
H^1 \bigl( \Cal O_{D_t}/I_t^2 \otimes \Cal O_{D_t}(L) \bigr) \\
\qquad \qquad \qquad \qquad \qquad \qquad
\qquad \qquad \qquad \qquad \, \, \,
\overset{\rho_2} \to {\longrightarrow}
H^1 \bigl( \Cal O_{l_t}(L) \bigr) \longrightarrow 0
\endcases
\tag 9.6.4
$$
Note that
each $\rho_m$ $(m \geq 2)$ is surjective.
\enddemo

\proclaim{Claim 2}\ \
Let $L$ be chosen so that
$(L\, . \,\,l_t)=-1$.
Then, $R^1\varphi_{t*}\Cal O_{D_t}(L) \not= 0$
if and only if
$N_{l_t/D_t}=\Cal O_{\Bbb P^1}(1) \oplus \Cal O_{\Bbb P^1}(-3)$.
(See (9.6.2).)
\endproclaim

\demo{Proof}\ \
If $N_{l_t/D_t} \simeq
\Cal O_{\Bbb P^1}(1) \oplus \Cal O_{\Bbb P^1}(-3)$
i.e.
$I_t/I_t^2 =\Cal O_{\Bbb P^1}(-1) \oplus \Cal O_{\Bbb P^1}(3)$,
then
$$
H^1\bigl( I_t/I_t^2 \otimes \Cal O_{D_t}(L)\bigr)
=H^1 \bigl(\Cal O_{\Bbb P^1}(-2) \oplus \Cal O_{\Bbb P^1}(2) \bigr)
\not= 0
$$
and $\rho_2$ in (9.6.4)
is not an isomorphism.
Thus $R^1\varphi_{t*}\Cal O_{D_t}(L) \not= 0$.

If $N_{l_t/D_t} \simeq \Cal O_{\Bbb P^1}(-1)^{\oplus 2}$
or $\Cal O_{\Bbb P^1}\oplus \Cal O_{\Bbb P^1}(-2)$
i.e.
$I_t/I_t^2 \simeq \Cal O_{\Bbb P^1}(1)^{\oplus 2}$
or $\Cal O_{\Bbb P^1} \oplus \Cal O_{\Bbb P^1}(2)$ respectively,
then
$$
H^1\bigl(\Cal O_{l_t}(L) \bigr) =
H^1\bigl( \Cal O_{\Bbb P^1}(-1) \bigr)=0
$$
and
$$
H^1 \bigl(I_t^{m-1}/I_t^m \otimes \Cal O_{D_t}(L) \bigr)
=H^1 \bigl( \text{S}^{m-1}(I_t/I_t^2) \otimes
\Cal O_{\Bbb P^1}(-1) \bigr) =0.
$$
Thus
$$
H^1 \bigl(\Cal O_{D_t}/I_t^m \otimes \Cal O_{l_t}(L) \bigr) =0
$$
for all $m$
by (9.6.4),
and hence
$$
R^1\varphi_{t*}\Cal O_{D_t}(L) =0.
\qed
$$
\enddemo

\proclaim{Claim 3}\ \
Let $L$ be chosen so that
$(L \, . \,\, l_t)=-2$.
Then

(i)\quad \,\, dim$\, R^1\varphi_{t*}\Cal O_{D_t}(L)=1$
\quad if \quad $N_{l_t/D_t} \simeq \Cal O_{\Bbb P^1}(-1)^{\oplus 2}$.

(ii)\quad dim$\, R^1\varphi_{t*}\Cal O_{D_t}(L) \geq 2$
\quad if \quad $N_{l_t/D_t} \simeq
\Cal O_{\Bbb P^1} \oplus \Cal O_{\Bbb P^1}(-2)$.
\endproclaim

\demo{Proof}\ \
First we note
$\dim \, H^1 \bigl( \Cal O_{l_t}(L) \bigr)
=h^1 \bigl(\Cal O_{\Bbb P^1}(-2) \bigr)=1$.
In particular,
$R^1 \varphi_{t*}$ $\Cal O_{D_t}(L) \not= 0$ by (9.6.4).
If $N_{l_t/D_t} \simeq \Cal O_{\Bbb P^1}(-1)^{\oplus 2}$
i.e. $I_t/I_t^2 \simeq \Cal O_{\Bbb P^1}(1)^{\oplus 2}$,
then
$$
H^1 \bigl( I_t^{m-1}/I_t^m \otimes \Cal O_{D_t}(L) \bigr)
=H^1\bigl(\text{S}^{m-1} (I_t/I_t^2)
\otimes \Cal O_{\Bbb P^1}(-2) \bigr)=0,
$$
and again by (9.6.4),
all $\rho_m$ ($m \geq 2$) is an isomorphism.
Thus (i) holds.

If $N_{L_t/D_t} =\Cal O_{\Bbb P^1} \oplus \Cal O_{\Bbb P^1}(-2)$
i.e. $I_t/I_t^2 \simeq \Cal O_{\Bbb P^1}\oplus \Cal O_{\Bbb P^1}(2)$,
then
$$
\dim \,H^1 \bigl( I_t/I_t^2 \otimes \Cal O_{D_t}(L) \bigr)
=h^1 \bigl(\Cal O_{\Bbb P^1}(-2) \oplus \Cal O_{\Bbb P^1}\bigr) =1,
$$
and
dim$\, H^1 \bigl(\Cal O_{D_t}/I_t^2\otimes \Cal O_{D_t}(L) \bigr) = 2$
in (9.6.4).
Thus (ii) holds.
\qed
\enddemo

\demo{Proof of Proposition 9.6 Continued}\ \
Claim 1, Claim 2, and Claim 3 implies that
the isomorphism type of
$N_{l_t/D_t} \simeq N_{E_i/X}\otimes \Cal O_{l_t}$
is constant on $t \in \varDelta$.
Since $l_0$ was arbitrarily chosen
so that $l_0 \not\subset \bigcup\limits_{j \not= i} E_j$,
we deduce that
there are at most finitely many lines
at which $N_{E_i/X}$ jumps.
We again put
$N_{E_i/X} \otimes \Cal O_l =
\Cal O_{\Bbb P^1}(a) \oplus \Cal O_{\Bbb P^1}(b)$
with $(a,b)=(-1, -1)$,
$(0, -2)$ or $(1, -3)$,
for a general line $l$ of $E_i$.
If $(a, b)=(0, -2)$ or $(1, -3)$,
then by ([OSS] p.205, Theorem 2.1.4)
which is a generalization of
Grauert-M\"ulich's theorem [GM],
$N_{E_i/X}$ in fact splits into line bundles:
$$
N_{E_i/X} \simeq \Cal O_{\Bbb P^2}(a) \oplus \Cal O_{\Bbb P^2}(b).
$$
In particular
$H^0(N_{E_i/X}) \not=0$
and $H^1(N_{E_i/X})=0$,
which implies that
$E_i$ deforms inside $X$,
a contradiction.
Thus $(a, b)$ must be $(-1, -1)$
and, in turn, by Van de Ven's theorem [V],
we have
$N_{E_i/X}\simeq \Cal O_{\Bbb P^2}(-1)^{\oplus 2}$,
and the proof of Proposition 9.6
is now completed.
\quad \qed
\enddemo

\proclaim{9.7}\ \ Proof of Theorem 9.1.
\endproclaim

We assume (9.2.1), and the results
Lemma 9.3, Proposition 9.6,
and will derive a contradiction.
Since $E_i \simeq \Bbb P^2$ and
$N_{E_i/X} \simeq \Cal O_{\Bbb P^2}(-1)^{\oplus 2}$
for each $E_i$ (Lemma 9.3 and Proposition 9.6),

\flushpar
(9.7.1) \ \ There is a birational morphism
$\varphi_i: X \to X_i'$ over $Y$
to some normal 4-fold $X_i'$ such that
$\varphi_i(E_i)$ is a point, say $Q_i$, and
$\varphi_i$ is an isomorphism elsewhere.

Let $\psi_i : X_i' \to Y$ be the induced morphism.
$g|_{X-E_i} = \psi_i|_{X_i'-Q_i}$
by construction.

First we claim that

\flushpar
(9.7.2) \, \, \,
$E$ is reducible.

In fact, if $E$ is irreducible,
then the fiber of $\psi := \psi_i : X' \to Y$
at $P=g(E)$
is the point $Q := Q_i$, while another fiber
of $\psi$ is a curve since
$g|_{X-E} = \psi|_{X'-Q}$, which
contradicts the upper semi-continuity of the fiber dimension.
Thus $E$ has at least two irreducible components.

Next,

\flushpar
(9.7.3) \, \, \,
For any $E_i$, $E_j$ $(i \not= j)$ with
$E_i \cap E_j \not= \emptyset$,
$E_i \cap E_j$ is a finite set of points.

In fact, if
$E_i \cap E_j$ contains a curve, say $C$,
then $\varphi_i(C)=\{Q_i\}$, and
$\varphi_i|_{E_j} : E_j \to \varphi_i(E_j)$
is an isomorphism outside $C$.
This is impossible, since
$E_j \simeq \Bbb P^2$. Hence (9.7.3) holds.

(9.7.2) and (9.7.3) imply
$$
\text{\Rh}\geq 2,
$$
which contradicts the assumption.
Hence the theorem.
\quad \qed

\head \S 10.\ The classification of $E$ in the case
$l_E(R)=1$ and $\text{\Rh}=1$
\endhead

In this section, we shall classify $E$
in the case $l_E(R)=1$ and $\text{\Rh}=1$.
This case is difficult, and we have at present
a partial answer to the classification:

\proclaim{Theorem 10.1}\ \
Assume that

\flushpar
\qquad \qquad \qquad \qquad \qquad \,
$l_E(R)=1$ \, and \, \Rh$=1$.

\flushpar
Then one of the followings holds:

\flushpar
(1) \, \, \, $E$ is irreducible,
$E \simeq \Bbb P^2$, and
$\Cal O_E(-K_X) \simeq \Cal O_{\Bbb P^2}(1)$.

\flushpar
(2) \, \, \, $E=E_1 \cup E_2$,
$E_1 \simeq E_2 \simeq \Bbb P^2$,
and $E_1 \cap E_2$ is a line of both $E_i$.

\flushpar
\qquad \qquad \qquad \qquad \qquad \quad \, \,
$\Cal O_{E_i}(-K_X) \simeq \Cal O_{\Bbb P^2}(1) \, \, \, \, (i=1,2)$.

\flushpar
(3) \, \, \, $E$ is irreducible.
Let $\nu : \widetilde{E} \to E$ be the normalization.
Then $\widetilde{E} \simeq S_m$, and

\flushpar
\qquad \qquad \qquad \qquad \qquad \qquad \quad
$\nu^* \Cal O_E(-K_X) \simeq \Cal O_{S_m}(1)$.

\flushpar
(4) \, \, \, $E = E_1 \cup E_2$.
Let $\nu_i : \widetilde{E}_i \to E_i$ be
the normalization of $E_i$ $(i=1,2)$. Then
$\widetilde{E}_i \simeq S_{m_i}$ for some
$m_i \geq 2$, and

\flushpar
\qquad \qquad \qquad \qquad \qquad \qquad \, \, \,
$\nu_i^* \Cal O_{E_i}(-K_X) \simeq \Cal O_{S_{m_i}}(1)$.

\flushpar
Moreover, let $v_i$
be the vertex of $\widetilde{E}_i$, then
$\nu_1(v_1) = \nu_2(v_2) =: Q$ in $E$.
$E_1 \cap E_2$ is a ruling of both $E_i$.
\endproclaim

First of all we shall prove the following:

\proclaim{Proposition 10.2}\ \
Under the Assumption 9.0, the followings hold:

\flushpar
(1) \, \, \, Let $l$ be any rational curve of length 1 in $E$.
Then each irreducible component $L$ of $L_U(l)$
(Definition 1.2) contains $E$ and
is numerically $g_U\text{-trivial}$:
$$
L \supset E, \, \, \, L \underset{g_U}\to \equiv 0.
$$

\flushpar
(2) \, \, \, Let $l$ and $l'$
be two rational curves of length 1 in $E$ shch that
$l \cup l'$ deforms outside $E$.
Then $l$ and $l'$ deform inside $E$ to each other.

\endproclaim

\demo{Proof}\ \

\flushpar
{\it (1) \/}
Since $L_X(l)$ is purely codimension 1 in $X$ (Theorem 9.1),
$L_U(l)$ is also purely codimension 1 in $U$. Since
each irreducible component
$L$ of $L_U(l)$ is disjoint from
a general fiber of $g_U$ and since $\rho(U/V)=1$,
we have (1).

\flushpar
{\it (2) \/}\ \
By the assumption, there is
an irreducible component
$L$ ({\it resp. \/} $L'$)
of $L_U(l)$ ({\it resp. \/} $L_U(l')$)
such that
$g_U(L) = g_U(L')$.
Then by (1),
$L=L'$. In particular, $l$ and $l'$
deform inside $U$, and hence inside $E$,
to each other.
\quad \qed
\enddemo

Next:

\proclaim{Proposition 10.3}\ \
Under the Assumption 9.0, $E$ never contains
any irreducible component whose normalization
is a geometrically ruled surface.
\endproclaim

\demo{Proof} \ \
Assume that $E$ contains an irreducible component,
say $E_1$, whose normalization $\widetilde{E}_1$
is isomorphic to
a geometrically ruled surface.
Then by Theorem 3.1,

\flushpar
(10.3.0) \qquad \qquad \qquad \qquad \qquad
$\widetilde{E}_1 \simeq \Sigma_1$ \, or \,
$\Bbb P^1 \times \Bbb P^1$,

\flushpar
as in the proof of Lemma 8.4. Let
$\{ \widetilde{M}, \widetilde{l} \}$
be the minimal section and a ruling
({\it resp. \/} two intersecting rulings)
of $\widetilde{E}_1$,
when $\widetilde{E}_1 \simeq \Sigma_1$
({\it resp. \/} $\widetilde{E}_1
\simeq \Bbb P^1 \times \Bbb P^1$), and
$$
M := \nu(\widetilde{M}) \, \text{ and } \, l := \nu(\widetilde{l}).
$$
Then

\flushpar
(10.3.1) \qquad \, \, \, \,
$\nu^* \Cal O_{E_1}(-K_X) \simeq
\Cal O_{\Sigma_1}(\widetilde{M} + 2\widetilde{l})$ \, \, \,
({\it resp. \/}
$\Cal O_{\Bbb P^1 \times \Bbb P^1}(\widetilde{M} + \widetilde{l})$).

\flushpar
Recall that $|-K_X|$ is $g\text{-free}$
(Theorem 4.1). Thus for a general smooth
$D \in |-K_U|$, $C_D := D \cap E_1$ is a
smooth irreducible curve
with $\nu^* C_D \sim
\widetilde{M} + 2\widetilde{l}$
({\it resp. \/} $\widetilde{M} + \widetilde{l}$),
which is ample.
If $E_1$ has a singular locus of dimension 1,
$C_D$ intersects with $\text{Sing }E_1$
for a general $D$, a contradiction.
Hence

\flushpar
(10.3.2) \qquad \qquad \qquad \qquad \qquad \quad
$\# \text{Sing }E_1 < +\infty$.

\flushpar
We shall prove that
\enddemo

\proclaim{claim}\ \
$M$ and $l$ never deform
inside $E$ to each other.
\endproclaim

\demo{Proof}\ \
Assume to the contrary, and let

\flushpar
(10.3.3) \qquad \qquad \qquad \qquad \qquad \quad \, \, \, \,
$E^- := \bigcup\limits_{i \not=1} E_i$.

\flushpar
Since $M$ and $l$ never deform inside $E_1$
to each other (10.3.2),

\flushpar
(10.3.4) \, \, \, There is a suitable
deformation $l'$ of $l$ inside $E_1$
({\it resp. \/} there are suitable deformations
$M'$ of $M$ and $l'$ of $l$ inside $E_1$)
such that $M$ and $l'$
({\it resp. \/} $M'$ and $l'$)
deform inside $E^-$ to each other.

\flushpar
Let us denote $l'$ ({\it resp. \/} $M'$ and $l'$)
again by $l$ ({\it resp. \/} $M$ and $l$), respectively,
by abuse of notations.
Then (10.3.4) is equivalent to say that
there is a family $\{m_t\}_{t \in T}$
of 1-cycles of $U$
with a 1-dimensional connected parameter space
$T$ (maybe reducible) such that

\flushpar
(10.3.5) \quad \, \,
$(-K_U \, . \, m_t)=1$, \, \, \,
$m_t \subset E^-$ \, \, \, $(\forall t \in T)$, \, \, \,
$m_{t_0} =M$, \, \, and \, \,
$m_{t_1} =l$.

\flushpar
By the freeness of $|-K_U|$
(Theorem 4.1),
we can take a general
smooth member $D \in |-K_U|$
so that
$$
\cases
& \! \! \! \! \! \!
(10.3.6) \, \, \,
\dim (D \cap E) = 1, \, H^1(\Cal O_{D \cap E})=0, \\
& \! \! \! \! \! \!
(10.3.7) \, \, \,
Q_0 := D \cap M \, \text{ and } \, Q_1 := D \cap l \,
\text{ are distinct points of } \, D \cap E_1, \, \text{ and} \\
& \! \! \! \! \! \!
(10.3.8) \, \, \,
\# \{t \in T \, | \, m_t \subset D\} < +\infty.
\endcases
$$

\flushpar
Let
$$
C_1 := D \cap E_1, \, \, \text{ and } \, \,
C^T := D \cap \bigcup\limits_{t \in T} m_t.
$$
Obviously

\flushpar
(10.3.9) \qquad \qquad \qquad \qquad \qquad \, \, \,
$Q_0, Q_1 \in C_1 \cap C^T$.

\flushpar
Moreover by (10.3.5) and (10.3.8),

\flushpar
(10.3.10) \, \, \,
$C^T$ is connected,
purely 1-dimensional, and is contained in $E^-$.

\flushpar
(10.3.7), (10.3.9) and (10.3.10) imply that
$D \cap E$ has two connected sub-1-cycles
$C_1$ and $C^T$ which
have no common irreducible components, and both of which
contains $Q_0$ and $Q_1$.
This contradicts
$H^1(\Cal O_{D \cap E}) = 0$ (10.3.6),
as in the proof of
Proposition 8.5 or 8.6.
Hence the claim. \quad \qed
\enddemo

\demo{Proof of Proposition 10.3 continued}

We shall get a contradiction
in the following cases, in view of (10.3.0):

\flushpar
Case (1) \qquad \qquad \qquad \qquad \qquad \qquad \, \,
$\widetilde{E}_1 \simeq \Sigma_1$.

\flushpar
Case (2) \qquad \qquad \qquad \qquad \qquad \quad \,
$\widetilde{E}_1 \simeq \Bbb P^1 \times \Bbb P^1$.

\flushpar
For {\it Case (1) \/} \, : \quad
Take any irreducible component
$L$ of $L_U(M)$ which contains $M$.
Then
by the claim above, together with
the fact that $M$ is rigid inside $E_1$,
we have $L \not\supset l$,
in particular $(L \, . \, l) >0$.
This contradicts $L \underset{g_U}\to \equiv 0$
(Proposition 10.2 (1)).
Hence $\widetilde{E}_1 \simeq \Sigma_1$
is impossible.

\flushpar
For {\it Case (2) \/} \, : \quad
Recall that
$$
\nu^* \Cal O_{E_1}(-K_X) \simeq
\Cal O_{\Bbb P^1 \times \Bbb P^1} (1,1)
$$
(10.3.1). Since any two points of $E_1$
are joined by a limit conic
(Theorem 3.1), and since
$\# \text{Sing }E_1 < +\infty$ (10.3.2),
it follows that

\flushpar
(10.3.13) \, \, \,
For some limit conic $C_0$ which is contained in $E_1$,
$$
\nu^* \Cal O_{E_1}(C_0) \simeq
\Cal O_{\Bbb P^1 \times \Bbb P^1} (1,1).
$$
Again by (10.3.2),
$$
\{ \nu^*C \, | \, C
\text{ is a deformation of } C_0 \text{ in } E_1 \}
= |\Cal O_{\Bbb P^1 \times \Bbb P^1} (1,1)|.
$$
Hence there is a deformation $C$
of a general fiber of $g$
which is contained in $E_1$ such that
$C$ is reducible: $C = M + l$, and
$$
\align
\nu^*M \in &|\Cal O_{\Bbb P^1 \times \Bbb P^1}(1,0)|, \\
\nu^*l \, \, \in &|\Cal O_{\Bbb P^1 \times \Bbb P^1}(0,1)|.
\endalign
$$
Then by Proposition 10.2 (2),
$M$ and $l$ deform inside $E$ to each other, which
contradicts the above claim.
Hence $\widetilde{E}_1 \simeq \Bbb P^1 \times \Bbb P^1$
is also impossible, and
the proof of Proposition 10.3 is
finished.
\quad \qed
\enddemo

By combining Theorem 2.8 and Proposition 10.3,
we immediately get:

\proclaim{Corollary 10.4}\ \ \
Assume

\flushpar
\qquad \qquad \qquad \qquad \qquad \,
$l_E(R)=1$ \, and \, \Rh$=1$.

\flushpar
Then
$\bigl( \widetilde{E}_i \, , \, \nu_i^* \Cal O_{E_i}(-K_X) \bigr)$
is isomorphic either to one of the followings:

\flushpar
(1) \, \, \,
$\bigl( \Bbb P^2 \, , \, \Cal O_{\Bbb P^2}(1) \bigr)$, \, or

\flushpar
(2) \, \, \,
$\bigl( S_m \, , \, \Cal O_{S_m}(1) \bigr)$ \, \, \, $(m \geq 2)$.
\quad \qed
\endproclaim

\proclaim{Corollary 10.5}\ \
Under the assumption \Rh$=1$,
Any two rational curves of length 1 deform inside $E$
to each other.
\endproclaim

\demo{Proof}\ \
By Corollary 10.4,

\flushpar
(10.5.1) \, \, \,
Any two rational curves of length 1
in $E_i$ deform inside $E_i$ to each other.

\flushpar
On the other hand,
let $l$ be any rational curve of length 1 in
any irreducible component $E_i$ of $E$.
Then

\flushpar
(10.5.2) \qquad \qquad \qquad \qquad \qquad \qquad \, \,
$L_U(l) \supset E$

\flushpar
(Proposition 10.2 (1)).
(10.5.1) and (10.5.2) immediately imply that
$l$ deforms to any rational curve
of length 1 in $E$.
\quad \qed
\enddemo

\proclaim{Proposition 10.6}\ \
If $E$ has an irreducible component, say $E_1$,
whose normalization $\widetilde{E}_1$
is isomorphic to $\Bbb P^2$,
then $E$ is either one of the followings:

\flushpar
(1) \, \, \, $E$ is irreducible and
is isomorphic to $\Bbb P^2$, \, or

\flushpar
(2) \, \, \, $E$ has two irreducible components:
$E = E_1 \cup E_2$,
both of which are isomorphic to $\Bbb P^2$,
and $E_1 \cap E_2$ is a line in $E_i$ \, $(i=1,2)$.
\endproclaim

\demo{Proof}\ \
Assume that $E_1$ say, is an irreducible component of
$E$ whose normalization $\widetilde{E}_1$
is isomorphic to
$\Bbb P^2$. Then since $|-K_X|$ is $g\text{-free}$
(Theorem 4.1),
exactly the same argument as in Lemma 9.3 shows that
$E_1 \simeq \Bbb P^2$.
Thus, by Corollary 10.4,

For each irreducible component $E_i$ of $E$,
one of the followings holds:
$$
\cases
(10.6.1) \, \, \, \, E_i \simeq \Bbb P^2, \, \text{ or} \\
(10.6.2) \, \, \, \text{ The normalization }
\widetilde{E}_i \text{ of } E_i
\text{ is isomorphic to } S_m
\, \, \, (m \geq 2).
\endcases
$$

Now assume that $E$ has irreducible components
one of which is of type (10.6.1) and another of type (10.6.2),
to get a contradiction.

Since any two rational curves of length 1 deform
inside $E$ to each other (Corollary 10.5),
we may assume that
there is a couple $\{E_i, E_j\}$ of
irreducible components of type (10.6.1), (10.6.2),
respectively, of $E$ such that

\flushpar
(10.6.3) \, \, \,
$E_i \cap E_j$ contains a rational curve $l_0$
which is a line
of $E_i$ and is
the image by $\nu_j$ of a ruling of $\widetilde{E}_j$.
In particular, $l_0$ contains the image
by $\nu_j$ of the vertex of $\widetilde{E}_j \simeq S_m$.

\flushpar
(10.6.4) \, \, \,
Take a general line $l \not= l_0$ of $E_i$
which passes through $v_j$, but is not contained in
$E_j$.

\flushpar
Since
$|-K_X|$ is $g\text{-free}$ and since
$\Cal O_{E_i}(-K_U) \simeq \Cal O_{\Bbb P^2}(1)$,
we can take $D \in |-K_U|$ such that
$D|_{E_i} = l$.
Such a $D$ is necessarily smooth,
by the same reason as in
(5.7.5) or (9.4.1).
In particular

\flushpar
(10.6.5) \qquad \qquad \qquad \qquad \qquad \quad \,
$H^1(\Cal O_{D \cap E}) = 0$.

\flushpar
On the other hand, since $D \supset l \ni v_j$ and since
$\nu_j^* \Cal O_{E_j}(-K_X) \simeq \Cal O_{S_m}(1)$,
it follows that
$D \cap E_j$ consists of $m$ rational curves
meeting at $v_j$;

\flushpar
(10.6.6) \qquad \qquad \qquad \qquad
$D \cap E_j = \bigcup\limits_{k=1}^m l_k, \quad
\bigcap\limits_{k=1}^m l_k \ni v_j$,

\flushpar
$(m \geq 2)$. Since $D \supset l$ and
$E_j \not\supset l$ (10.6.4),
it follows that
$D \cap E$ contains $m+1$ rational curves meeting at
$v_j$, which contradicts (10.6.5).
Hence we proved that

\flushpar
(10.6.7) \, \, \, Either all irreducible components
of $E$ are of type (a), or all irreducible components
are of type (b), in (10.6.1).

The rest things we have to prove are the followings:

\flushpar
(10.6.8) \, \, \,
If $E$ has an irreducible component of type (a) above,
then $E$ has at most two irreducible components.

\flushpar
(10.6.9) \, \, \,
Furthermore if $E$ is reducible, then the intersection
is a line of each $\Bbb P^2$.

For (10.6.8): \quad
Assume that $E$ has at least three
$\Bbb P^2\text{'s}$,
to get a contradiction.
By Corollary 10.5 again,
$E$ has three irreducible components
$E_1$, $E_2$, and $E_3$ such that
$E_1 \cap E_2$ contains a line $l_{12}$ of both $E_1$ and $E_2$,
and
$E_2 \cap E_3$ contains a line $l_{13}$ of both $E_2$ and $E_3$.
Then by considering a general $D \in |-K_U|$
passing through the point $l_{12} \cap l_{23}$,
we get a contradiction as exactly the same argument as above.
Hence (10.6.8).

For (10.6.9): \quad
By (10.6.8), $E = E_1 \cup E_2$,
$E_1 \simeq E_2 \simeq \Bbb P^2$.
First if $E_1 \cap E_2$
contains a curve other than a line of $E_1$,
then for a general $D \in |-K_U|$, $D \cap E$
has two irreducible curves
intersecting with
at least two distinct points to each other,
a contradiction.
Next if $E_1 \cap E_2$ has a 0-dimensional
connected component, say $\{Q\}$, then consider a
general line $l$ in $E_1$
passing through $Q$.
Then we can take a smooth $D \in |-K_U|$
such that $D |_{E_1} = l$,
as in (5.7.5) or (9.4.1).
Since $E_1 \cap E_2$ contains a line of $E_1$ which does
not pass through $Q$,
it follows that $D \cap E$
has two irreducible components
meeting at two distinct points,
again a contradiction.
Hence (10.6.9).
\quad \qed
\enddemo

\proclaim{Proposition 10.7}\ \
If $E$ has an irreducible component
$E_i$ whose normalization $\widetilde{E}_i$
is isomorphic to $S_m$ $(m \geq 2)$,
then $E$ is either one of the followings:

\flushpar
(1) \, \, \, $E$ is irreducible:
$\widetilde{E} \simeq S_m$, \, or

\flushpar
(2) \, \, \, $E$ has two irreducible components:
$E = E_1 \cup E_2$,
$\widetilde{E}_1 \simeq S_m$, \,
$\widetilde{E}_2 \simeq S_{m'}$,
and $E_1 \cap E_2$ is the image by $\nu_i$
of a ruling of $\widetilde{E}_i$ for $i=1,2$.
The images of the vertices are the same point.
\endproclaim

\demo{Proof}\ \
This proof is inspired by Iskovskikh [Is], Takeuchi [T].
({\it cf.\/} [Sh1], [R1], [MoMu1,2].)
We consider a
relative and 4-dimensional analogue of
the projection $U \dashrightarrow U'$ of Fano-Iskovskikh-Takeuchi
in the following:

If $E$ is irreducible, then there is nothing to prove.

Assume that $E$ is reducible.
By Corollary 10.4 and Proposition 10.6,
each $E_i$ has the normalization
$\widetilde{E}_i \simeq S_{m_i}$ \, $(m_i \geq 2)$.
Let
$l_i$ and $v_i$ be the image by $\nu_i$
of an arbitrary ruling and the vertex of
$\widetilde{E}_i$, respectively.
Since $l_i\text{'s}$ are the only
rational curves of length 1 in $E_i$,
$v_i$ $(i=1, \dots, n)$ must be the same point
of $E$, by Theorem 3.1.
We denote the point $v_1 = \dots = v_n$
simply by $v \in E$.

Let $\Delta^{\circ} \subset V \! - \! P$
be the discriminant divisor [Bea]
of the conic bundle
$g_U|_{U - E} : U \! - \! E \to V \! - \! P$,
and let $\Delta := \Delta^{\circ} \cup P$.
Moreover let $S_v \ni v$ be the connected component
of the intersection of
two general $g_U\text{-very}$ ample
divisors both of which pass through $v$
so that
$$
\cases
S_v \text{ is a smooth surface proper over } V, \\
g_U|_{S_v - v} : S_v \! - \! v \to g_U(S_v) \! - \! P
\text{ is an isomorphism, and} \\
\dim (g_U (S_v) \cap \Delta) \leq 1.
\endcases
\tag 10.7.1
$$
Consider the blow-up $\varphi : \overline{U} \to U$
with center $S_v$.
Let $\overline{E}$, $\overline{E_i}$, and
$\overline{l_i}$ be the proper transform
of $E$, $E_i$, and $l_i$, respectively,
and let $m := \varphi^{-1}(v) \simeq \Bbb P^1$.
Since $\overline{U}$ is also a smooth 4-fold,

\flushpar
(10.7.2) \qquad \qquad \qquad \qquad \qquad \quad \, \, \, \, \,
$\overline{E_i} \supset m$ \, \, \, $(\forall i)$.

\flushpar
Since $\rho(U/V)=1$,
an easy calculation similar to (5.5.5) shows
$$
\cases
(-K_{\overline{U}} \, \, . \, \overline{l_i})=0, \, \, \,
(-K_{\overline{U}} \, \, . \, m)=1, \, \,
\text{ and} \\
\overline{NE}(U/V) = \Bbb R_{\geq 0} [\overline{l_i}]
+ \Bbb R_{\geq 0} [m].
\endcases
\tag 10.7.3
$$
Consider the contraction morphism
$\psi : \overline{U} \to U'$
associated to the half-line
$\Bbb R_{\geq 0} [\overline{l_i}]$.
By (10.7.3),
$$
\cases
-K_{\overline{U}} \underset{\psi}\to {\equiv} 0, \, \text{ and} \\
B := \psi(m) \text{ is not a point.}
\endcases
\tag 10.7.4
$$
Since $l_i\text{'s}$ cover the whole $\overline{E_i}$,
and these are all contracted by $\psi$,
it follows that $\dim \psi(\overline{E_i}) \leq 1$.
{}From this, (10.7.2), and (10.7.4),

\flushpar
(10.7.5) \qquad \qquad \qquad \qquad \qquad \, \, \, \, \,
$\psi(\overline{E_i}) =B$ \, \, \, $(\forall i)$.

On the other hand, since
$g_U|_{U - E} : U \! - \! E \to V \! - \! P$
is a conic bundle and
$\varphi|_{\overline{U} - \overline{E}}$
is the blow-up of a codimension 2 subsection
(10.7.1),

\flushpar
(10.7.6) \, \, \,
Any irreducible curve in $\overline{U} \! - \! \overline{E}$
which is contracted
by $\psi$ must be lying over
a point of $g_U(S_v) \cap \Delta$
via $g_U \circ \varphi$.

\flushpar
Since $\dim (g_U(S_v) \cap \Delta) \leq 1$ (10.7.1),
it follows that
$\dim \text{Exc }(\psi|_{\overline{U} \! - \! \overline{E}}) =2$,
and hence $\dim \text{Exc }\psi =2$. Namely,
$\psi$ is a flopping contraction.
{}From this and (10.7.5),

\flushpar
(10.7.7) \, \, \,
$\overline{E_i}\text{'s}$
are irreducible components of
$\text{Exc }\psi$, with
$\dim \psi(\overline{E_i})=1$.

\flushpar
(10.7.8) \, \, \,
Now we shall prove that $E$ has
at most two irreducible components:
$$
E = E_1 \cup E_2.
$$

Assume that $E$ has $n \geq 3$
irreducible components, to get a contradiction.
Consider the fiber $C$ of $\psi$
over a general point of $B$ (10.7.5).
Then

\flushpar
(10.7.9) \, \, \,
$C$ has $n$ irreducible components
meeting at a point of $m$.

\flushpar
On the other hand,
by (10.7.4) and [KaMaMa],
$R^1 \psi_* \Cal O_{\overline{U}} =0$, and
in particular

\flushpar
(10.7.10) \qquad \qquad \qquad \qquad \qquad \quad \, \,
$H^1(\Cal O_C)=0$.

\flushpar
These contradict with each other. Hence (10.7.8).

Finally, by Corollary 10.5,

\flushpar
(10.7.11) \, \, \, There exists
a rational curve $l$ of length 1
which is contained in $E_1 \cap E_2$.

Assume that $E_1 \cap E_2 \supsetneq l$,
to get a contradiction.

If $E_1 \cap E_2$ contains a curve
other than $l$, then for a general
smooth $D \in |-K_X|$ (Theorem 4.1),
$D \cap E$ has two irreducible components
which intersect
at two distinct points
with each other, a contradiction to
$H^1(\Cal O_{D \cap E})= 0$.

If $E_1 \cap E_2$ contains a point,
say $Q$, as a connected component,
then let us go back to $\psi : \overline{U} \to U'$,
and let $\overline{Q} := \varphi^{-1}(Q)$.
Then the fiber $C^Q$ of $\psi$
containing $\overline{Q}$
has two irreducible components
which meet at both $\overline{Q}$ and
$C^Q \cap m$, which contradicts
$H^1(\Cal O_{C^Q})=0$ as in (10.7.10).

Hence $E_1 \cap E_2$ is exactly
the image by $\nu_i$ of some ruling in
$\widetilde{E}_i \simeq S_{m_i}$
$(i=1,2)$ (10.7.11), and we get (2) of the proposition.
\quad \qed
\enddemo

\proclaim{10.8}\ \ Proof of Theorem 10.1.
\endproclaim

Now the Theorem 10.1 is the combination of
Proposition 10.6 and 10.7.
\quad \qed

\proclaim{10.9}\ \ Conclusion.
\endproclaim

Theorem 0.5 is a combination of Theorem 3.1, 4.1, and 7.1.

Theorem 0.6 comes directly from Theorem 5.1.
(The assumption of Theorem 5.1 is weaker than
that of Theorem 0.6.)

Theorem 0.7 and 0.8 are exactly Theorem 8.1 and 10.1 with 9.1,
respectively.

\head \S 11.\ Examples.
\endhead

\definition{Example 11.1}\ \
(Beltrametti ([Bel] Example 3.6), Mukai, Reid)

\flushpar
Let $(x_0 : x_1 : x_2)$,
$(y_0 : y_1 : y_2 : y_3)$ be the
homogeneous coordinate of $\Bbb P^2$,
$\Bbb P^3$, respectively, and let
$$
X := \{ x_0 y_0 + x_1 y_1 + x_2 y_2 = 0 \} \subset
\Bbb P^2 \times \Bbb P^3.
$$
$X$ is a Fano 4-fold of Picard number 2
({\it cf. \/} [Wil], [Mu]).
Let $g : X \to Y = \Bbb P^3$ be the second projection.
Then $X$ is a smooth projective 4-fold,
and $g$ is the contraction of an extremal ray of $X$.
$g$ has a fiber $E \simeq \Bbb P^2$
at $P=(0:0:0:1)$, while the other
fibers are all lines in $\Bbb P^2$.
$l_E(R)=2$.
\enddefinition

\definition{Example 11.2}\ \
Let
$g : X \to \Bbb P^3$ be
as in Example 11.1.
Let $D$ be a general (2,2)-divisor of $X$,
and $\beta : X' \to X$
the double cover branched over $D$.
Then $g \circ \beta : X' \to Y$
has a fiber $E' \simeq \Bbb P^1 \times \Bbb P^1$.
$l_{E'}(R)=1$ and
$\rho^{\text{an}}
(X' \! \! \supset \! \! E' / \, Y \! \! \ni \! \! P) \, =2$.
\enddefinition

\definition{Example 11.3}\ \
Let
$g : X \to \Bbb P^3$ be
as in Example 11.1.
Let $D$ be a general (2,2)-divisor of $X$
whose intersection
with $E \simeq \Bbb P^2$ is a reducible conic in $E$.
Let $\beta : \widetilde{X} \to X$ be
the double cover branched over $D$.
Then $g \circ \beta : X' \to Y$
has a fiber $E' \simeq S_2$.
$l_{E'}(R)=1$ and
$\rho^{\text{an}}
(X' \! \! \supset \! \! E' / \, Y \! \! \ni \! \! P) \, =1$.
\enddefinition

\definition{Example 11.4}\ \
Let
$g : X \to \Bbb P^3$ be
as in Example 11.1.
Let $D$ be a general (2,2)-divisor of $X$
whose intersection
with $E \simeq \Bbb P^2$ is a double line in $E$.
Let $\beta : \widetilde{X} \to X$ be
the double cover branched over $D$.
Then $g \circ \beta : X' \to Y$
has a reducible fiber $E'=E_1 \cup E_2$,
$E_1 \simeq E_2 \simeq \Bbb P^2$.
$l_{E'}(R)=1$ and
$\rho^{\text{an}}
(X' \! \! \supset \! \! E' / \, Y \! \! \ni \! \! P) \, =1$.
\enddefinition

\definition{Example 11.5}\ \
Let $(x_0 : x_1 : x_2)$ and
$(y_0 : y_1 : y_2 : y_3)$ be
as in Example 11.1. Let
$$
X := \{ x_0^2 y_0 + x_1^2 y_1 + x_2^2 y_2 = 0 \} \subset
\Bbb P^2 \times \Bbb P^3,
$$
and $g : X \to \Bbb P^3$ the second projection.
Then
$g$ has a fiber $E \simeq \Bbb P^2$
at $y=(0:0:0:1)$, while the other
fibers are all conics in $\Bbb P^2$.
$l_E(R)=1$ and \Rh$=1$.
\enddefinition

Shepherd-Barron told us the following:

\definition{Example 11.6}\ \ (Shepherd-Barron)

Let $E_0 \simeq S_3 \subset \Bbb P^4$ be
the cone over a twisted cubic curve
in $\Bbb P^3$.
Let $C \subset E_0$ be
a general quadric hypersurface section in $\Bbb P^4$,
and $\pi : X \to \Bbb P^4$
the blow-up with center $C$.
Let $E \subset X$ be the proper transform of
$E_0$. Then $X$ is a Fano 4-fold of Picard number 2.
Let $g : X \to Y$ be the contraction
of the extremal ray $R$ other than
$\pi$. Then $Y \simeq \Bbb P^3$, and
$g$ has the fiber $E \simeq S_3$.
$l_E(R)=1$, and \Rh$=1$.
\enddefinition

Mukai told us the following:

\definition{Example 11.7}\ \ (Mukai)

Let $E_0 \simeq \Sigma_1 \subset \Bbb P^4$ be
a cubic surface in $\Bbb P^4$.
Let $C \subset E_0$ be
a general quadric hypersurface section,
$\pi : X \to \Bbb P^4$
the blow-up with center $C$, and
$E \subset X$ the proper transform of
$E_0$. Then $X$ is a Fano 4-fold of Picard number 2, and
the other contraction $g : X \to Y \simeq \Bbb P^3$
has the fiber $E \simeq \Sigma_1$.
$l_E(R)=1$, and \Rh$=2$.
\enddefinition

Although we do not know at present any example
of Mukai-Wi\'sniewski type, Wi\'sniewski
told us the following, which
satisfies all of the condition (MW)
except $\rho(X/Y)=1$.
Mukai also constructed a similar example.

\definition{Example 11.8}\ \ (Wi\'sniewski) \, (locally of
Mukai-Wi\'sniewski type)

Let $Y$ be a quadric hypersurface
in $\Bbb P^4$ with an isolated vertex
$P$. Let $h : W \to Y$ be the blow-up
with the center $P$, and
$F := \text{Exc }h \simeq \Bbb P^1 \times \Bbb P^1$.
Then $W$ has a structure of a $\Bbb P^1\text{-bundle}$
$$
W \simeq \Bbb P \bigl( \Cal O_{\Bbb P^1 \times \Bbb P^1}(-1,0)
\oplus \Cal O_{\Bbb P^1 \times \Bbb P^1}(0,-1) \bigr)
\overset p\to\longrightarrow \Bbb P^1 \times \Bbb P^1 =: B,
$$
and $F$ is its section.
Let $\Cal E := p^* \bigl( \Cal O_{\Bbb P^1 \times \Bbb P^1}
\oplus \Cal O_{\Bbb P^1 \times \Bbb P^1}(0,-1) \bigr)$,
and
$$
\phi_1 : Z^{(1)} := \Bbb P(\Cal E) \to W
$$
be the associated $\Bbb P^1\text{-bundle}$.
Then the tautological section
$S \simeq \Bbb P^1 \times \Bbb P^1$
of $\Cal E |_F$ has the normal bundle
$\Cal O_{\Bbb P^1 \times \Bbb P^1}(0,-1)
\oplus \Cal O_{\Bbb P^1 \times \Bbb P^1}(-1,-1)$
in $Z^{(1)}$. Blow-up $Z^{(1)}$ with the center
$S$:
$$
\phi_2 : Z \to Z^{(1)}.
$$
Let $\psi := h \circ \phi_1 \circ \phi_2 : Z \to Y$,
then
$-K_Z$ is $\psi\text{-ample}$ and
$\psi^{-1}(P)$ is an `alternate union'
$D_1 \cup D_2$:
$$
D_1 \simeq \Bbb P^1 \times \Sigma_1, \,
D_2 \simeq \Sigma_1 \times \Bbb P^1,
\, \text{ with } \, D_1 \cap D_2 =: S' \simeq
\Bbb P^1 \times \Bbb P^1.
$$
Let $\{l_1, l_2\}$ be a couple of
intersecting rulings in $S'$.
Then $\Bbb R_{\geq 0}[l_1]$
and $\Bbb R_{\geq 0}[l_2]$ are
extremal rays of $\overline{NE}(Z/Y)$.
Let
$$
\pi_2 := \text{contr }_{\Bbb R_{\geq 0}[l_2]}
: Z \to X^{(1)},
$$
say, and let $l' := \pi_2(l_1)$. Then
$\pi_2$ is a divisorial contraction
which contracts $D_2$ to $\Sigma_1$,
and
$X^{(1)}$ is smooth [A].
The fiber of $X^{(1)} \to Y$ at $P$
is $\Sigma_1 \cup (\Bbb P^2 \times \Bbb P^1)$
intersecting along $l'$ with each other,
which is the negative section
of $\Sigma_1$ and is a section of
$\Bbb P^2 \times \Bbb P^1 \to \Bbb P^1$.
Finally let
$$
\pi_1 := \text{contr }_{\Bbb R_{\geq 0}[l']}
: X^{(1)} \to X.
$$
Then
$\pi_1$ is a divisorial contraction
which contracts $\Bbb P^2 \times \Bbb P^1$
to $\Bbb P^2$,
and $X$ is smooth [loc.cit].
$g : X \to Y$ has the fiber
$E := g^{-1}(P) \simeq \Bbb P^2 \cup \Bbb P^2$
intersecting at a single point,
and is a $\Bbb P^1\text{-bundle}$ elsewhere.
Namely, $g$ is, locally near $P$,
of Mukai-Wi\'sniewski type.
Unfortunately, we can easily check
that $\rho(X/Y)=2$.
$$
\CD
Z @>{\pi_2}>> X^{(1)} @>{\pi_1}>> X \\
@V{\phi_2}VV @. @VVgV \\
Z^{(1)} @>>{\phi_1}> W @>>h> Y
\endCD
$$
$$
\CD
(\Bbb P^1 \! \times \! \Sigma_1)
\cup (\Sigma_1 \! \times \! \Bbb P^1)
@>>> \Sigma_1 \cup (\Bbb P^2 \! \times \! \Bbb P^1)
@>>> \Bbb P^2 \! \cup \! \Bbb P^2 \\
@VVV @. @VVV \\
\Bbb P^1 \! \times \! \Sigma_1 @>>>
\Bbb P^1 \! \times \! \Bbb P^1 @>>> P
\endCD
$$
\enddefinition

Finally we shall give an example of
extremal contraction of smooth 4-folds which
is a Del Pezzo fibration, whose
special fiber is irreducible and irrational,
in connection with Corollary 0.11 ({\it cf.\/} Fujita [F]).

\definition{Example 11.9}\ \
(Inspired by Hidaka-Oguiso,
after Shepherd-Barron's 11.6 above)

Let $E_0 \subset \Bbb P^4$ be
the cone over an
elliptic curve which is of
(2,2)-complete intersection
in $\Bbb P^3$.
Let $C \subset E_0$ be
a general quadric hypersurface section in $\Bbb P^4$,
and $\pi : X \to \Bbb P^4$
the blow-up with center $C$.
Let $E \subset X$ be the proper transform of
$E_0$. Then $X$ is a Fano 4-fold of Picard number 2,
and the other extremal ray defines
$g : X \to \Bbb P^2$.
$g$ is a Del Pezzo fibration of degree 4,
while it has a fiber $E$, which is
an irrational surface.
\enddefinition

\enddocument